\documentclass{article}

\usepackage[affil-it]{authblk}
\usepackage{amsmath}
\usepackage{graphicx} 
\usepackage{bm}
\usepackage{amsthm}
\usepackage{amssymb}
\usepackage{a4wide}
\usepackage{color}
\usepackage[dvipsnames, table, xcdraw]{xcolor}
\usepackage{hyperref}   
\hypersetup{                       
        colorlinks=true,               
        linkcolor=blue,                
        filecolor=magenta,             
        urlcolor=cyan,                 
        citecolor=ForestGreen          
}
\usepackage[utf8]{inputenc}

\newcommand{\CG}{\text{CG}}

\graphicspath{{./Figures/PDFs}}

\newtheorem{problem}{Problem}
\newtheorem{definition}{Definition}
\newtheorem{theorem}{Theorem}
\newtheorem{proposition}{Proposition}
\newtheorem{remark}{Remark}

\newcommand\iprod[3][]{(#2, #3)_{#1}}
\newcommand\norm[2][]{\|#2\|_{#1}}

\newcommand{\Ltwo}{L^2(\Omega)}
\newcommand{\Hone}{H^1(\Omega)}
\newcommand{\Honez}{H^1_0(\Omega)}
\newcommand{\Hdiv}{H(\mathrm{div},\Omega)}
\newcommand{\Hdivz}{H_0(\mathrm{div},\Omega)}

\newcommand{\XSp}{X_\mathrm{p}}
\newcommand{\XSph}{X_{\mathrm{p},h}}
\newcommand{\XSd}{X_\mathrm{d}}
\newcommand{\XSdh}{X_{\mathrm{d},h}}
\newcommand{\YSp}{Y_\mathrm{p}}
\newcommand{\YSd}{Y_\mathrm{d}}
\newcommand{\ZSp}{Z_\mathrm{p}}
\newcommand{\ZSd}{Z_\mathrm{d}}

\newcommand{\Pp}{P_\mathrm{p}}
\newcommand{\Pph}{P_{\mathrm{p},h}}
\newcommand{\Bp}{B_\mathrm{p}}
\newcommand{\Bsp}{B_\mathrm{stab,p}}
\newcommand{\Sp}{S_\mathrm{p}}
\newcommand{\Sph}{S_{\mathrm{p},h}}
\newcommand{\ap}{a_\mathrm{p}}
\newcommand{\bp}{b_\mathrm{p}}

\newcommand{\Pfem}{P}
\newcommand{\Pperp}{P^\perp}

\newcommand{\Pd}{P_\mathrm{d}}
\newcommand{\Pdh}{P_{\mathrm{d},h}}
\newcommand{\Bd}{B_\mathrm{d}}
\newcommand{\Bsd}{B_\mathrm{stab,d}}
\newcommand{\Sd}{S_\mathrm{d}}
\newcommand{\Sdh}{S_{\mathrm{d},h}}
\newcommand{\ad}{a_\mathrm{d}}
\newcommand{\bd}{b_\mathrm{d}}

\newcommand{\kad}{k_{\ad}}
\newcommand{\kbd}{k_{\bd}}

\newcommand{\zero}{\bm{0}}
\newcommand{\normal}{\hat{\bm{n}}}
\DeclareMathOperator{\Image}{Im}
\DeclareMathOperator{\Kernel}{Ker}

\newcommand{\uF}{u}
\newcommand{\eF}{\bm{e}}
\newcommand{\sF}{\bm{s}}
\newcommand{\lF}{\lambda}
\newcommand{\mF}{\bm{\mu}}
\newcommand{\nF}{\nu}

\newcommand{\eFt}{\tilde{\eF}}
\newcommand{\sFt}{\tilde{\sF}}

\newcommand{\uFh}{\uF_h}
\newcommand{\eFh}{\eF_h}
\newcommand{\sFh}{\sF_h}
\newcommand{\lFh}{\lF_h}
\newcommand{\mFh}{\mF_h}

\newcommand{\xF}{\bm{x}}

\newcommand{\xFh}{\xF_h}

\newcommand{\yF}{\bm{y}}

\newcommand{\zF}{\bm{z}}

\newcommand{\tu}{\tau_{\uF}}
\newcommand{\te}{\tau_{\eF}}
\newcommand{\ts}{\tau_{\sF}}
\newcommand{\tl}{\tau_{\lF}}
\newcommand{\tm}{\tau_{\mF}}

\newcommand{\ttu}{\tau_1}
\newcommand{\tte}{\tau_2}
\newcommand{\tts}{\tau_3}
\newcommand{\ttl}{\tau_4}
\newcommand{\ttm}{\tau_5}
\newcommand{\ttdu}{\tau_6}
\newcommand{\ttds}{\tau_7}
\newcommand{\ttdl}{\tau_8}
\newcommand{\ttdm}{\tau_9}

\newcommand{\ku}{k_{\uF}}
\newcommand{\ke}{k_{\eF}}
\newcommand{\ks}{k_{\sF}}
\newcommand{\kl}{k_{\lF}}
\newcommand{\km}{k_{\mF}}

\newcommand{\kYus}{k_1}
\newcommand{\kYlm}{k_2}
\newcommand{\kPu}{k_3}
\newcommand{\kPl}{k_4}

\newcommand{\length}{\ell_\Omega}
\newcommand{\CPF}{C_\mathrm{PF}}
\newcommand{\Cinv}{C_\mathrm{inv}}

\newcommand{\Pt}{P'}
\newcommand{\Ph}{P_h}

\newcommand{\IN}{\text{ in }}

\newcommand{\dimo}{{d_\Omega}}

\newcommand{\measbd}{\mathcal{H}^{d_\Omega-1}(\partial\Omega)}

\title{A variational multiscale approach to\\
PDE-constrained optimization problems arising in\\
Data-Driven Computational Mechanics}

\author[1,2]{Ramon Codina}
\author[3]{Roberto F. Ausas}
\author[3]{Pedro B. Bazon}
\author[4]{Cristian G. Gebhardt\thanks{Corresponding author. \textit{E-mail-address:} \href{mailto:cristian.gebhardt@uib.no}{cristian.gebhardt@uib.no}}}

\affil[1]{Universitat Polit\`{e}cnica de Catalunya, Jordi Girona 1-3, Edifici C1, 08034 Barcelona, Spain}
\affil[2]{Centre Internacional de M\`{e}todes Num\`{e}rics en Enginyeria (CIMNE), Campus Nord UPC, Gran Capità s/n, 08034 Barcelona, Spain}
\affil[3]{Instituto de Ci\^{e}ncias Matem\'{a}ticas e de Computac\~{a}o, Universidade
de S\~{a}o Paulo, Av. Trab. S\~{a}o Carlense 400, 13566-590 S\~{a}o
Carlos-SP, Brazil}
\affil[4]{University of Bergen, Geophysical Institute and Bergen Offshore Wind Centre (BOW), All\'{e}gaten 70, 5007 Bergen, Norway}

\date{December 23, 2025}

\begin{document}

\maketitle

\begin{abstract}
\noindent We consider the primal and dual forms of the optimality conditions for PDE-contrained optimization problems arising in \textit{Data-Driven Computational Mechanics} when specialized to the reaction-diffusion context. Starting with the continuous setting, we establish well-posedness of such concomitant formulations. Then, we propose stable and consistent finite element approximations for these underlying primal and dual problems relying on the \textit{Variational MultiScale} framework. For quasi-uniform finite element partitions, we investigate approximations' general properties and establish well-posedness for two canonical choices of the sub-grid scales, i.e., the \textit{Algebraic Sub-Grid Scale} and \textit{Orthogonal Sub-Grid Scale}. Moreover, for continuous finite element functions, we are able to move back and forth between the discrete primal and dual formulations only by changing the design of the stabilization parameters. To conclude, we stress-test the proposed approximations through a series of progressively sophisticated cases, providing both a comparative and qualitative assessment of their numerical performance.
\vspace{5mm}


\end{abstract}

\begin{center}    
\textit{In memoriam of Prof. Gustavo C. Buscaglia [1964-2025],\\ a dear friend and inspiring colleague.}
\end{center}

\section{Introduction}

\textit{Data-Driven Computational Mechanics} (DDCM) constitutes the most recent paradigm shift in continuum mechanics, with certainly broad applicability. Specifically for the reaction-diffusion problem, DDCM seeks to find the concentration gradient and the associated flux vector field that best match an empirical extension of gradient-flux pairs. This is achieved by formulating an optimization problem where discrete and continuous fields are simultaneously sought to minimize their discrepancy, provided the continuous fields satisfy the physical conservation and geometrical compatibility laws as well as the Second Law of Thermodynamics.

A basic question that arises is the existence of optimal continuous fields when the discrete fields are held constant. This ``optimization subproblem,'' which is a continuous problem, can be analyzed and solved computationally using standard methods from the linear theory of \textit{Partial Differential Equations} (PDEs). Addressing this specific subproblem is the primary objective of this work. To set the stage for our approach, we first review key related works, differentiate our methodology from existing literature and state the novelty of our present contribution.

DDCM is a highly active research area whose origins date back less than ten years. Its core methodologies have predominantly been established within solid and structural mechanics, as exemplified by its application to the static analysis of multi-bar structures \cite{Kirchdoerfer2016}. Similar ideas led to the development of a data-driven magnetostatic solver \cite{Galetzka2020}. This initial approach was later modified to improve robustness in the presence of constitutive data outliers \cite{Kirchdoerfer2017} and subsequently extended to the dynamic analysis of structures featuring configuration-independent mass matrices \cite{Kirchdoerfer2018}. In the interest of full rigor, we note that the question of existence of solutions has also been addressed in the context of data-driven elasticity by several studies \cite{Conti2018,Conti2020,Gebhardt2026}.
Other research explores a hybrid strategy involving the offline reconstruction of a smooth constitutive manifold from empirical data. This concept was first implemented for time-independent analyses \cite{Gebhardt:2020a}, then adapted to solve time-dependent problems \cite{Gebhardt:2020b}. The approach has been further generalized to tackle challenging \textit{Mathematical Programs with Complementarity Constraints} (MPCCs), particularly relevant to contact mechanics simulations \cite{Gebhardt2024a}. Recently in  \cite{Gebhardt2024b}, the structural properties of the DDCM problem discretized by means of the \textit{Finite Element Method} (FEM) is investigated and its global solvability is established. Based on the same principles, DDCM is applied to discrete electrical circuits in \cite{Gebhardt2025}. Remarkably, the present work addresses a significant gap, as there is only a single antecedent of DDCM in the context of the diffusion-reaction problem, developed previously by some of the authors of this contribution \cite{Bazon2025} where a \textit{Galerkin/Least-Squares} (GLS) method for the stabilization of the resulting finite element discretization is proposed.

All previously introduced works that numerically solve the discrete-continuous optimization problem in DDCM rely on computing a solution to a ``optimization subproblem''. This is because those works adopt an \textit{Alternating Direction Method} (ADM). This work constitutes a crucial major improvement towards applying DDCM in the diffusion-reaction problem. Moreover, our analysis considers two or more spatial dimensions in contrast to contributions that focus on problems defined on closed intervals of the real line \cite{Gebhardt:2020a, Gebhardt:2020b, Gebhardt2024a, Gebhardt2024b}. We introduce novel finite element formulations and establish their well-posedness for the chosen discrete spaces. The formulations proposed are consistent, and therefore convergence is basically a consequence of stability, although we shall not pursue a convergence analysis of these formulations. The required stability properties are achieved by directly approximating the underlying saddle-point problem, an approach related to established methodologies such as those found in \cite{Masud2002, Bochev2006, Badia2010, Burman2023}.

Starting with the continuous setting, we present the primal and dual formulations for the PDE-constrained optimization problem and establish their well-posedness. In addition and by leveraging the necessary regularity of the problem's functional setting, we successfully obtain the minimal irreducible forms of these two formulations. Then, we step into the discrete setting by adopting a finite element discretization. To ensure stability for both the discretized primal and dual problems when using continuous arbitrary approximations, we adopt the \textit{Variational MultiScale} (VMS) framework to systematically derive the necessary stabilization terms \cite{Hughes1995,Hughes1998,Codina2017}. In particular, we analyze some approximation's general properties and establish the consistency and well-posedness of the resulting schemes when considering both the \textit{Algebraic Sub-Grid Scale} (ASGS) and \textit{Orthogonal Sub-Grid Scale} (OSGS) methods~\cite{codina2000, codina2002}. We show that one can move from the primal to the dual formulation just by changing appropriately the design of the stabilization terms, as in \cite{Badia2009} for Darcy's problem or in \cite{Badia2014} for the wave equation. We then stress-test the proposed discrete approximations.

The remainder of this article is structured as follows. Section 2 briefly establishes the context by reviewing the DDCM problem and derives through two main simplifications the PDE-constrained optimization problem targeted by the present contribution. Then, it culminates, in the continuous setting, with the formulation and analysis of the associated primal and dual forms. Section 3 formulates stable finite element approximations for the primal and dual forms and establishes their well-posedness. Section 4 details the numerical implementation, realized using the Firedrake finite element platform, and presents a series of progressively sophisticated examples to illustrate the proposed setting's properties. Section 5 provides concluding remarks, discusses limitations and suggests possible avenues for future work.

\section{Theory}

\subsection{Data-Driven Computational Mechanics}

Given a domain $\Omega\subset\mathbb{R}^\dimo$ such that $\dimo\in\{2,3\}$ (a bounded and connected open subset with Lipschitz-continuous boundary $\partial\Omega$), consider a physical system that occupies the closure $\bar{\Omega}$ and that is subject to body sources, i.e., $q:\Omega\to\mathbb{R}$. All necessary conditions required to formulate the conservation and geometric laws of such a physical system up to constitutive dependencies can be described with the help of: \textit{three primal fields}, e.g., a concentration or temperature $\uF:\Omega\to\mathbb{R}$, the concentration or temperature gradient $\eF:\Omega\to\mathbb{R}^\dimo$ and the flux $\sF:\Omega\to\mathbb{R}^\dimo$; and, \textit{two dual fields} $\lF:\Omega\to\mathbb{R}$ and $\mF:\Omega\to\mathbb{R}^\dimo$, where the first one corresponds to the enforcement of the conservation law in itself and the second one to the enforcement of the geometric law relating the primal fields $\uF$ and $\eF$. In this work, we only consider those conditions to formulate the conservation and geometric laws up to constitutive dependencies that can be formulated as
$$
\Phi((\lF,\mF),(\uF,\eF,\sF);q,\Omega) :=
\iprod[]{\lF}{\zeta\uF+\nabla\cdot\sF-q}+
\iprod[]{\mF}{\eF-\nabla\uF} = 0 \quad\mathrm{in}\quad\Omega,
$$
where $\iprod[]{\cdot}{\cdot}$ represents a suitable inner product, the differential operator $\nabla$ is such that the identity $\iprod[]{\nabla\uF}{\sF} = - \iprod[]{\uF}{\nabla\cdot\sF}$ is identically satisfied for any admissible pair $(\uF,\sF)$ under appropriate boundary conditions \cite[Thm. 6.2]{Kurula_Zwart:2012} and $\zeta\in L^\infty(\Omega)$ is a reaction coefficient as considered by \cite{Bazon2025}, which gives origin to the coupling of all optimality conditions. To simplify the exposition, we shall consider $\zeta$ as a non-negative physical constant in what follows. In addition, achieving physical admissibility of the primal fields $(\eF,\sF)$ requires considering the Second Law of Thermodynamics, i.e., $\eF\cdot\sF\leq 0$ almost everywhere in $\Omega$. 

To include the constitutive dependencies, let us consider the set $\mathcal{N}_m = \lbrace n\rbrace_{n=1}^{N_m}$, where $N_m$ stands for the number of independent experimental measurements available and
$$
\mathfrak{D}_m =\left\lbrace\left(\mathfrak{e}_n,\mathfrak{s}_n\right)\right\rbrace_{n\in\mathcal{N}_m}
$$
represents the set of all experimental pairs that lay, up to a certain tolerance, on the underlying constitutive manifold whose local or global representation is unknown.
For $\mathcal{N}_a\subseteq\mathcal{N}_m$, the set that lists all ``active'' elements of $\mathfrak{D}_m$ (which are to be found subject to $\Phi((\lF,\mF),(\uF,\eF,\sF);q,\Omega)$ vanishes identically) such that $|\mathcal{N}_a|\leq|\mathcal{N}_m|$, let $ P_a(\Omega) = \left\lbrace\Omega_n\right\rbrace_{n\in \mathcal{N}_a}$ be a partition of the domain $\Omega$ (which is also to be sought simultaneously) such that auxiliary fields $(\eFt,\sFt)$ can be constructed, from the active constitutive pairs, over the whole domain in the form of a simple function, namely
$$
(\eFt,\sFt)(\omega) = \sum_{n\in\mathcal{N}_a}\chi_{\Omega_n}(\omega)(\mathfrak{e}_n,\mathfrak{s}_n), \quad\forall\omega\in\Omega,
$$
where $\chi_{\Omega_n}:\Omega\to\mathbb{R}$ is the characteristic function of $\Omega_n$ defined as
$$
\chi_{\Omega_n}(\omega) = 
\left\lbrace
\begin{array}{rl}
1, & \omega\in\Omega_n,\\
0, & \omega\notin\Omega_n .
\end{array}
\right.
$$
This construction for $(\eFt,\sFt)$ is a good approximation in the $L^2(\Omega)$-sense.

By choosing proximity measures $\varphi_e$ and $\varphi_s$ that satisfy $\varphi_e(0) = \varphi_s(0) = 0$, a possible realization of the DDCM approach relies on solving a \textit{Discrete-Continuous Nonlinear Optimization Problem} (DCNLP) given by
\begin{multline*}
\inf_{(\uF,\eF,\sF,\mathcal{N}_a,P_a(\Omega),\nF)}\sup_{(\lF,\mF)}\quad\left\lbrace\varphi_e(\eF-\eFt)+\varphi_s(\sF-\sFt)+\Phi((\lF,\mF),(\uF,\eF,\sF);q,\Omega)\right.\\
\left.\:\:\mathrm{s.t.}\:-\eF\cdot\sF\geq 0,\:\:\nF\geq 0,\:\:\nF\eF\cdot\sF = 0\right\rbrace.
\end{multline*}
Standard choices of the proximity measures are those induced by the $L^2(\Omega)$-norm, e.g.,
$$
\varphi_e(\eF-\eFt)+\varphi_s(\sF-\sFt) = \frac12\norm[]{\eF-\eFt}^2+\frac{\kappa}{2}\norm[]{\sF-\sFt}^2,
$$
where $\kappa \in (0,\infty)$ is a constant providing unit consistency, required to be small enough when necessary (cf. Theorem~\ref{thm:thm2} below).

We can start with our argument, only after formally introducing the PDE-constrained optimization problem at given $(\eFt,\sFt)$ and $q$, which is defined as follows:  
\begin{problem}[$P_0$]
Let $(\eFt,\sFt)$ and $q$ be given. Find $(\yF^*,\nu^*,\zF^*) = ((\uF^*,\eF^*,\sF^*),\nu^*,(\lF^*,\mF^*))\in Y \times N \times Z$ such that
    \begin{equation}
        (\yF^*,\nu^*,\zF^*) =\arg\left\lbrace
        \inf_{\yF\in Y, \nu\in N}\sup_{\zF \in Z}
        \mathcal{L}(\yF,\zF)
        \:\:\mathrm{s.t.}\:
        -\eF\cdot\sF\geq0,\:\:\nu\geq 0,\:\: \nu\eF\cdot\sF = 0
        \right\rbrace,
    \end{equation}
where
$$
\mathcal{L}(\yF,\zF) =         \frac12\norm[]{\eF-\eFt}^2+\frac{\kappa}{2}\norm[]{\sF-\sFt}^2+\iprod[]{\lF}{\nabla\cdot\sF+\zeta\uF-q}+\iprod[]{\mF}{\eF-\nabla\uF}
$$
and subject to the specification of appropriate function spaces and boundary conditions.
\end{problem}

\begin{remark}
Our description of $P_0$ will be completed by specifying: \textit{i}) appropriate function spaces $U,E,S$ ($Y= U\times E\times S$), $L,M$ $(Z = L\times M)$, $N$; and, \textit{ii}) appropriate boundary conditions for the fields $\uF$ and/or $\sF$; whichever these conditions are, in $P_0$ it is assumed that they are homogeneous, for simplicity.  
\end{remark}
The weak form of the optimality conditions for $P_0$ reads
\begin{equation}
\left\lbrace
\begin{array}{rcll}
    \iprod[]{\delta\uF}{\zeta\lF+\nabla\cdot\mF} &=& 0                                & \forall \delta\uF\in U,\\
    \iprod[]{\delta\eF}{\eF+\mF-\nu\sF}                 &=& \iprod[]{\delta\eF}{\eFt}       & \forall \delta\eF\in E,\\   
    \iprod[]{\delta\sF}{\kappa\sF-\nabla\lF-\nu\eF}     &=& \iprod[]{\delta\sF}{\kappa\sFt} & \forall \delta\sF\in S,\\
    \iprod[]{\delta\lF}{\zeta\uF+\nabla\cdot\sF} &=& \iprod[]{\delta\lF}{q}          & \forall \delta\lF\in L,\\
    \iprod[]{\delta\mF}{-\nabla\uF+\eF}          &=& 0                            & \forall \delta\mF\in M,\\
    -\eF\cdot\sF &\geq& 0 &\mathrm{a.e.}\IN\Omega,\\
    \nu &\geq& 0 & \mathrm{a.e.}\IN\Omega,\\
    \nu\eF\cdot\sF &=& 0 &\mathrm{a.e.}\IN\Omega,
\end{array}
\right.    
\end{equation}
provided that the chosen function spaces and considered boundary conditions are such that all boundary terms identically vanish, i.e., at least
$$
\int_{\partial \Omega}\delta\uF (\mF\cdot\normal)\: d\measbd = 0,
\:\:
\int_{\partial \Omega}(\delta \sF\cdot\normal)\lF\: d\measbd = 0
\quad\text{and}\quad
\int_{\partial \Omega}(\delta\mF\cdot\normal)\uF\: d\measbd = 0,
$$
where $\measbd$ stands for the Hausdorff measure on the boundary of $\Omega$.

\begin{remark}\label{rem2}
The consistency in units of $P_0$ is achieved by noting that the dimensional relationship between $\zeta$ and $\kappa$ is given by
$$
[\zeta] = [\kappa]^{-1/2}[\length]^{-2},
$$
where $\length$ represents the characteristic length of the problem and $[\cdot]$ stands for physical units. In particular, since we assume $\zeta$ is a constant we may take $\ell_\Omega^4 \zeta^2 \kappa$ as a dimensionless constant, which we will consider $\ell_\Omega^4 \zeta^2 \kappa > 0$. In the case $\zeta = 0$ (or $\zeta \to 0$) our arguments can be easily adapted. Therefore, the cost function conjugated with constraints' enforcement terms will be dimension free.      
\end{remark}

In the remainder of this work, we specialize DDCM to the reaction-diffusion problem. Our analysis proceeds via two sequential levels of simplification to handle the resulting PDE-constrained optimization problem:
\begin{itemize}
    \item We assume the thermodynamically-consistent auxiliary fields $(\eFt,\sFt)$ are prescribed, which simplifies the original DCNLP problem to a MPCC.
    \item We further simplify the analysis by omitting the enforcement of the Second Law of Thermodynamics for the fields $(\eF,\sF)$. This step linearizes the MPCC, yielding a simpler problem in saddle-point form. Nevertheless, the explicit inclusion of this Second Law that we have proposed in problem $P_0$ should serve to design more robust optimization schemes. We shall not explore this fact in this paper.
\end{itemize}
Although the linear saddle-point problem does not intrinsically guarantee thermodynamic consistency for $(\eF,\sF)$, we verify \textit{a posteriori} that any computed numerical solution fulfills the Second Law of Thermodynamics.

Finally, with the simplified setting adopted, we examine under which conditions this underlying linear problem guarantees a unique global minimizer $\yF^*=(\uF^*,\eF^*,\sF^*)$ with unique multiplier $\zF^
*=(\lF^*,\mF^*)$, considering not only the continuous primal and dual formulations, but also their finite element discretizations. Crucially, we introduce and analyze stabilizations via the VMS approach, which offers an alternative to the GLS technique utilized in a previous work \cite{Bazon2025}.
 
\subsection{Primal form}

The primal formulation of the current problem is given then by the following equations:
$$
\left\lbrace
\begin{array}{rcll}
    \zeta\lF+\nabla\cdot\mF &=& 0          & \IN U^*,\:     U = \Honez,\\
    \eF+\mF                 &=& \eFt       & \IN E^* \simeq E = [\Ltwo]^\dimo,\\   
    \kappa\sF-\nabla\lF     &=& \kappa\sFt & \IN S^* \simeq S = [\Ltwo]^\dimo,\\
    \zeta\uF+\nabla\cdot\sF &=& q          & \IN L^*,\:     L = \Honez,\\
    -\nabla\uF+\eF          &=& \zero      & \IN M^* \simeq M = [\Ltwo]^\dimo,
\end{array}
\right.
$$
where $\simeq$ stands for isomorphism between functional spaces. The standard notation is used for the Lebesgue and Sobolev spaces.

By taking $M = E$, $E = \nabla U$, and $S = \nabla L$, we arrive at
\begin{align}
\left\lbrace
\begin{array}{rcll}
    \zeta\lF-\Delta\uF &=& -\nabla\cdot\eFt& \IN U^*,\: U = \Honez,\\
    -\kappa\zeta\uF-\Delta\lF  &=& \kappa \nabla\cdot\sFt-\kappa q & \IN L^*,\: L = \Honez.
\end{array}
\right. \label{eq:primal-u-l}
\end{align}
Provided higher regularity in the function spaces involved and the domain, we have the reduced strong form
$$
\kappa\zeta^2 u+\Delta \Delta \uF = \Delta\nabla\cdot \eFt-\zeta\kappa\nabla\cdot\sFt+\zeta\kappa q \quad\IN \Ltwo
,$$
where
$$
\uF\in H^4_0(\Omega),
\quad
\eFt \in [H^3(\Omega)]^\dimo,
\quad
\sFt \in [H^1(\Omega)]^\dimo
\quad\text{and}\quad
q\in L^2(\Omega), 
$$
which is a well-posed elliptic problem. In matrix-vector representation, problem~\eqref{eq:primal-u-l} takes the form
$$
-
\begin{bmatrix}
    1 & 0\\
    0 & 1
\end{bmatrix}
\begin{pmatrix}
    \Delta \uF \\
    \Delta \lF 
\end{pmatrix}
+
\begin{bmatrix}
    0 & \zeta\\
    -\zeta\kappa & 0
\end{bmatrix}
\begin{pmatrix}
    \uF \\
    \lF 
\end{pmatrix}
=
\begin{pmatrix}
    \nabla\cdot\eFt\\
    \kappa\nabla\cdot\sFt-\kappa q 
\end{pmatrix}.
$$
The eigenvalues of the ``reaction'' matrix are $\pm \,{\hbox{i}}\, \sqrt{\kappa}\zeta$, showing that the problem cannot be uncoupled (for $\zeta\neq 0$) and behaves as a fourth order problem. For the primal formulation, the fields $\sF$ and $\mF$ are unnecessary. However, for the dual formulation, keeping them makes sense, as we shall see. Furthermore, $\eF$ or $\mF$ could be eliminated one in terms of the other, but this is not possible in the dual formulation.

Finally, the governing equations can be expressed in their weak form as
\begin{align*}
\Bp(\xF,\delta \xF)
&:=     \zeta\iprod[]{      \lF}{\delta\uF}-\iprod[]{      \mF}{\nabla\delta\uF}\\
&\quad      +\iprod[]{      \eF}{\delta\eF}+\iprod[]{      \mF}{      \delta\eF}\\
&\quad+\kappa\iprod[]{      \sF}{\delta\sF}-\iprod[]{\nabla\lF}{      \delta\sF}\\
&\quad+ \zeta\iprod[]{      \uF}{\delta\lF}-\iprod[]{      \sF}{\nabla\delta\lF}\\
&\quad      -\iprod[]{\nabla\uF}{\delta\mF}+\iprod[]{      \eF}{      \delta\mF}\\
&= \iprod[]{\eFt}{\delta\eF}+\kappa\iprod[]{\sFt}{\delta\sF}+\iprod[]{q}{\delta\lF} =: L(\delta \xF). 
\end{align*}
Note that the space $E$
acts as pivot space in this primal formulation.

\begin{definition}

Let the space $\YSp$ be defined as
$$\YSp = \left\lbrace\yF = (\uF,\eF,\sF)\in\Honez\times[\Ltwo]^\dimo\times[\Ltwo]^\dimo\right\rbrace$$
and let us equip $\YSp$ with the norm
$$
\norm[Y_\textrm{p}]{\yF} = \sqrt{
\frac{1}{\length^2}\norm[]{\uF}^2+
                   \norm[]{\nabla \uF}^2+
                   \norm[]{\eF}^2+
             \kappa\norm[]{\sF}^2}.
$$
In addition, let the space $\ZSp$ be defined as
$$\ZSp = \left\lbrace\zF = (\lF,\mF)\in\Honez\times[\Ltwo]^\dimo\right\rbrace$$
and let us equip $\ZSp$ with the norm
$$
\norm[\ZSp]{\zF} = \sqrt{
\frac{1}{\kappa\length^2}\norm[]{\lF}^2+
\frac{1}{\kappa}         \norm[]{\nabla\lF}^2+
                         \norm[]{\mF}^2}.
$$
Finally, let us define $\ap: \YSp\times\YSp\rightarrow \mathbb{R}$ and $\bp : \YSp\times \ZSp \rightarrow \mathbb{R}$ as
\begin{align}
\ap(\yF_1,\yF_2) & = \iprod[]{\eF_1}{\eF_2} + \kappa \iprod[]{\sF_1}{\sF_2}  \label{eq:defap}\\
\bp(\yF,\zF) & = -\iprod[]{\nabla\lF}{\sF}+\zeta\iprod[]{\lF}{\uF}+\iprod[]{\mF}{\eF-\nabla\uF}.\label{eq:defbp}
\end{align}
\end{definition}

\begin{problem}[$\Pp$] Find $\xF = (\yF,\zF)\in\XSp = \YSp\times\ZSp$ such that
\begin{equation}
\Bp(\xF,\delta\xF) = L(\delta\xF),\quad\forall\delta\xF = (\delta\yF, \delta \zF)\in\XSp, 
\end{equation} 
or, equivalently,
\begin{alignat}{3}
\ap(\yF,\delta\yF) + \bp(\delta \yF,\zF) & = 
\iprod[]{\eFt}{\delta\eF}+\kappa\iprod[]{\sFt}{\delta\sF}\qquad && \forall \delta \yF\in \YSp,\nonumber\\
\bp(\yF,\delta\zF) & = \iprod[]{q}{\delta\lF} && \forall \delta \zF\in \ZSp.\nonumber
\end{alignat}
\end{problem}

\begin{theorem}
    Problem $\Pp$ is well posed.
\end{theorem}

\begin{proof}
The bilinear and linear forms are continuous in $\XSp$.
From \cite[Thm. 3.1]{Bazon2025}, the necessary requirements on $\ap$ and $\bp$ to have a well-posed mixed saddle-point problem are met. Finally, the solution continuously depends on $(\eFt,\sFt,q)$.
\end{proof}

\subsection{Dual form}

The dual formulation of the current problem is given then by the following equations:
$$
\left\lbrace
\begin{array}{rcll}
    \zeta\lF+\nabla\cdot\mF &=& 0          & \IN U^* \simeq U = \Ltwo,\\
    \eF+\mF                 &=& \eFt       & \IN E^* \simeq E = [\Ltwo]^\dimo,\\   
    \kappa\sF-\nabla\lF     &=& \kappa\sFt & \IN S^*,       S = \Hdivz,\\
    \zeta\uF+\nabla\cdot\sF &=& q          & \IN L^* \simeq L = \Ltwo,\\
    -\nabla\uF+\eF          &=& \zero      & \IN M^*,       M = \Hdivz,
\end{array}
\right.
$$
where $\Hdivz$ is the space of functions in $[L^2(\Omega)]^\dimo$ with divergence in $L^2(\Omega)$ and vanishing normal trace on $\partial\Omega$.

By taking $M = E$, $U = \sqrt{\kappa}\length^2\nabla \cdot S$, and $L = \sqrt{\kappa}\length^2\nabla \cdot M$, we arrive at
\begin{equation}
\left\lbrace
\begin{array}{rcll}
    \zeta\mF-\nabla(\nabla\cdot\sF) &=& \zeta\eFt-\nabla q& \IN S^*,\: S = \Hdivz,\\
    -\zeta\kappa\sF-\nabla(\nabla\cdot\mF)  &=& -\zeta\kappa\sFt & \IN M^*,\: M = \Hdivz.
\end{array}
\right. \label{eq:dual-mu-s}
\end{equation}
Provided higher regularity in the function spaces involved and the domain, we have the reduced strong form
$$
\zeta^2\kappa \sF+\nabla(\nabla\cdot(\nabla(\nabla\cdot\sF))) = -\zeta \nabla(\nabla\cdot\eFt)+\zeta^2\kappa\nabla\cdot\sFt+\nabla(\nabla\cdot(\nabla q)) \quad\IN [\Ltwo]^\dimo,$$
where we may take
$$\sF \in [H^4(\Omega)]^\dimo,
\quad
\eFt \in [H^2(\Omega)]^\dimo,
\quad
\sFt \in [H^1(\Omega)]^\dimo
\quad\text{and}\quad
q\in H^3(\Omega), 
$$
together with the boundary conditions 
$$\nabla\cdot(\nabla(\nabla\cdot\sF))=0
\quad\text{and}\quad(\nabla(\nabla\cdot\sF))\cdot\normal=0$$
almost everywhere on all $\partial\Omega$, which is a well-posed elliptic problem.
In matrix-vector representation, problem \eqref{eq:dual-mu-s} takes the form
$$
-
\begin{bmatrix}
    I & 0\\
    0 & I
\end{bmatrix}
\begin{pmatrix}
\nabla(\nabla\cdot\sF)\\
\nabla(\nabla\cdot\mF)
\end{pmatrix}
+
\begin{bmatrix}
    0 & \zeta I\\
    -\zeta\kappa I & 0
\end{bmatrix}
\begin{pmatrix}
    \sF \\
    \mF 
\end{pmatrix}
=
\begin{pmatrix}
    \zeta\eFt-\nabla q\\
    -\zeta\kappa\sFt 
\end{pmatrix}.
$$
The eigenvalues of the ``reaction'' matrix are $\pm\, \hbox{i}\, \sqrt{\kappa}\zeta$, showing that the problem cannot be uncoupled (for $\zeta\neq 0$) and behaves as a fourth order problem. For the dual formulation, the fields $\uF$ and $\lF$ are unnecessary. However, for the primal formulation, keeping them makes sense as seen above.

Finally, the governing equations can be expressed in their weak form as
\begin{align*}
\Bd(\xF,\delta\xF)
&:=     \zeta\iprod[]{\lF}{           \delta\uF}+\iprod[]{\nabla\cdot\mF}{           \delta\uF}\\
&\quad      +\iprod[]{\eF}{           \delta\eF}+\iprod[]{           \mF}{           \delta\eF}\\
&\quad+\kappa\iprod[]{\sF}{           \delta\sF}+\iprod[]{           \lF}{\nabla\cdot\delta\sF}\\
&\quad+ \zeta\iprod[]{\uF}{           \delta\lF}+\iprod[]{\nabla\cdot\sF}{           \delta\lF}\\
&\quad      +\iprod[]{\uF}{\nabla\cdot\delta\mF}+\iprod[]{           \eF}{           \delta\mF}\\
&= \iprod[]{\eFt}{\delta\eF}+\kappa\iprod[]{\sFt}{\delta\sF}+\iprod[]{q}{\delta\lF} =: L(\delta \xF). 
\end{align*}
In principle, one could consider $E = \Hdivz$, but as we shall see there is no way to control the divergence of $\eF$. Thus, keeping $\eF$ and $\mF$ as unknowns simultaneously makes sense. Again, note that the space $E$ acts as pivot space in this dual formulation.

\begin{definition} 
Let the space $\YSd$ be defined as
$$
\YSd = \left\lbrace\yF = (\uF,\eF,\sF)\in\Ltwo\times[\Ltwo]^\dimo\times\Hdivz\right\rbrace
$$
and let us equip $\YSd$ with the norm
$$
\norm[\YSd]{\yF} = \sqrt{
\frac{1}{\length^2}\norm[]{\uF}^2+
                   \norm[]{\eF}^2+
             \kappa\norm[]{\sF}^2+
    \kappa\length^2\norm[]{\nabla\cdot\sF}^2}.
$$
In addition, let the space $\ZSd$ be defined as
$$
\ZSd = \left\lbrace\zF = (\lF,\mF)\in\Ltwo\times\Hdivz\right\rbrace
$$
and let us equip $\ZSd$ with the norm
$$
\norm[\ZSd]{\zF} = \sqrt{
\frac{1}{\kappa\length^2}\norm[]{\lF}^2+
                         \norm[]{\mF}^2+
\length^2                \norm[]{\nabla\cdot\mF}^2}.
$$
Finally, let us define $\ad: \YSd\times\YSd\rightarrow \mathbb{R}$ and $\bd : \YSd\times \ZSd \rightarrow \mathbb{R}$ as
\begin{align}
\ad(\yF_1,\yF_2) & = \iprod[]{\eF_1}{\eF_2} + \kappa \iprod[]{\sF_1}{\sF_2}  \label{eq:defad}\\
\bd(\yF,\zF) & =
\iprod[]{\lF}{\nabla\cdot\sF+\zeta\uF}+\iprod[]{\mF}{\eF}+\iprod[]{\nabla\cdot\mF}{\uF}.\label{eq:defbd}
\end{align}
\end{definition}

\begin{problem}[$\Pd$] Find $\xF = (\yF,\zF)\in\XSd = \YSd\times\ZSd$ such that
\begin{equation}
\Bd(\xF,\delta\xF) = L(\delta\xF),\quad\forall\delta\xF = (\delta\xF,\delta\yF)\in\XSd, 
\end{equation} 
or, equivalently,
\begin{alignat}{3}
\ad(\yF,\delta\yF) + \bd(\delta \yF,\zF) & = 
\iprod[]{\eFt}{\delta\eF}+\kappa\iprod[]{\sFt}{\delta\sF}\qquad && \forall \delta \yF\in \YSd,\nonumber\\
\bd(\yF,\delta\zF) & = \iprod[]{q}{\delta\lF} && \forall \delta \zF\in \ZSd.\nonumber
\end{alignat}
\end{problem}

\begin{theorem} \label{thm:thm2}
    For $\kappa$ small enough, problem $\Pd$ is well posed.
\end{theorem}

\begin{proof}
The bilinear forms in \eqref{eq:defad}-\eqref{eq:defbd} and the linear form $L$ are continuous. Moreover, from \eqref{eq:defad} we have that
$$\ad(\yF,\yF) = \norm[]{\eF}^2+\kappa\norm[]{\sF}^2.$$

Let us proof the $\inf$-$\sup$ condition for \eqref{eq:defbd}. Let us consider $(\uF_{\mF},\sF_{\lF},\eF_{\mF})$ such that
$$\uF_{\mF} = \length^2\nabla\cdot\mF,\:\:
\nabla\cdot\sF_{\lF} = \frac{1}{\kappa\length^2}\lF\:\: (\textrm{no }\sF_{\lF}\cdot\normal\textrm{ prescribed on }\partial\Omega)
\:\:\mathrm{and}\:\:
\eF_{\mF} = \mF.$$
Since the divergence operator
$\nabla\cdot:\Hdivz\to\Ltwo$ is surjective, i.e., its closed image is precisely the entire target space or equivalently
$\Image(\nabla\cdot) = \left\lbrace\nabla\cdot\sF:\sF\in\Hdivz\right\rbrace= \Ltwo$, the existence of $\sF_{\lF}$ is warranted.
Then, by inserting
$\yF_{\zF} = (\uF_{\mF},\eF_{\mF},\sF_{\lF})$
in $\bd(\yF,\zF)$ we obtain
$$\bd(\yF_{\zF},\zF) =
\frac{1}{\kappa\length^2}\norm[]{\lF}^2
+\norm[]{\mF}^2
+\zeta\iprod[]{\lF}{\uF_{\mF}}
+\length^2\norm[]{\nabla\cdot\mF}^2.
$$
Now, we use Young's inequality as follows
$$
\zeta \iprod[]{\lF}{\uF_{\mF}} = 
\length^2\iprod[]{\zeta\lF}{\nabla\cdot\mF} \geq -\frac{\length^2\zeta^2}{2}\norm[]{\lF}^2-\frac{\length^2}{2}\norm[]{\nabla\cdot\mF}^2
$$
to arrive at
$$\bd(\yF_{\zF},\zF) \geq
\frac{2-\kappa\zeta^2\length^4}{2\kappa\length^2}\norm[]{\lF}^2
+\norm[]{\mF}^2
+\frac{\length^2}{2}\norm[]{\nabla\cdot\mF}^2
$$
and therefore,
$$
\inf_{(\lF,\mF)\in \ZSd}\sup_{(\uF,\eF,\sF)\in \YSd}\frac{\iprod[]{\lF}{\nabla\cdot\sF+\zeta\uF}+\iprod[]{\mF}{\eF}+\iprod[]{\nabla\cdot\mF}{\uF}}{\norm[\ZSd]{(\lF,\mF)}\norm[\YSd]{(\uF,\eF,\sF)}}
\geq \kbd
$$
with 
$$\kbd = \min\left\lbrace \frac{2-\kappa\zeta^2\length^4}{2}, \frac12\right\rbrace > 0,\:\:\forall\kappa\in\left(0,\frac{2}{\zeta^2\length^4}\right).$$

Let us consider now the coercivity of $\ad$ in the kernel of $\bd$. Let 
\begin{align}
\Kernel(\bd(\cdot,\zF)) = \left\lbrace \yF = (\uF,\eF,\sF)\in\YSd : 
\bd(\yF,\zF) = 0 ~ \forall \zF\in \ZSd\right\rbrace.\label{eq:kerB}
\end{align}
Let us pick $\yF_\infty = (\uF_\infty,\eF_\infty,\sF_\infty)$, with $\uF_\infty \in {\cal C}^\infty_0(\Omega)$. Then, $\yF_\infty\in \Kernel(\bd(\cdot,\zF))$ if $\nabla\cdot\sF_\infty = -\zeta\uF_\infty$, $\eF_\infty = \nabla\uF_\infty$, $\uF_\infty = 0$ on $\partial\Omega$. Using Poincaré's inequality, i.e.,
$$ \norm[]{\nabla\uF_\infty}^2 \geq \frac12\norm[]{\nabla\uF_\infty}^2+\frac{1}{2\CPF^2\length^2}\norm[]{\uF_\infty}^2,
$$
with $\CPF$ a dimensionless constant, we arrive at
$$\ad(\yF_\infty,\yF_\infty) = \norm[]{\eF_\infty}^2+\kappa\norm[]{\sF_\infty}^2 \geq \frac{1}{8\CPF^2\length^2}\norm[]{\uF_\infty}^2+\frac12\norm[]{\eF_\infty}^2+\kappa\norm[]{\sF_\infty}^2+\frac{1}{8\zeta^2\CPF^2\length^2}\norm[]{\nabla\cdot\sF_\infty}^2,$$
or in a more compact form
\begin{equation}
\ad(\yF_\infty,\yF_\infty) \geq \kad\norm[\YSd]{(\uF_\infty,\eF_\infty,\sF_\infty)}^2, \label{eq:coer-ker}
\end{equation}
with
$$\kad = \min\left\lbrace\frac{1}{8\CPF^2},\frac12,\frac{1}{8\CPF^2\kappa\zeta^2\length^4}\right\rbrace > 0, \:\:\forall\frac{1}{\kappa}\in\left(\frac{\zeta^2\length^4}{2},\infty\right),$$
which is consistent with the range of $\kappa$ required to satisfy the previous $\inf$-$\sup$ condition. Recall that we consider $\kappa\zeta^2\length^4 > 0$ (cf. Remark~\ref{rem2}).

The continuity of $\bd$ and a standard compactness argument allows us to establish that \eqref{eq:coer-ker} holds on the whole space in \eqref{eq:kerB}. We can conclude that
$$\kad = \min\left\lbrace\frac12,\frac{1}{16\CPF^2}\right\rbrace > 0.$$
Then, all the necessary requirements on $\ad$ and $\bd$ to have a well-posed mixed saddle-point problem are met.

Finally, the solution continuously depends on $(\eFt,\sFt,q)$.
\end{proof}

\section{Finite element approximation}

For the spatial discretization of both the primal and dual problems, we employ in the following the finite element method. Let us consider then a finite element partition ${\cal P}_h = \{ K\}$ of the domain $\Omega$. Our objective is to consider the simplest possible setting, avoiding technicalities and focusing on the key difficulties of the problem. Thus, we shall only consider {\em quasi-uniform} finite element partitions of diameter $h$. Furthermore, the finite element spaces to be constructed from ${\cal P}_h$ will be considered to be made of ${\cal C}^0(\Omega)$ functions, although the extension to spaces of discontinuous functions is not difficult. All finite element spaces in the following are based on ${\cal P}_h$, although their order can be different. We shall denote the finite element approximation to spaces $U$, $E$, $S$, $L$ and $M$ as $U_h$, $E_h$, $S_h$, $L_h$ and $M_h$, respectively.

As usual, finite element functions and spaces will be indicated with the subscript $h$. This subscript will be used also for the broken $L^2(\Omega)$-inner product and norm, respectively denoted by $\iprod[h]{\cdot}{\cdot}$ and $\norm[h]{\cdot}$.

We will make use of the inverse inequality
\begin{align}
\norm[h]{\nabla v_h} \leq \frac{\Cinv}{h}\norm[h]{v_h},\label{eq:inv-ineq}
\end{align}
where $v_h$ is any finite element function and $\Cinv$ a dimensionless constant.

The problem to be approximated involves five variables. If the standard Galerkin method is used, the associated finite element spaces need to satisfy the global inf-sup condition, which translates into different (`little') inf-sup conditions between the spaces in play. To achieve stable formulations that are independent of the particular finite element spaces chosen, thus avoiding these inf-sup conditions, we are required to include stabilization terms. To this end, we rely on the VMS framework to systematically accomplish this task. The key idea is to split the unknown as $\xF = \xFh+\xF'$, where $\xFh$ is the finite element solution to be computed and $\xF'= (\uF',\eF',\sF',\lF',\mF')$ is the sub-grid scale, which cannot be captured and needs to be approximated in terms of $\xFh$. We shall not describe the rationale to model the sub-grid scales (see for example \cite{Codina2017}). The model we propose is
\begin{equation}
\left\lbrace
\begin{array}{rl}
    \uF' &= \tu\Pt(-\zeta\lFh-\nabla\cdot\mFh),\\
    \eF' &= \te\Pt(\eFt-\eFh-\mFh),\\
    \sF' &= \ts\Pt(\kappa\sFt-\kappa\sFh+\nabla\lFh),\\
    \lF' &= \tl\Pt(q-\zeta\uFh-\nabla\cdot\sFh),\\
    \mF' &= \tm\Pt(\nabla\uFh-\eFh),
    \end{array}
\right.
\label{eq:subscales}
\end{equation}
where $\tu$, $\te$, $\ts$, $\tl$ and $\tm$ are the stabilization parameters defined later and $\Pt$ is a projection operator. Two canonical choices are: \textit{i}) $\Pt$ being the identity in $\Ltwo$, leading to the ASGS method; and, \textit{ii}) $\Pt$ being 
the projection on the orthogonal complement (with respect to the $\Ltwo$-norm) of the finite element space, leading to the OSGS formulation.

Following the VMS approach, we can systematically derive a stabilized bilinear form of the problem associated to the projection $\Pt$, which will be of the form:
$$
B'_\mathrm{stab}(\xFh, \delta \xFh) := 
B(\xFh, \delta \xFh) + S_h'(\xFh, \delta \xFh), 
$$
where $B(\xFh, \delta \xFh)$ is the original bilinear form associated to the continuous problem and $S_h'(\xFh, \delta \xFh)$  is the stabilization that depends, in principle, on the problem, i.e., either primal or dual form, and the particular operator $\Pt$ selected. We shall write it as $S_h' = S_h$ when $\Pt = I$ and as $S_h' = S_h^\bot$ when $\Pt = \Pperp = I - \Ph$, where $\Ph$ is the projection onto the appropriate finite element space. Appropriate subscripts will be used for the primal and the dual problems. To simplify the notation, we will not specify the space onto which $\Pt$ is applied, being it clear by the context.

Likewise, the linear form of the right-hand-side will also be modified to $L_\mathrm{stab}(\delta \xFh)$. This modification has the same expression for the primal and the dual problems, and is given by
\begin{align*}
L'_\mathrm{stab}(\delta \xFh)
& = L(\delta \xFh)
- (\tl\Pt(q), \zeta\delta \uFh + \nabla\cdot\delta\sFh) 
- (\te\Pt(\eFt), \delta \eFh + \delta \mFh) \\
& - (\ts\Pt(\kappa\sFt) , \kappa \delta \sFh - \nabla \delta \lFh).
\end{align*}

\begin{remark} When applied to {\emph{continuous}} finite element functions it is easily checked that
$$\Bp(\xFh, \delta \xFh) = \Bd(\xFh, \delta \xFh).$$
Therefore, the difference between the primal and the dual formulation is only in the expression of $S(\xFh, \delta \xFh)$ and, since the element residuals are identical, in the stabilization parameters. In other words, we will be able to move from the primal to the dual formulations \emph{only changing the design of the stabilization parameters}.
\end{remark}

\subsection{Discretized primal problem} \label{sec:anal_primal}

The stabilized bilinear form for the discrete primal formulation can be written as
$$
\Bsp'(\xFh,\delta\xFh) := 
\Bp(\xFh, \delta\xFh)+\Sph'(\xFh, \delta\xFh), 
$$
where $\Sph'(\xFh, \delta\xFh)$ is, in principle, a mesh-dependent term providing the desired stability. 

In addition, let us consider:
$$
    \tu = \ku h^2, \quad 
    \te = \ke,\quad 
    \ts = \frac{\ks}{\kappa},\quad 
    \tl = -\frac{\kl h^2}{\zeta^2\length^4},\quad 
    \tm = -\km,
$$
where $\ku$, $\ke$, $\ks$, $\kl$ and $\km$ are algorithmic (dimensionless) constants. Then, the stabilization term is as follows:
\begin{align*}
& \Sph'(\xFh,\delta\xFh)\\
&~ :=-\ku h^2\left[\zeta^2\iprod[]{\Pt\lFh}{\delta\lFh}+\zeta\left(\iprod[h]{\Pt\nabla\cdot\mFh}{\delta\lFh}+\iprod[h]{\Pt\lFh}{\nabla\cdot\delta\mFh}\right)+\iprod[h]{\Pt\nabla\cdot\mFh}{\nabla\cdot\delta \mFh}\right]\\
&\quad\:\:-\ke\left[\iprod[]{\Pt\eFh}{\delta\eFh}+\iprod[]{\Pt\mFh}{\delta\eFh}+\iprod[]{\Pt\eFh}{\delta \mFh}+\iprod[]{\Pt\mFh}{\delta\mFh}\right]\\
&\quad\:\:-\frac{\ks}{\kappa}\left[\iprod[]{\Pt\nabla\lFh}{\nabla\delta\lFh}-\kappa\iprod[]{\Pt\sFh}{\nabla\delta\lFh}-\kappa\iprod[]{\Pt\nabla\lFh}{\delta\sFh}+\kappa^2\iprod[]{\Pt\sFh}{\delta\sFh}\right]\\
&\quad\:\:+\frac{\kl h^2}{\zeta^2\length^4}\left[\zeta^2\iprod[]{\Pt\uFh}{\delta\uFh}+\zeta\left(\iprod[h]{\Pt\nabla\cdot\sFh}{\delta\uFh}+\iprod[h]{\Pt\uFh}{\nabla\cdot\delta\sFh}\right)+\iprod[h]{\Pt\nabla\cdot\sFh}{\nabla\cdot\delta \sFh}\right]\\
&\quad\:\:+\km\left[\iprod[]{\Pt\nabla\uFh}{\nabla\delta\uFh}-
\iprod[]{\Pt\eFh}{\nabla\delta\uFh}-\iprod[]{\Pt\nabla\uFh}{\delta\eFh}+\iprod[]{\Pt\eFh}{\delta\eFh}\right]. 
\end{align*}
Note that for the continuous finite element interpolations we are considering, we could in fact replace $\iprod[h]{\cdot}{\cdot}$ by $\iprod[]{\cdot}{\cdot}$.

Particularity, for $\Pt = \Pperp$ we have that
\begin{align*}
\Sph^\perp(\xFh,\delta\xFh)
&=\km \iprod[]{\Pperp(\nabla\uFh)}{\Pperp(\nabla\delta\uFh)}\\
&\:\:\:\:\:
+\frac{\kl h^2}{\zeta^2\length^4}\iprod[h]{\Pperp(\nabla\cdot\sFh)}{\Pperp(\nabla\cdot\delta\sFh)}\\
&\:\:\:\:\:-\frac{\ks}{\kappa}\iprod[]{\Pperp(\nabla\lFh)}{\Pperp(\nabla\delta\lFh)}\\
&\:\:\:\:\:
-\ku h^2\iprod[h]{\Pperp(\nabla\cdot\mFh)}{\Pperp(\nabla\cdot\delta\mFh)}.
\end{align*}

\begin{problem}[$P_{\mathrm{p},h}$] Find $\xFh \in X_h$ such that
\begin{equation}
\Bsp'(\xFh,\delta \xFh) = L'_{\mathrm{stab}}(\delta\xFh),\quad\forall\delta\xFh\in\XSph. 
\end{equation} 
\end{problem}

\begin{definition} Let $\XSph$ be defined as
$$\XSph =  U_h\times E_h\times S_h\times L_h\times M_h$$
and equipped with the norm defined by
\begin{align*}
\norm[\XSph]{\xFh}^2 & =
\frac{1}{\length^2}\norm[]{\uFh}^2+
\norm[]{\nabla\uFh}^2+
\norm[]{\eFh}^2+
\kappa\norm[]{\sFh}^2+
\kappa h^2\norm[h]{\nabla\cdot\sFh}^2\\
& + \frac{1}{\kappa\length^2}\norm[]{\lFh}^2+
\frac{1}{\kappa}\norm[]{\nabla\lFh}^2+
\norm[]{\mFh}^2+
h^2\norm[h]{\nabla\cdot\mFh}^2.
\end{align*}
\end{definition}

\begin{theorem}\label{thm:thm3}
The discretized primal problem $\Pph$ is well posed for $\Pt$ being either the identity (ASGS) or the $L^2(\Omega)$-projection on the orthogonal complement of the finite element space chosen (OSGS). In addition, it consistently approximates the continuous primal problem $\Pp$.   
\end{theorem}

\begin{proof}
The continuity of $\Bsp'(\xFh,\delta\xFh)$ and $L'_{\mathrm{stab}}(\delta\xFh)$ in the norm of $\XSph$ is straightforward. For the inf-sup stability of $\Bsp(\xFh,\delta\xFh)$, 
let $\Pt$ be $I$ and set 
$\delta\xFh^0 = (\uFh,\eFh,\sFh,-\lFh,-\mFh)$.
In addition, let us consider the following inequalities:
\begin{align*}
& -2\km\iprod[]{\eFh}{\nabla\uFh} \geq -\frac12\norm[]{\eFh}^2-4\km^2\norm[]{\nabla\uFh}^2, \quad\text{(Young)}, \\
& \frac{2}{\zeta}\iprod[h]{\nabla\cdot\sFh}{\uFh} \geq -\frac{\kYus}{\zeta^2}\norm[h]{\nabla\cdot \sFh}^2-\frac{1}{\kYus}\norm[]{\uFh}^2,\quad\forall\kYus > 0, \quad\text{(Young)},\\
& 2\zeta\iprod[h]{\nabla\cdot\mFh}{\lFh} \geq -\kYlm\norm[h]{\nabla\cdot\mFh}^2-\frac{\zeta^2}{\kYlm}\norm[]{\lFh}^2,\quad \forall\kYlm > 0,\quad\text{(Young)},\\
& \norm[]{\nabla\uFh}^2 \geq \kPu\norm[]{\nabla\uFh}^2+\frac{(1-\kPu)}{\CPF^2\length^2}\norm[]{\uFh}^2,\quad \forall\kPu\in(0, 1),\quad\text{(Poincaré)},\\
& \norm[]{\nabla\lFh}^2 \geq \kPl\norm[]{\nabla\lFh}^2+\frac{(1-\kPl)}{\CPF^2\length^2}\norm[]{\lFh}^2,\quad \forall\kPl\in(0, 1),\quad\text{(Poincaré)}.
\end{align*}
After some algebraic manipulations, we arrive at
\begin{align*}
\Bsp(\xFh,\delta\xFh^0) \geq 
&\: \ttu\norm[]{\uFh}^2+
    \tte\norm[]{\eFh}^2+
    \tts\norm[]{\sFh}^2+
    \ttl\norm[]{\lFh}^2+
    \ttm\norm[]{\mFh}^2+
    \\    
&\: \ttdu\norm[]{\nabla\uFh}^2+
    \ttds\norm[h]{\nabla\cdot\sFh}^2+ 
    \ttdl\norm[]{\nabla\lFh}^2+
    \ttdm\norm[h]{\nabla\cdot\mFh}^2,
\end{align*}
where
\begin{align*}
    \ttu &= \left(1-\frac{1}{\kYus}\right)\frac{\kl h^2}{\length^4}+\frac{1-\kPu}{\CPF^2\length^2}\km(1-4\km),\\
    \tte &= \km
    -\ke+\frac34,\\
    \tts &= \kappa(1-\ks),\\
    \ttl &= \left(1-\frac{1}{\kYlm}\right)\ku h^2\zeta^2+\frac{1-\kPl}{\CPF^2\length^2}\frac{\ks}{\kappa},\\
    \ttm &= \ke,\\
    \ttdu &= \kPu\km(1-4\km),\\
    \ttds &= (1-\kYus)\frac{\kl h^2}{\length^4\zeta^2},\\
    \ttdl &= \kPl\frac{\ks}{\kappa},\\
    \ttdm &= (1-\kYlm)\ku h^2.\\
\end{align*}
Expression $\Bsp(\xFh,\delta\xFh^0)$ is positive provided that
\begin{align*}
& \ku \geq 0, \ke\in(0,\ke^*),\ks\in(0, 1),\kl\geq 0,\km\in\left(0,\frac14\right),\\
& \kYus\in(\kYus^*,1], \quad \kYlm\in(\kYlm^*,1],\quad \kPu\in(0,1),\quad \kPl\in(0,1),
\end{align*}
where
$$
\ke^* = \km+\frac34,\quad
\kYus^* = \frac{1}{1+\frac{\km(1-4\km)}{\kl}\frac{(1-\kPu)}{\CPF^2}\frac{\length^2}{h^2}},\quad
\kYlm^* = \frac{1}{1+\frac{\ks}{\ku}\frac{(1-\kPl)}{\CPF^2}\frac{1}{\kappa h^2\zeta^2\length^2}}.
$$

For $\Pt = \Pperp$, i.e., the projector onto the $L^2(\Omega)$-orthogonal complement of the chosen particular finite element space, we have
\begin{multline*}
\Bsp^\perp(\xFh,\delta\xFh^0) \geq 
    \norm[]{\eFh}^2+
    \kappa\norm[]{\sFh}^2+\\
    \km\norm[]{\Pperp(\nabla\uFh)}^2+
    \frac{\kl h^2}{\zeta^2\length^4}\norm[h]{\Pperp(\nabla\cdot\sFh)}^2+ 
    \frac{\ks}{\kappa}\norm[]{\Pperp(\nabla\lFh)}^2+
    \ku h^2\norm[h]{\Pperp(\nabla\cdot\mFh)}^2.
\end{multline*}
To obtain full control on $\nabla\uFh$ and $\nabla\lFh$ can be done using the argument briefly summarized next. The task we need to accomplish is to get control over $\norm[]{P(\nabla\uFh)}$ and $\norm[]{P(\nabla\lFh)}$. For such an end, let us take $\delta\xFh^1 = (0,\zero,-\ks\kappa^{-1}\Pfem(\nabla\lFh),\zero,-\km\Pfem(\nabla\uFh))$, which yields 
\begin{align*}
\Bsp^\perp(\xFh,\delta\xFh^1)
&=
+\km\norm[]{\Pfem(\nabla\uFh)}^2 
+\frac{\ks}{\kappa}\norm[]{\Pfem(\nabla\lFh)}^2\\
&\:\:\:\:\:
-\km\iprod[]{\eFh}{\Pfem(\nabla\uFh)}
-\ks\iprod[]{\sFh}{\Pfem(\nabla\lFh)}\\
&\:\:\:\:\:
+\ku\km h^2\iprod[h]{\Pperp(\nabla\cdot\mFh)}{\nabla\cdot(\Pfem(\nabla\uFh))}\\
&\:\:\:\:\:
-\frac{\ks\kl h^2}{\kappa\zeta^2\length^4}\iprod[h]{\Pperp(\nabla\cdot\sFh)}{\nabla\cdot(\Pfem(\nabla\lFh))}.
\end{align*}
By employing Young's inequality for the third and fourth terms and by combining Young's and inverse inequalities for the fifth and sixth terms, namely
$$+\ku\km h^2\iprod[h]{\Pperp(\nabla\cdot\mFh)}{\nabla\cdot(\Pfem(\nabla\uFh))}
\geq
-\frac{k_7\Cinv^2 h^2\ku^2\km}{2}\norm[h]{\Pperp(\nabla\cdot\mFh)}^2
- \frac{\km}{2 k_7}\norm[]{\Pfem(\nabla\uFh)}^2$$
and
$$-\frac{\ks\kl h^2}{\kappa\zeta^2\length^4}\iprod[h]{\Pperp(\nabla\cdot\sFh)}{\nabla\cdot(\Pfem(\nabla\lFh))}
\geq
-\frac{k_8\Cinv^2 h^2\kl^2\ks}{2\zeta^8\length^8\kappa}\norm[h]{\Pperp(\nabla\cdot\sFh)}^2
-\frac{\ks}{2 k_8\kappa}\norm[]{\Pfem(\nabla\lFh)}^2,$$
we arrive at
\begin{align*}
\Bsp^\perp(\xFh,\delta\xFh^1)
&\geq
+\left(1-\frac{1}{2 k_5}-\frac{1}{2 k_7}\right)\km\norm[]{\Pfem(\nabla\uFh)}^2
\\
&\:\:\:\:\:
+\left(1-\frac{1}{2 k_6}-\frac{1}{2 k_8}\right)\frac{\ks}{\kappa}\norm[]{\Pfem(\nabla\lFh)}^2\\
&\:\:\:\:\:
-\frac{k_5\km}{2}\norm[]{\eFh}^2
-\frac{k_6\ks\kappa}{2}\norm[]{\sFh}^2\\
&\:\:\:\:\:
-\frac{k_7\ku^2\km\Cinv^2 h^2}{2}\norm[h]{\Pperp(\nabla\cdot\mFh)}^2\\
&\:\:\:\:\:
-\frac{k_8\ks\kl^2\Cinv^2 h^2}{2 \zeta^4\length^8\kappa}\norm[h]{\Pperp(\nabla\cdot\sFh)}^2\\
&\gtrsim
+\norm[]{\Pfem(\nabla\uFh)}^2
+\frac{1}{\kappa}\norm[]{\Pfem(\nabla\lFh)}^2\\
&\:\:\:\:\:
-\norm[]{\eFh}^2
-\kappa\norm[]{\sFh}^2\\
&\:\:\:\:\:
-\kappa h^2\norm[h]{\Pperp(\nabla\cdot\sFh)}^2
-h^2\norm[h]{\Pperp(\nabla\cdot\mFh)}^2,
\end{align*}
where the control over $\norm[]{P(\nabla\uFh)}$ and $\norm[]{P(\nabla\lFh)}$ is warranted for sufficiently large $k_i$ (for all $i$ in $\{5,6,7,8\}$).
The additional control over $\norm[]{\uFh}$ and $\norm[]{\lFh}$ is achieved by employing Poincaré's inequality.
Finally, setting
$\delta\xFh^\beta = \delta\xFh^0+\beta_1\delta\xFh^1$
and taking $\beta_1$ sufficiently small yields
$$
\Bsp^\perp(\xFh,\delta\xFh^\beta) \gtrsim \norm[\XSph]{\xFh}^2,
$$
which holds true due to conjugation of the estimates for $\Bsp^\perp(\xFh,\delta\xFh^0)$ and $\Bsp^\perp(\xFh,\delta\xFh^1)$. Here and below, $\gtrsim$ stands for $\geq$ up to constants, and likewise for $\lesssim$. Then, the desired $\inf$-$\sup$ condition follows after checking that $\norm[\XSph]{\xFh}\gtrsim\norm[\XSph]{\delta\xFh^\beta}$.
\end{proof}

\begin{proposition}
Let $\ku = \kl = 0$, $\Pt = I$ and $M_h = E_h$. Then, the following residuum orthogonalities hold:
\begin{alignat*}{3}
& \iprod[]{\kappa\sF-\nabla\lFh-\kappa\sFt}{\delta\sFh} = 0\quad && (\text{orthogonality w.r.t. } S_h),\\
& \iprod[]{\eFh+\mFh-\eFt                 }{\delta\eFh} = 0\quad && (\text{orthogonality w.r.t. } E_h),\\
& \iprod[]{-\nabla\uFh+\eFh               }{\delta\mFh} = 0\quad && (\text{orthogonality w.r.t. } M_h).
\end{alignat*}

\end{proposition}

\begin{proof}
First, by setting $\ku = \kl = 0$, it follows
$$\iprod[]{\kappa\sF+\kappa\ts\Pt(\kappa\sFt-\kappa\sFh-\nabla\lFh)+ \nabla\lFh-\kappa\sFt}{\delta\sFh} = 0,$$
$$\iprod[]{\eFh+\te\Pt(\eFt-\eFh-\mFh)+\mFh+\tm\Pt(\nabla\uFh-\eFh)-\eFt}{\delta\eFh} = 0$$
and
$$\iprod[]{-\nabla\uFh+\eFh+\te\Pt(\eFt-\eFh-\mFh)}{\delta\mFh} = 0.$$
Second, by setting $\Pt = I$, we have
$$(1-\ts)\iprod[]{\kappa\sF-\nabla\lFh-\kappa\sFt}{\delta\sFh} = 0,$$
which is proposition's first orthogonality result,%
$$\iprod[]{(\te-1)(\eFt-\eFh-\mFh)+\tm(\nabla\uFh-\eFh)}{\delta\eFh} = 0$$
and
$$\iprod[]{-(\nabla\uFh-\eFh)+\te(\eFt-\eFh-\mFh)}{\delta\mFh} = 0.$$
By choosing
$$\delta\mFh = \tm\delta\eFh$$
and after some algebraic manipulations on the two last equations, we arrive at
$$\iprod[]{\eFh+\mFh-\eFt}{\delta\eFh} = 0,$$
which is proposition's second orthogonality result.
Similarly, by choosing
$$\delta\eFh = \frac{1}{\tm}\delta\mFh$$
and after some algebraic manipulations, we arrive at
$$\iprod[]{-\nabla\uFh+\eFh}{\delta\mFh} = 0,$$
which is proposition's third orthogonality result.%
\end{proof}

\begin{proposition}[minimally stabilized discrete primal formulation]
Let $h$ be small enough with respect to $\length$. For $\ku=\kl=0$ and $M_h = E_h$, the following assertions hold true: \textit{i}) the stabilization term is mesh-independent, i.e., $\Sph(\xFh,\delta\xFh)$ can be written as $\Sp(\xFh,\delta\xFh)$, with
\begin{align*}
\Sp(\xFh,\delta\xFh)
&:=-\ke\left(\iprod[]{\eFh}{\delta\eFh}+\iprod[]{\mFh}{\delta\eFh}+\iprod[]{\eFh}{\delta \mFh}+\iprod[]{\mFh}{\delta\mFh}\right)\\
&\quad\:\:-\frac{\ks}{\kappa}\left(\iprod[]{\nabla\lFh}{\delta\nabla\lFh}-\kappa\iprod[]{\sFh}{\nabla\delta\lFh}-\kappa\iprod[]{\nabla\lFh}{\sFh}+\kappa^2\iprod[]{\sFh}{\delta\sFh}\right)\\
&\quad\:\:+\km\left(\iprod[]{\nabla\uFh}{\nabla\delta\uFh}-
\iprod[]{\eFh}{\nabla\delta\uFh}-\iprod[]{\nabla\uFh}{\delta\eFh}+\iprod[]{\eFh}{\delta\eFh}\right); 
\end{align*}
and, \textit{ii}) the ASGS and the OSGS formulations are equivalent.
\end{proposition}
\begin{proof}
    The first assertion follows intermediately. The second one is a direct consequence of the fact that all subscales
    are (by design) either zero or orthogonal to the corresponding finite element space (due to the previous proposition).
\end{proof}

\subsection{Discretized dual problem} \label{sec:anal_dual}

The stabilized bilinear form for the discrete dual formulation is
$$
\Bsd'(\xFh,\delta\xFh) := 
\Bd(\xFh,\delta\xFh)+\Sd'(\xFh,\delta\xFh), 
$$
where $\Sdh'(\xFh, \delta\xFh)$ is, in principle, a mesh-dependent term providing the desired stability. In addition, let us consider:
$$
    \tu = \ku\length^2,\quad
    \te = \ke,\quad
    \ts = \frac{\ks h^2}{\kappa\length^2},\quad
    \tl = -\frac{\kl}{\zeta^2\length^2},\quad
    \tm = -\frac{\km h^2}{\length^2}.
$$
Then, the stabilization term is as follows:
\begin{align*}
& \Sdh'(\xFh,\delta\xFh)\\
&~ :=-\ku \length^2\left[\zeta^2\iprod[]{\Pt\lFh}{\delta\lFh}+\zeta\left(\iprod[]{\Pt\nabla\cdot\mFh}{\delta\lFh}+\iprod[]{\Pt\lFh}{\nabla\cdot\delta\mFh}\right)+\iprod[]{\Pt\nabla\cdot\mFh}{\nabla\cdot\delta\mFh}\right]\\
&\quad\:\:-\ke\left[\iprod[]{\Pt\eFh}{\delta\eFh}+\iprod[]{\Pt\mFh}{\delta\eFh}+\iprod[]{\Pt\eFh}{\delta \mFh}+\iprod[]{\Pt\mFh}{\delta\mFh}\right]\\
&\quad\:\:-\frac{\ks h^2}{\kappa\length^2}\left[\iprod[h]{\Pt\nabla\lFh}{\delta\nabla\lFh}-\kappa\iprod[h]{\Pt\sFh}{\nabla\delta\lFh}-\kappa\iprod[h]{\Pt\nabla\lFh}{\delta\sFh}+\kappa^2\iprod[]{\Pt\sFh}{\delta\sFh}\right]\\
&\quad\:\:+\frac{\kl}{\zeta^2\length^2}\left[\zeta^2\iprod[]{\Pt\uFh}{\delta\uFh}+\zeta\left(\iprod[]{\Pt\nabla\cdot\sFh}{\delta\uFh}+\iprod[]{\Pt\uFh}{\nabla\cdot\delta\sFh}\right)+\iprod[]{\Pt\nabla\cdot\sFh}{\nabla\cdot\delta\sFh}\right]\\
&\quad\:\:+\frac{\km h^2}{\length^2}\left[\iprod[h]{\Pt\nabla\uFh}{\nabla\delta\uFh}-
\iprod[h]{\Pt\eFh}{\nabla\delta\uFh}-\iprod[h]{\Pt\nabla\uFh}{\delta\eFh}+\iprod[]{\Pt\eFh}{\delta\eFh}\right]. 
\end{align*}
Particularity, for $\Pt = \Pperp$ we have that
\begin{align*}
\Sdh^\perp(\xFh,\delta\xFh)
&= +\frac{\km h^2}{\length^2}\iprod[h]{\Pperp(\nabla\uFh)}{\Pperp(\nabla\delta\uFh)}\\
&\:\:\:\:\:
+\frac{\kl}{\zeta^2\length^2}\iprod[]{\Pperp(\nabla\cdot\sFh)}{\Pperp(\nabla\cdot\delta\sFh)}\\
&\:\:\:\:\:-\frac{\ks h^2}{\kappa\length^2}\iprod[h]{\Pperp(\nabla\lFh)}{\Pperp(\nabla\delta\lFh)}\\
&\:\:\:\:\:
-\ku\length^2\iprod[]{\Pperp(\nabla\cdot\mFh)}{\Pperp(\nabla\cdot\delta\mFh)}.
\end{align*}

\begin{problem}[$\Pdh$] Find $\xFh\in\XSdh$ such that
\begin{equation}
\Bsd'(\xFh,\delta\xFh) = L'_{\mathrm{stab}}(\delta\xFh),\quad\forall\delta\xFh\in\XSdh. 
\end{equation} 
\end{problem}

\begin{definition} Let $\XSdh$ be defined as
$$\XSdh =  U_h \times E_h \times S_h \times L_h \times M_h $$
and equipped with the norm defined by
\begin{align*}
\norm[\XSdh]{\xFh}^2 & =
\frac{1}{\length^2}\norm[]{\uFh}^2+
\frac{h^2}{\length^2}\norm[h]{\nabla\uFh}^2+
\norm[]{\eFh}^2+
\kappa\norm[]{\sFh}^2+
\kappa\length^2\norm[]{\nabla\cdot\sFh}^2\\
& + \frac{1}{\kappa\length^2}\norm[]{\lFh}^2+
\frac{h^2}{\kappa\length^2}\norm[h]{\nabla\lFh}^2+
\norm[]{\mFh}^2+
\length^2\norm[]{\nabla\cdot\mFh}^2.
\end{align*}
\end{definition}

\begin{theorem}\label{thm:thm4}
The discretized dual problem ($P_{\mathrm{d},h}$) is well posed for $\Pt$ being either the identity (ASGS) or the $L^2(\Omega)$-projection on the orthogonal complement of the finite element space chosen (OSGS). In addition, it consistently approximates the continuous problem $P_{\mathrm{d}}$.      
\end{theorem}

\begin{proof}
The continuity of $\Bsd'(\xFh,\delta \xFh)$ and $L'_{\mathrm{stab}}(\delta\xFh)$ in the norm of $\XSdh$ is straightforward. For the inf-sup stability of $\Bsd(\xFh,\delta\xFh)$, let $\Pt$ be equal to $I$ and set 
$\delta \xFh^0 = (\uFh,\eFh,\sFh,-\lFh,-\mFh)$.
In addition, let us consider the following inequalities:
\begin{align*}
& -2 \km \iprod[h]{\eFh}{\nabla\uFh} \geq -\frac12\norm[]{\eFh}^2-4 \km^2\norm[h]{\nabla\uFh}^2, \quad\text{(Young)}, \\
& \frac{2}{\zeta}\iprod[]{\nabla\cdot\sFh}{\uFh} \geq -\frac{\kYus}{\zeta^2}\norm[]{\nabla\cdot\sFh}^2-\frac{1}{\kYus}\norm[]{\uFh}^2,\quad\forall\kYus > 0,\quad\text{(Young)},\\
& 2\zeta\iprod[]{\nabla\cdot\mFh}{\lFh} \geq -\kYlm\norm[]{\nabla\cdot\mFh}^2-\frac{\zeta^2}{\kYlm}\norm[]{\lFh}^2,\quad\forall\kYlm > 0,\quad\text{(Young)},\\
&\norm[h]{\nabla\uFh}^2 \geq \kPu\norm[h]{\nabla\uFh}^2+\frac{(1-\kPu)}{\CPF^2\length^2}\norm[]{\uFh}^2,\quad\forall\kPu\in(0,1),\quad\text{(Poincaré)},\\
& \norm[h]{\nabla\lFh}^2 \geq \kPl\norm[h]{\nabla\lFh}^2+\frac{(1-\kPl)}{\CPF^2\length^2}\norm[]{\lFh}^2,\quad\forall\kPl\in(0,1),\quad\text{(Poincaré)}.
\end{align*}
After some algebraic manipulations, we arrive at
\begin{align*}
\Bsd(\xFh,\delta\xFh^0) \geq 
&\: \ttu\norm[]{\uFh}^2+
    \tte\norm[]{\eFh}^2+
    \tts\norm[]{\sFh}^2+
    \ttl\norm[]{\lFh}^2+
    \ttm\norm[]{\mFh}^2
    \\    
&\: + \ttdu\norm[h]{\nabla\uFh}^2+
    \ttds\norm[]{\nabla\cdot\sFh}^2+ 
    \ttdl\norm[h]{\nabla\lFh}^2+
    \ttdm\norm[]{\nabla\cdot\mFh}^2,
\end{align*}
where
\begin{align*}
    \ttu&= \left(1-\frac{1}{\kYus}\right)\frac{\kl}{\length^2}+\frac{(1-\kPu)}{\CPF^2\length^2}\km(1-4\km)\frac{h^2}{\length^2},\\
    \tte &= \left(\km-\frac14\right)\frac{h^2}{\length^2}
    -\ke+1,\\
    \tts &= \kappa\left(1-\ks\frac{h^2}{\length^2}\right),\\
    \ttl &= \left(1-\frac{1}{\kYlm}\right)\ku\length^2\zeta^2+\frac{(1-\kPl)}{\CPF^2\length^2}\frac{\ks h^2}{\kappa\length^2},\\
    \ttm &= \ke,\\
    \ttdu &= \kPu\km(1-4\km)\frac{h^2}{\length^2},\\
    \ttds &= (1-\kYus)\frac{\kl}{\length^2\zeta^2},\\
    \ttdl &= \kPl\frac{\ks h^2}{\kappa\length^2},\\
    \ttdm &= (1-\kYlm)\ku \length^2.
\end{align*}
Expresssion $\Bsd(\xFh,\delta\xFh^0)$ is positive provided that
\begin{align*}   
& \ku> 0,\ke\in(0,\ke^*),\ks\in\left[0,\ks^*\right),\kl > 0,\km\in\left[0,\frac14\right],\\
& \kYus\in(\kYus^*,1),\quad 
\kYlm\in(\kYlm^*,1),\quad
\kPu,\in(0,1),\quad
\kPl\in(0, 1),
\end{align*}
where
\begin{align*}
& \ke^* = 1+\left(\km-\frac14\right)\frac{h^2}{\length^2},\quad
\:\:
\ks^* = \frac{\length^2}{h^2}, \\
& \kYus^* = \frac{1}{1+\frac{\km (1-4\km)}{\kl}\frac{(1-\kPu)}{\CPF^2}\frac{h^2}{\length^2}}, 
\:\:
\kYlm^* = \frac{1}{1+\frac{\ks}{\ku}\frac{(1-\kPl)}{\CPF^2}\frac{h^2}{\kappa\zeta^2\length^6}}.
\end{align*}
For $\Pt = \Pperp$, i.e., the projector onto the $L^2(\Omega)$-orthogonal complement of the chosen particular finite element space, we have
\begin{multline*}
\Bsd^\perp(\xFh,\delta\xFh^0) \geq 
    \norm[]{\eFh}^2+
    \kappa\norm[]{\sFh}^2\\
    +\frac{\km h^2}{\length^2}\norm[h]{\Pperp(\nabla\uFh)}^2+
    \frac{\kl}{\zeta^2\length^2}\norm[]{\Pperp(\nabla\cdot\sFh)}^2  
    + \frac{\ks h^2}{\kappa\length^2}\norm[h]{\Pperp(\nabla\lFh)}^2+
    \ku\length^2\norm[]{\Pperp(\nabla\cdot\mFh)}^2.
\end{multline*}
Stability on the whole norm of $\nabla\cdot\sFh$ and $\nabla\cdot\mFh$ is obtained using an argument similar to that of Theorem~\ref{thm:thm3}.
The first task we need to accomplish is to get control over $\norm[]{\Pfem(\nabla\cdot\sFh)}$ and $\norm[]{\Pfem(\nabla\cdot\mFh)}$. For such an end, let us take $\delta\xFh^1 = (\ku\length^2\Pfem(\nabla\cdot\mFh),\zero,\zero,\kl\length^{-2}\zeta^{-2}\Pfem(\nabla\cdot\sFh),\zero)$, which yields 
\begin{align*}
\Bsd^\perp(\xFh,\delta\xFh^1)
&=
+\frac{\kl}{\length^2\zeta^2}\norm[]{\Pfem(\nabla\cdot\sFh)}^2
+\ku\length^2\norm[]{\Pfem(\nabla\cdot\mFh)}^2
\\
&\:\:\:\:\:
+\frac{\kl}{\length^2\zeta}\iprod[]{\uFh}{\Pfem(\nabla\cdot\sFh)}
+\ku\zeta\length^2\iprod[]{\lFh}{\Pfem(\nabla\cdot\mFh)}
\\
&\:\:\:\:\:
+\ku\km h^2\iprod[h]{\Pperp(\nabla\uFh)}{\nabla(\Pfem(\nabla\cdot\mFh))}
\\
&\:\:\:\:\:
-\frac{\ks\kl h^2}{\kappa\zeta^2\length^4}\iprod[h]{\Pperp(\nabla\lFh)}{\nabla(\Pfem(\nabla\cdot\sFh))}.
\end{align*}
By employing Young's inequality for the third and fourth terms and by combining Young's and inverse inequalities for the fifth and sixth terms, namely,
$$+\ku\km h^2\iprod[h]
{\Pperp(\nabla\uFh)}{\nabla(\Pfem(\nabla\cdot\mFh))}
\geq
-\frac{k_7\Cinv^2 h^2\ku\km^2}{2\length^2}\norm[h]{\Pperp(\nabla\uFh)}^2
- \frac{\ku\length^2}{2 k_7}\norm[h]{\Pfem(\nabla\cdot\mFh)}^2$$
and
$$-\frac{\ks\kl h^2}{\kappa\zeta^2\length^4}\iprod[h]
{\Pperp(\nabla\lFh)}{\nabla(\Pfem(\nabla\cdot\sFh))}
\geq
-\frac{k_8\Cinv^2 h^2\ks^2\kl}{2\zeta^2\length^6\kappa^2}\norm[h]{\Pperp(\nabla\lFh)}^2
-\frac{\kl}{2 k_8\length^2\zeta^2}\norm[h]{\Pfem(\nabla\cdot\sFh)}^2,$$
we arrive at
\begin{align*}
\Bsd^\perp(\xFh,\delta\xFh^1)
&\geq
+\left(1-\frac{1}{2 k_5}-\frac{1}{2 k_8}\right)\frac{\kl}{\length^2\zeta^2}\norm[]{\Pfem(\nabla\cdot\sFh)}^2\\
&\:\:\:\:\:
+\left(1-\frac{1}{2 k_6}-\frac{1}{2 k_7}\right)\ku\length^2\norm[]{\Pfem(\nabla\cdot\mFh)}^2\\
&\:\:\:\:\:
-\frac{k_5\kl}{2\length^2}\norm[]{\uFh}^2
-\frac{k_7\Cinv^2 h^2\km^2\ku}{2\length^2}\norm[h]{\Pperp(\nabla\uFh)}^2
\\
&\:\:\:\:\:
-\frac{k_6\ku\length^2\zeta^2}{2}\norm[]{\lFh}^2
-\frac{k_8\Cinv^2 h^2\ks^2\kl}{2\kappa^2\length^6\zeta^2}\norm[h]{\Pperp(\nabla\lFh)}^2
\\
&\gtrsim
+\kappa\length^2\norm[]{\Pfem(\nabla\cdot\sFh)}^2
+\length^2\norm[]{\Pfem(\nabla\cdot\mFh)}^2
\\
&\:\:\:\:\:
-\frac{1}{2\length^2}\norm[]{\uFh}^2
-\frac{h^2}{\length^2}\norm[h]{\Pperp(\nabla\uFh)}^2
\\
&\:\:\:\:\:
-\frac{1}{\kappa\length^2}\norm[]{\lFh}^2
-\frac{h^2}{\kappa\length^2}\norm[h]{\Pperp(\nabla\lFh)}^2,
\end{align*}
where the control over $\norm[]{\Pfem(\nabla\cdot\sFh)}$ and $\norm[]{\Pfem(\nabla\cdot\mFh)}$ is warranted for sufficiently large $k_i$ (for all $i$ in $\{5,6,7,8\}$).

The second task we need to accomplish, is to get control over $\norm[]{\Pfem_{L_h}(\uFh)}$ and $\norm[]{\Pfem_{U_h}(\lFh)}$, where $\Pfem_{L_h}$ is the $L^2(\Omega)$-projection onto $L_h$ and $\Pfem_{U_h}$ is the $L^2(\Omega)$-projection onto $U_h$, which now need to be distinguished. The natural choice would be to take $\sqrt{\kappa}U_h = L_h$ (recall that both spaces are made of continuous functions), but for generality we consider the case in which these spaces are different. For such an end, let us take
$\delta \xFh^2 = (\zeta^{-1}\length^{-2}\kappa^{-1}\Pfem_{U_h}(\lFh)
,\zero,\zero,\zeta^{-1}\length^{-2}\Pfem_{L_h}(\uFh)
,\zero)$, which yields
\begin{align*}
\Bsd^\perp(\xFh,\delta\xFh^2)
&=
+\frac{1}{\length^2}\norm[]{\Pfem_{L_h}(\uFh)}^2
+\frac{1}{\kappa\length^2}\norm[]{\Pfem_{U_h}(\lFh)}^2
\\
&\:\:\:\:\:
+\frac{\km h^2}{\kappa\length^4\zeta}\iprod[h]{\Pperp(\nabla\uFh)}{\nabla(\Pfem_{U_h}(\lFh))}
-\frac{\ks h^2}{\kappa\length^4\zeta}\iprod[h]{\nabla(\Pfem_{L_h}(\uFh))}{\Pperp(\nabla\lFh)}\\
&\:\:\:\:\:
+\frac{1}{\length^2\zeta}\iprod[]{\Pfem_{L_h}(\uFh)}{\nabla\cdot\sFh}
+\frac{1}{\kappa\length^2\zeta}\iprod[]{\Pfem_{U_h}(\lFh)}{\nabla\cdot\mFh}.
\end{align*}
By combining Young's and inverse inequalities for the third and fourth terms and employing Young's inequality to the fifth and sixth terms, we arrive at
\begin{align*}
\Bsd^\perp(\xFh,\delta\xFh^2)
&\geq
+\left(1-\frac{1}{2 k_9}-\frac{\ks}{2 k_{11}\kappa\length^4\zeta^2}\right)\frac{1}{\length^2}\norm[]{\Pfem_{L_h}(\uFh)}^2
-\frac{k_{11} \km h^2}{2\length^2}\norm[h]{\Pperp(\nabla\uFh)}^2\\
&\:\:\:\:\:
+\left(1-\frac{1}{2 k_{10}}-\frac{\km}{2k_{12}\kappa\length^4\zeta^2}\right)\frac{1}{\kappa\length^2}\norm[]{\Pfem_{U_h}(\lFh)}^2
-\frac{k_{12}\ks h^2}{2\kappa\length^2}\norm[h]{\Pperp(\nabla\lFh)}^2\\
&\:\:\:\:\:
-\frac{k_9}{2\length^2\zeta^2}\norm[]{\nabla\cdot\sFh}^2
-\frac{k_{10}}{2\kappa\length^2\zeta^2}\norm[]{\nabla\cdot\mFh}^2\\
&\gtrsim
+\frac{1}{\length^2}\norm[]{\Pfem_{L_h}(\uFh)}^2
-\frac{h^2}{\length^2}\norm[h]{\Pperp(\nabla\uFh)}^2\\
&\:\:\:\:\:
+\frac{1}{\kappa\length^2}\norm[]{\Pfem_{U_h}(\lFh)}^2
-\frac{h^2}{\kappa\length^2}\norm[h]{\Pperp(\nabla\lFh)}^2\\
&\:\:\:\:\:
-\kappa\length^2\norm[]{\nabla\cdot\sFh}^2
-\length^2\norm[]{\nabla\cdot\mFh}^2,
\end{align*}
where the control over $\norm[]{\Pfem_{L_h}(\uFh)}$ and $\norm[]{\Pfem_{U_h}(\lFh)}$ is warranted for sufficiently large $k_i$ (for all $i$ in $\{9,10,11,12\}$) and sufficiently small $\ks$ and $\km$.

Finally, setting
$\delta\xFh^\beta = \delta\xFh^0+\beta_1\delta\xFh^1+\beta_2\delta\xFh^2$
and taking $\beta_1$ and $\beta_2$ positive and sufficiently small yields
$$
\Bsd^\perp(\xFh,\delta\xFh^\beta) \gtrsim \norm[\XSdh^\ast]{\xFh}^2,
$$
which holds true due to the conjugation of the estimates for $\Bsd^\perp(\xFh,\delta\xFh^i)$ (for all $i$ in $\{0,1,2\}$), and the norm $\norm[\XSdh^\ast]{\cdot}$ is the same as $\norm[\XSdh]{\cdot}$ replacing the $L^2(\Omega)$-norm of $\uFh$ and $\lFh$ by $P_{L_h}(\uFh)$ and $P_{U_h}(\lFh)$, respectively. Then, the desired $\inf$-$\sup$ condition follows after checking that $\norm[\XSdh^\ast]{\xFh}\gtrsim\norm[\XSdh^\ast]{\delta\xFh^\beta}$.
\end{proof}

\begin{remark}
When $\sqrt{\kappa}U_h\not = L_h$ the norm $\norm[\XSdh^\ast]{\cdot}$ is slightly weaker than $\norm[\XSdh]{\cdot}$. However, this does not alter the well-posedness result stated in Theorem~\ref{thm:thm4}.
\end{remark}

\begin{proposition}
Let $\ks = \km = 0$, $\Pt = I$ and $\sqrt{\kappa}U_h = L_h$. Then, the following residuum orthogonalities hold: 
\begin{alignat*}{3}
& \iprod[]{\zeta\lFh+\nabla\cdot\mFh}{\delta\uFh} =  0\quad && (\text{orthogonality w.r.t. }U_h),\\
&\iprod[]{\eFh+\mFh-\eFt                 }{\delta\eFh} = 0\quad && (\text{orthogonality w.r.t. } E_h),\\
&\iprod[]{\zeta\uFh+\nabla\cdot\sFh-q}{\delta\lFh} = 0\quad && (\text{orthogonality w.r.t. }L_h) .
\end{alignat*}
\end{proposition}

\begin{proof}
First, by setting $\ks = \km = 0$, it follows that
$$\iprod[]{\zeta\lFh+\zeta\tl\Pt(q-\zeta\uFh-\nabla\cdot \sFh)+\nabla\cdot\mFh}{\delta\uFh} = 0$$
and
$$\iprod[]{\zeta\uFh+\zeta\tu\Pt(-\zeta\lFh-\nabla\cdot\mFh)+\nabla\cdot\sFh-q}{\delta\lFh} = 0.$$
Second, by setting $\Pt = I$, we have
$$\iprod[]{\zeta\lFh+\zeta\tl(q-\zeta\uFh-\nabla\cdot\sFh)+\nabla\cdot\mFh}{\delta\uFh} = 0$$
and
$$\iprod[]{\zeta\uFh+\zeta\tu(-\zeta\lFh-\nabla\cdot\mFh)+\nabla\cdot\sFh-q}{\delta\lFh} = 0.$$
By choosing
$$\delta\lFh = \zeta\tl\delta\uFh$$
and after some algebraic manipulations on the two last equations, we arrive at
$$\iprod[]{\zeta\lFh+\nabla\cdot\mFh}{\delta\uFh} =  0,$$
which is proposition’s first orthogonality result.
Similarly, by choosing
$$\delta\uFh = \zeta\tu\delta\lFh$$
and after some algebraic manipulations, we arrive at
$$\iprod[]{\zeta\uFh+\nabla\cdot\sFh-q}{\delta\lFh} = 0,$$
which is proposition’s third orthogonality result. Finally, considering $P_{\mathrm{d},h}$ and noting that $\mF' = 0$ (by design), we arrive at
$$
\iprod[]{\eFh-\mFh-\eFt}{\delta\eFh}=0,$$
which is proposition's second orthogonality result.
\end{proof}
\begin{proposition}[minimally stabilized discrete dual formulation]
Let $h$ be small enough with respect to $\length$. For $\ks=\km=0$ and $U_h = L_h$, the following assertions hold true: \textit{i}) the stabilization term is mesh-independent, i.e., $\Sdh(\xFh,\delta\xFh)$ can be written as $\Sd(\xFh,\delta\xFh)$ with
\begin{align*}
\Sd(\xFh,\delta\xFh)
&:=-\ku \length^2\left(\zeta^2\iprod[]{\lFh}{\delta\lFh}+\zeta\left(\iprod[]{\nabla\cdot\mFh}{\delta\lFh}+\iprod[]{\lFh}{\nabla\cdot\delta\mFh}\right)+\iprod[]{\nabla\cdot\mFh}{\nabla\cdot\delta \mFh}\right)\\
&\quad\:\:-\ke\left(\iprod[]{\eFh}{\delta\eFh}+\iprod[]{\mFh}{\delta\eFh}+\iprod[]{\eFh}{\delta \mFh}+\iprod[]{\mFh}{\delta\mFh}\right)\\
&\quad\:\:+\frac{\kl}{\zeta^2\length^2}\left(\zeta^2\iprod[]{\uFh}{\delta\uFh}+\zeta\left(\iprod[]{\nabla\cdot\sFh}{\delta\uFh}+\iprod[]{\uFh}{\nabla\cdot\delta\sFh}\right)+\iprod[]{\nabla\cdot\sFh}{\nabla\cdot\delta \sFh}\right);
\end{align*}
and, \textit{ii}) the ASGS and the OSGS formulations are equivalent.
\end{proposition}
\begin{proof}
    The first assertion follows intermediately. The second is a direct consequence of the fact that all subscales are (by design) either zero  or orthogonal to the corresponding finite element space (due to the previous proposition).
\end{proof}

\section{Numerical results}

This section presents a numerical assessment of both the primal and dual formulations introduced above. Two numerical examples are considered. The first aims to evaluate the convergence properties of the methods by computing relevant error norms against a smooth manufactured solution. These results will be compared with those from a \textit{Natural} formulation, a \textit{Fully stabilized} equal-order formulation based on GLS and an \textit{Unstabilized} equal-order formulation, all recently discussed in \cite{Bazon2025}.

The second example addresses a more challenging problem featuring an internal reaction boundary layer, in order to assess the approximation properties of the proposed methods in the presence of large solution gradients. This case clearly illustrates key differences between the primal and dual formulations and further allows us to examine the effect of different choices for the data fields $(\eFt,\sFt)$ making this final experiment particularly relevant for DDCM applications.

We consider a sequence of conforming, shape-regular finite element partitions $\mathcal{P}_h$ of the computational domain, with characteristic mesh size $h$. A primary motivation for the VMS approach is the ability to use equal-order approximations. Accordingly, we employ the following choice of discrete finite element spaces:
\begin{equation*}\label{eq:equal_order_func_spaces}
	\left(\uF_{h}, \eF_{h}, \sF_{h}, \lF_{h}, \mF_{h}\right)
	\in
	\CG_{k} \times \left[\CG_{k}\right]^{\dimo}
	\times \left[\CG_{k}\right]^{\dimo}
	\times \CG_{k}
	\times \left[\CG_{k}\right]^{\dimo},
\end{equation*}
where $\dimo$ denotes the spatial dimension of the computational domain $\Omega$, $k>0$ and 
\begin{equation*}
	\quad \CG_{k} := \left\{ v \in C^0(\Omega) \;\middle|\; v|_K \in \mathbb{P}_k(K),\ \forall K \in \mathcal{P}_h \right\}	
\end{equation*}
where $\mathbb{P}_k(K)$ denotes the space of polynomials of degree $k$ defined on each element $K \in \mathcal{P}_h$.
In the ASGS method, being $\Pt$ the identity operator, no additional unknowns are
required, since the subscales are simply evaluated locally at the element level. In contrast, 
for the OSGS method, with $\Pt$ denoting the $\Ltwo$-projector onto the orthogonal complement of the 
finite element space, one additional unknown per field must be considered.
That is, if we denote by $w \in W$ any of the unknown fields, such that $w = w_h + w'$ with $w_h\in W_h \subsetneq W$, 
the subscale component $w'$ is in this case given by
\begin{equation}
w' = \tau_w \, \Pt \left ( r(w) \right ) = \tau_w \, \left [ r_w - P_{W_{h}} \left ( r_w \right ) \right ]
\end{equation}
where $r_w$ is the algebraic residual equation of the variable $w$ and $P_{W_{h}}$ stands for the $\Ltwo$-projection onto $W_h$.
Since this leads to an excessive computational burden, in practice the explicit computation 
of the subscales is avoided by applying some sort of fixed-point iteration scheme. 
In our implementation, these variables are included as unknowns and solved for 
monolithically. This approach yields cleaner code, simplifies the implementation 
and eliminates the numerical error associated with an additional iterative procedure. 
We anticipate, however, that no significant differences in results were observed between the 
ASGS and OSGS methods for either the primal or dual formulations. 
Consequently, all results presented hereafter are obtained by using the ASGS method. Finally, the implementation of both formulations, including the two variants,
ASGS and OSGS, was carried out using the finite element platform Firedrake \cite{FiredrakeUserManual} and is available in \url{https://github.com/pedrobbazon/DDCM.git}.

\subsection{Example 1: Convergence to a manufactured solution}

Let us first describe the construction of a manufactured solution with respect to which the relevant error norms will be computed. Its construction deserves some attention, although it essentially follows the idea presented in \cite{Bazon2025}. By examining the optimality conditions in their strong form, it is straightforward to devise a manufactured solution that resembles a relevant scenario for DDCM.

\begin{figure}[ht!]
	\centering
	\caption{Manufactured solution for the example 1 showing the primal and dual fields as well as the data fields $\sFt$ and $\eFt$.}
	~\\
	\def\svgwidth{1.0\linewidth}
	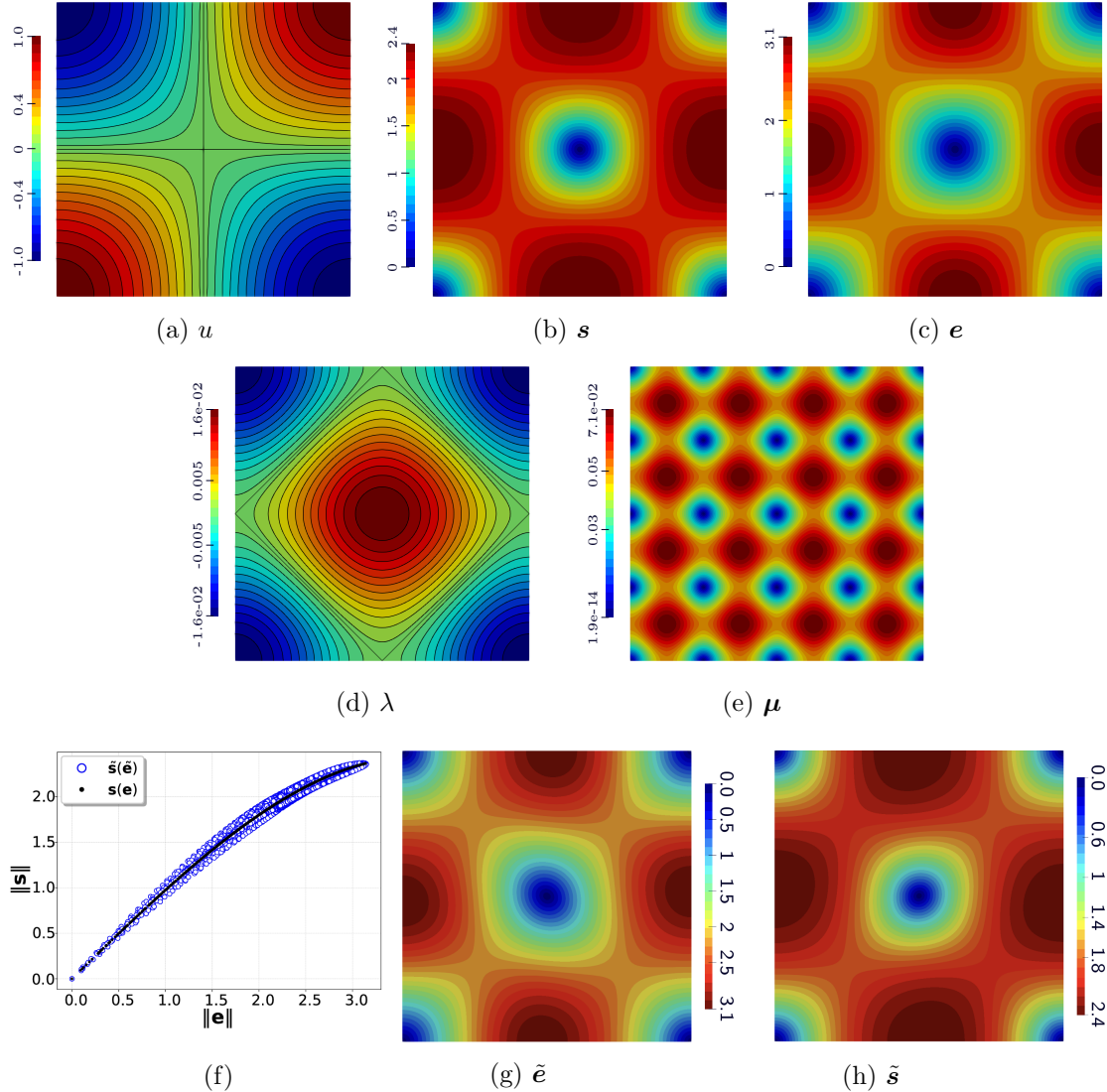
	\label{fig:reference_fields_solution}
\end{figure}

The manufactured solution is illustrated in Figure~\ref{fig:reference_fields_solution}, which
shows all the primal fields $(\uF, \sF, \eF)$, the dual fields $(\lF,\mF)$, and the data fields
$(\sFt,\eFt)$. To construct this solution, we begin by defining the primal scalar field
\begin{equation}
 \uF(x,y) = \cos(\pi x)\cos(\pi y),
\end{equation}
whose gradient is given by
\begin{equation}
	\eF(x,y) = -\pi \,[\sin(\pi x)\cos(\pi y),\, \cos(\pi x)\sin(\pi y)]^{\top}.
\end{equation}
We then define $\eFt$ as a perturbation of $\eF$ as follows
\begin{equation}
\eFt(x,y) = \eF(x,y) + \frac{1}{20} \left[\sin\left(4\pi x \right), \sin\left(4 \pi y \right) \right]^{\top}
\end{equation}
Then, for the primal flux $\sF$ we consider 
\begin{equation}\label{eq:nonlinear_fluxlaw}
	\sF = -\frac{\eF}{\|\eF\|}\left(\|\eF\| - \frac{1}{40}\,\|\eF\|^{3}\right).
\end{equation}
such as to mimic a nonlinear constitutive relation between $\sF$ and $\eF$, as plotted in Figure~\ref{fig:reference_fields_solution}~(f), which shows the corresponding relation between the magnitude of $\sF$ and the magnitude of $\eF$ obtained by sampling the solution over a computational mesh.
For the dual scalar field $\lambda$ we take the function
\begin{equation}
	\lF(x,y) = -\frac{1}{40\pi}\,[\cos(2\pi x) + \cos(2\pi y)],
\end{equation}
and finally from the equation for $\eF$ we simply obtain $\mF$ as
\begin{equation}\label{eq:formu}
\mF = \eFt - \eF
\end{equation}
We need to introduce a (dual) source term on the right-hand-side of the equation for $u$, which must satisfy
\begin{equation}
	\zeta \lF + \nabla \cdot \mF = f,
\end{equation}
but note that this is done solely for the purpose of assessing the convergence 
properties of the method using such a manufactured solution.
Finally, $\sFt$ is obtained from
\begin{equation}
	\sFt = \sF - \kappa^{-1} \, \nabla{\lF}
\end{equation}
and the source term $q$ is obtained by substituting the proposed fields into
the balance equation
\begin{equation}
	q = \zeta \uF + \nabla \cdot \sF.
\end{equation}
Figure~\ref{fig:reference_fields_solution}~(f) also displays, using blue dots, the sampled values of the magnitude of $\sFt$ as a function of the magnitude of $\eFt$, taken over the computational mesh. The computational domain $\Omega$ is chosen as the unit square $[0,1]^2$, with the scale parameter and reaction coefficient both set to unity, i.e., $\kappa = 1$ and $\zeta = 1$.

\subsubsection*{Boundary conditions}

Regarding boundary conditions on $\partial \Omega$, we consider two scenarios to ensure a fair comparison between the primal and dual formulations:
In the first case, we assume the entire boundary is of Dirichlet type and strongly impose the boundary conditions on the scalar unknowns $\uF$ and $\lF$, i.e.,\begin{equation}
	\left\lbrace
	\begin{array}{rclll}
		\uF &=& \mathsf{g}_{\uF} &\mbox{on}&~\partial \Omega, \\
		\lF &=& \mathsf{g}_{\lF} &\mbox{on}&~\partial \Omega, \\
	\end{array}
	\right.
\end{equation}
while in the latter case, we assume the entire boundary is of Neumann type and strongly impose the boundary conditions on the normal component of the vector unknowns $\sF$ and $\mF$, i.e.,
\begin{equation}
	\left\lbrace
	\begin{array}{rclll}
		\sF \cdot \normal  &=& \mathsf{h}_{\sF} &\mbox{on}&~\partial \Omega, \\
		\mF \cdot \normal &=& \mathsf{h}_{\mF} &\mbox{on}&~\partial \Omega, \\
	\end{array}
	\right.
\end{equation}
where the values on the right-hand-side are evaluated  from the manufactured fields
at the locations of the corresponding degrees of freedom.

\subsubsection*{Choice of estabilization parameters}

From Sections \ref{sec:anal_primal} and \ref{sec:anal_dual}, we observe considerable flexibility in selecting the stabilization parameters for both the primal and dual formulations. 
In this example, for the primal formulation, the stabilization parameters are chosen as
$$
k_{\uF} = k_{\lF} = 0, 
\quad 
k_{\eF} = k_{\sF} = k_{\mF} = \tfrac{1}{8}.
$$
This choice is convenient, as it reduces the number of parameters to be selected to just three. It also exploits the residual orthogonality property established earlier, ensuring that the results are independent of the choice of the projector operator $\Pt$.
On the other hand, for the dual formulation we take
$$
k_{\uF} = k_{\lF} = k_{\eF} = k_{\mF} = \tfrac{1}{4}, 
\quad 
k_{\sF} = 1.
$$
Although we could also have selected $k_{\mF} = k_{\sF} = 0$, the choice $k_{\mF} = \tfrac{1}{4}$ makes the selection of $k_{\eF}$ mesh-independent. This is convenient for assessing the convergence rates of the method across the full range of mesh refinements tested.
As a final remark, the convergence studies presented below employ linear and quadratic polynomial approximations (i.e., $k \in \{1,2\}$).
Therefore, in the following discussion, we consider the algebraic implementation of the methods, which we refer to as \textit{Primal-ASGS} and \textit{Dual-ASGS}.

\begin{figure}[ht!!]
	\centering
	\caption{Contours of fields $\uF$, $\sF$, $\eF$, $\lF$ and $\mF$ 
	for the \textit{Primal-ASGS} formulation for different levels of mesh refinement.}
	~\\
	\def\svgwidth{0.7\columnwidth}
	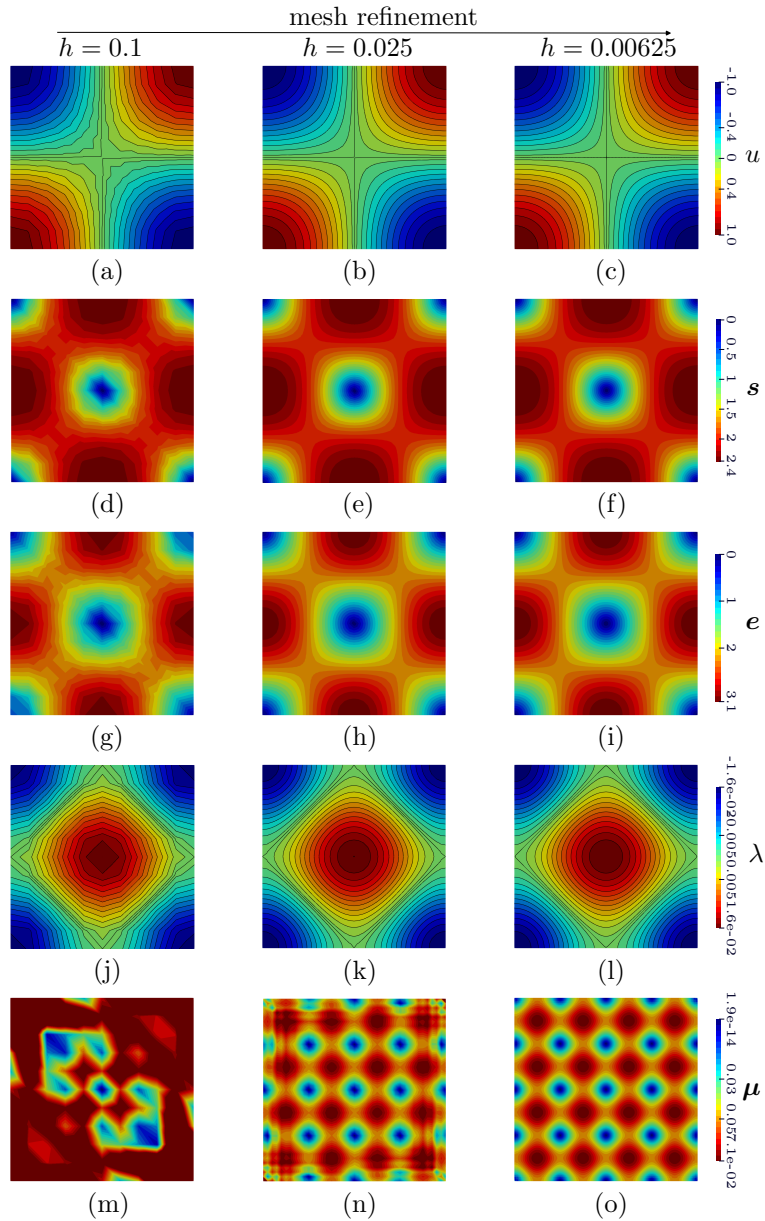
	\label{fig:primal_formulation_fields_solution}
\end{figure}

\begin{figure}[ht!!]
	\centering
	\caption{Contours of fields $\uF$, $\sF$, $\eF$, $\lF$ and $\mF$ 
	for the \textit{Dual-ASGS} for different levels of mesh refinement.}
	~\\
	\def\svgwidth{0.7\columnwidth}
	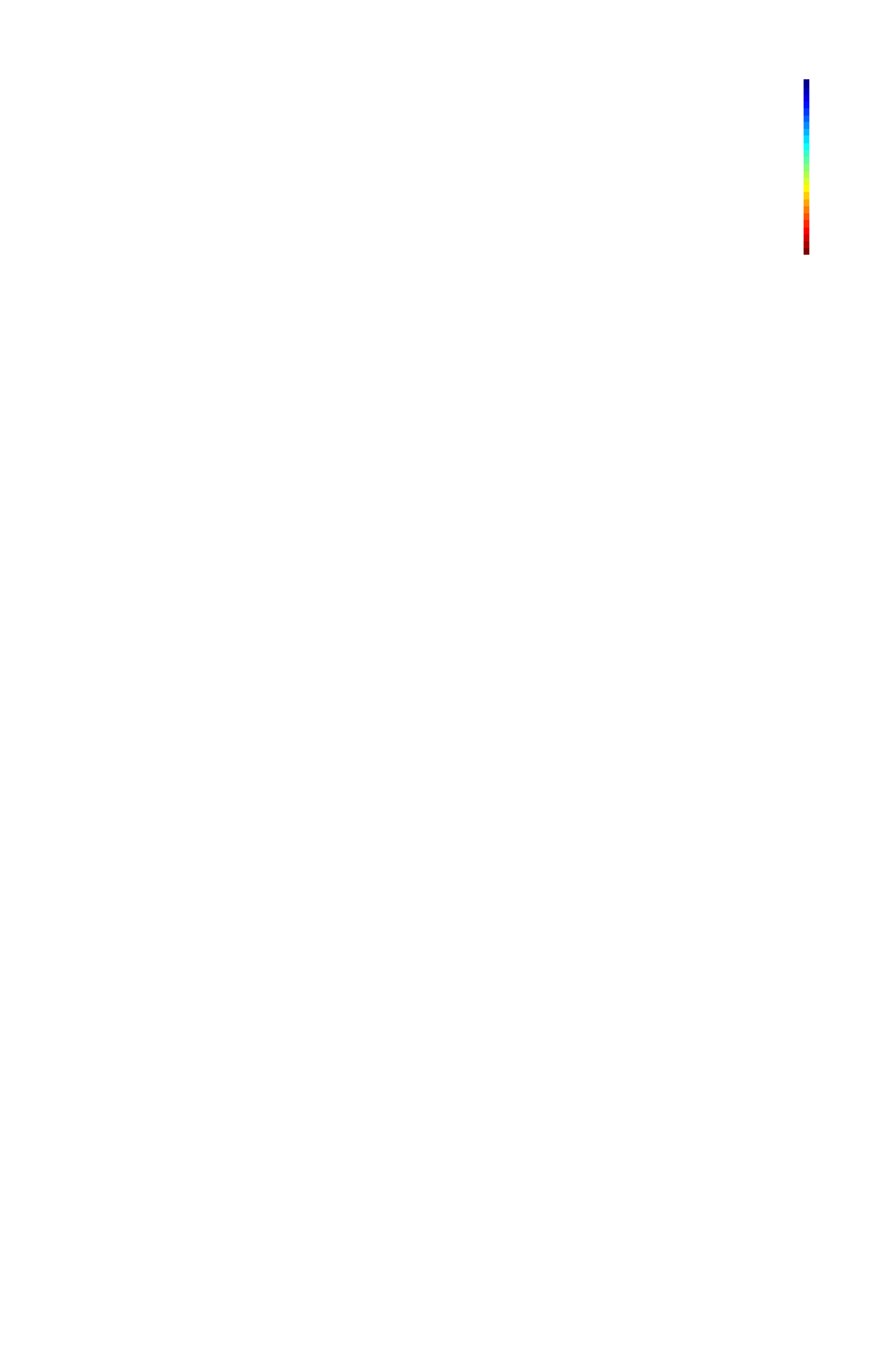
	\label{fig:dual_formulation_fields_solution}
\end{figure}

We begin the numerical analysis with a qualitative comparison of the unknown fields across successive mesh refinements. Three refinement levels are examined: a very coarse mesh ($h = 0.1$), an intermediate mesh ($h = 0.025$), and a very fine mesh ($h = 0.00625$). The solutions for the \textit{Primal-ASGS} and \textit{Dual-ASGS} formulations are shown in Figures~\ref{fig:primal_formulation_fields_solution} and~\ref{fig:dual_formulation_fields_solution}, respectively. These results should be compared against the reference solution in Figure~\ref{fig:reference_fields_solution}.
A eye-catching observation from these plots is that both formulations perform well on the finest mesh. The only minor exception is the field $\lF$, which exhibits slight numerical artifacts near the boundaries, see Figure~\ref{fig:dual_formulation_fields_solution}~(l). These artifacts, however, diminish with further refinement and do not impact the convergence of the method, as will be demonstrated. For the coarsest and intermediate mesh levels, the differences between the two formulations become more evident. While the primal formulation does a reasonable job of capturing the scalar fields, it struggles to accurately resolve particularly the $\mF$ field, Figure~\ref{fig:primal_formulation_fields_solution}~(m--n). As expected, the situation is reversed for the dual formulation, which performs significantly better in capturing
$\mF$ but exhibits poorer accuracy for the scalar fields, see Figure~\ref{fig:dual_formulation_fields_solution}~(a,~j--k), where this behavior is most evident).
This qualitative assessment suggests that both formulations capture the field $\sF$ well, a finding later corroborated by the $\Ltwo$-error norm in the numerical convergence analysis. However, the situation differs when the error is measured in the $\Hdiv$-norm, a distinction not apparent from this preliminary evaluation.

\begin{figure}[h!]
	\centering
	\caption{Convergence of scalar fields $\uF$ and $\lF$ for $k = 1$. Figures (a--b) show the results for $\uF$ and (c--d) for $\lF$, considering all formulations: \textit{Natural}, equal-order \textit{Unstabilized}, equal-order \textit{Fully stabilized}, \textit{Primal-ASGS} and \textit{Dual-ASGS}. Errors are evaluated in the $\Ltwo$- and $\Hone$-norms, as indicated. For the \textit{Dual-ASGS} formulation, only the $\Ltwo$-norm is assessed.}
		~\\
	\def\svgwidth{1.0\columnwidth}
\begingroup%
  \makeatletter%
  \providecommand\color[2][]{%
    \errmessage{(Inkscape) Color is used for the text in Inkscape, but the package 'color.sty' is not loaded}%
    \renewcommand\color[2][]{}%
  }%
  \providecommand\transparent[1]{%
    \errmessage{(Inkscape) Transparency is used (non-zero) for the text in Inkscape, but the package 'transparent.sty' is not loaded}%
    \renewcommand\transparent[1]{}%
  }%
  \providecommand\rotatebox[2]{#2}%
  \newcommand*\fsize{\dimexpr\f@size pt\relax}%
  \newcommand*\lineheight[1]{\fontsize{\fsize}{#1\fsize}\selectfont}%
  \ifx\svgwidth\undefined%
    \setlength{\unitlength}{584.02459573bp}%
    \ifx\svgscale\undefined%
      \relax%
    \else%
      \setlength{\unitlength}{\unitlength * \real{\svgscale}}%
    \fi%
  \else%
    \setlength{\unitlength}{\svgwidth}%
  \fi%
  \global\let\svgwidth\undefined%
  \global\let\svgscale\undefined%
  \makeatother%
  \begin{picture}(1,0.83006369)%
    \lineheight{1}%
    \setlength\tabcolsep{0pt}%
    \put(0.25747727,0.00265722){\color[rgb]{0,0,0}\makebox(0,0)[lt]{\lineheight{1.25}\smash{\begin{tabular}[t]{l}(c)\end{tabular}}}}%
    \put(0.25684251,0.42916618){\color[rgb]{0,0,0}\makebox(0,0)[lt]{\lineheight{1.25}\smash{\begin{tabular}[t]{l}(a)\end{tabular}}}}%
    \put(0.76340203,0.00265722){\color[rgb]{0,0,0}\makebox(0,0)[lt]{\lineheight{1.25}\smash{\begin{tabular}[t]{l}(d)\end{tabular}}}}%
    \put(0.76340203,0.42915962){\color[rgb]{0,0,0}\makebox(0,0)[lt]{\lineheight{1.25}\smash{\begin{tabular}[t]{l}(b)\end{tabular}}}}%
    \put(0,0){\includegraphics[width=\unitlength,page=1]{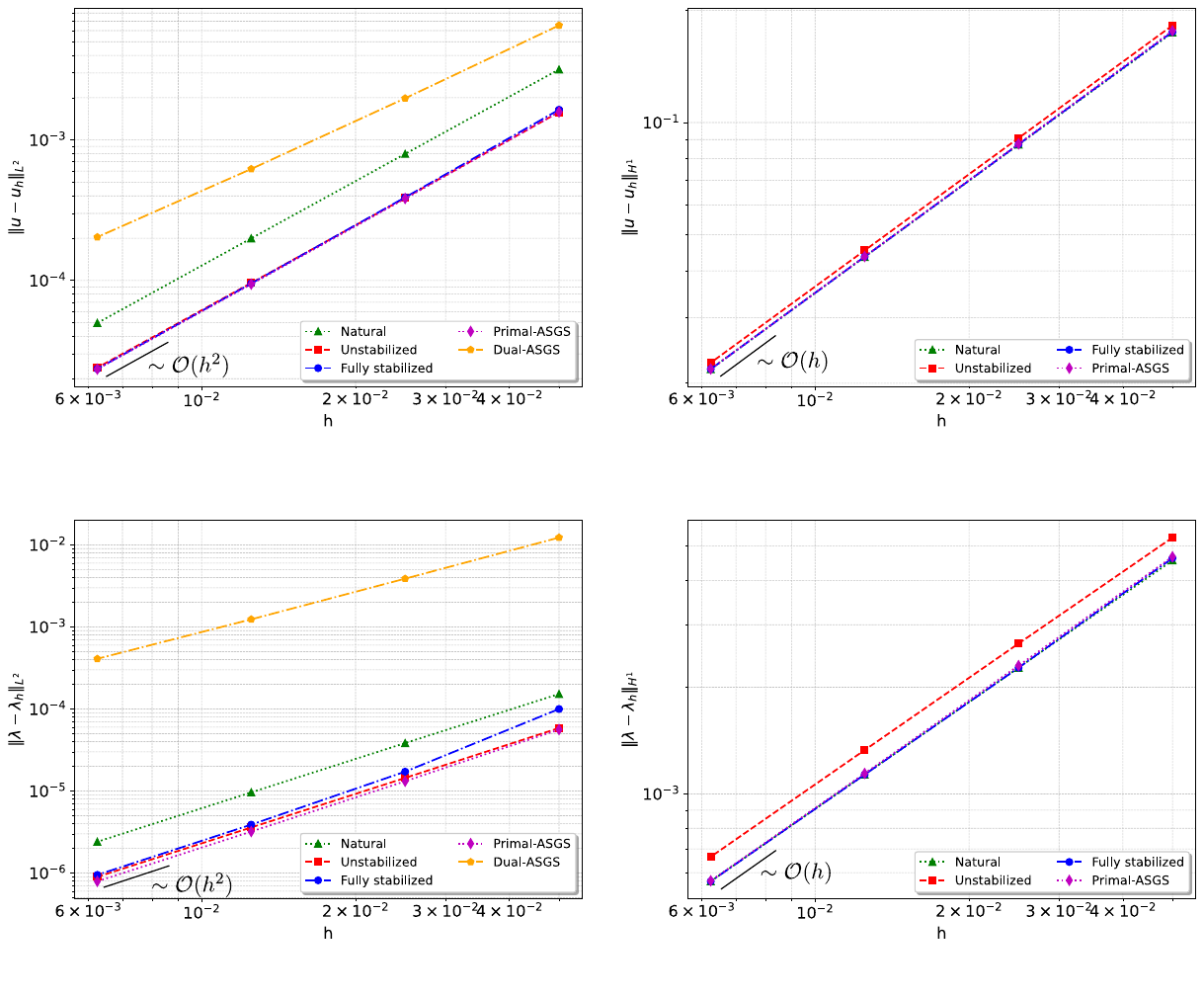}}%
  \end{picture}%
\endgroup%

	\label{fig:convergence_scalar_linear_art02}
\end{figure}

\begin{figure}[h!]
	\centering
	\caption{Convergence of vector fields $\sF$, $\eF$ and $\mF$ with $k=1$, evaluated in the $\Ltwo$-norm considering all formulations: \textit{Natural}, equal-order \textit{Unstabilized}, equal-order \textit{Fully stabilized}, \textit{Primal-ASGS} and \textit{Dual-ASGS}.
	}
	~\\
	\def\svgwidth{1.0\columnwidth}
\begingroup%
  \makeatletter%
  \providecommand\color[2][]{%
    \errmessage{(Inkscape) Color is used for the text in Inkscape, but the package 'color.sty' is not loaded}%
    \renewcommand\color[2][]{}%
  }%
  \providecommand\transparent[1]{%
    \errmessage{(Inkscape) Transparency is used (non-zero) for the text in Inkscape, but the package 'transparent.sty' is not loaded}%
    \renewcommand\transparent[1]{}%
  }%
  \providecommand\rotatebox[2]{#2}%
  \newcommand*\fsize{\dimexpr\f@size pt\relax}%
  \newcommand*\lineheight[1]{\fontsize{\fsize}{#1\fsize}\selectfont}%
  \ifx\svgwidth\undefined%
    \setlength{\unitlength}{581.83230423bp}%
    \ifx\svgscale\undefined%
      \relax%
    \else%
      \setlength{\unitlength}{\unitlength * \real{\svgscale}}%
    \fi%
  \else%
    \setlength{\unitlength}{\svgwidth}%
  \fi%
  \global\let\svgwidth\undefined%
  \global\let\svgscale\undefined%
  \makeatother%
  \begin{picture}(1,0.8161779)%
    \lineheight{1}%
    \setlength\tabcolsep{0pt}%
    \put(0.49834612,0.00266711){\color[rgb]{0,0,0}\makebox(0,0)[lt]{\lineheight{1.25}\smash{\begin{tabular}[t]{l}(c)\end{tabular}}}}%
    \put(0.26257038,0.41898157){\color[rgb]{0,0,0}\makebox(0,0)[lt]{\lineheight{1.25}\smash{\begin{tabular}[t]{l}(a)\end{tabular}}}}%
    \put(0.77387368,0.41897168){\color[rgb]{0,0,0}\makebox(0,0)[lt]{\lineheight{1.25}\smash{\begin{tabular}[t]{l}(b)\end{tabular}}}}%
    \put(0,0){\includegraphics[width=\unitlength,page=1]{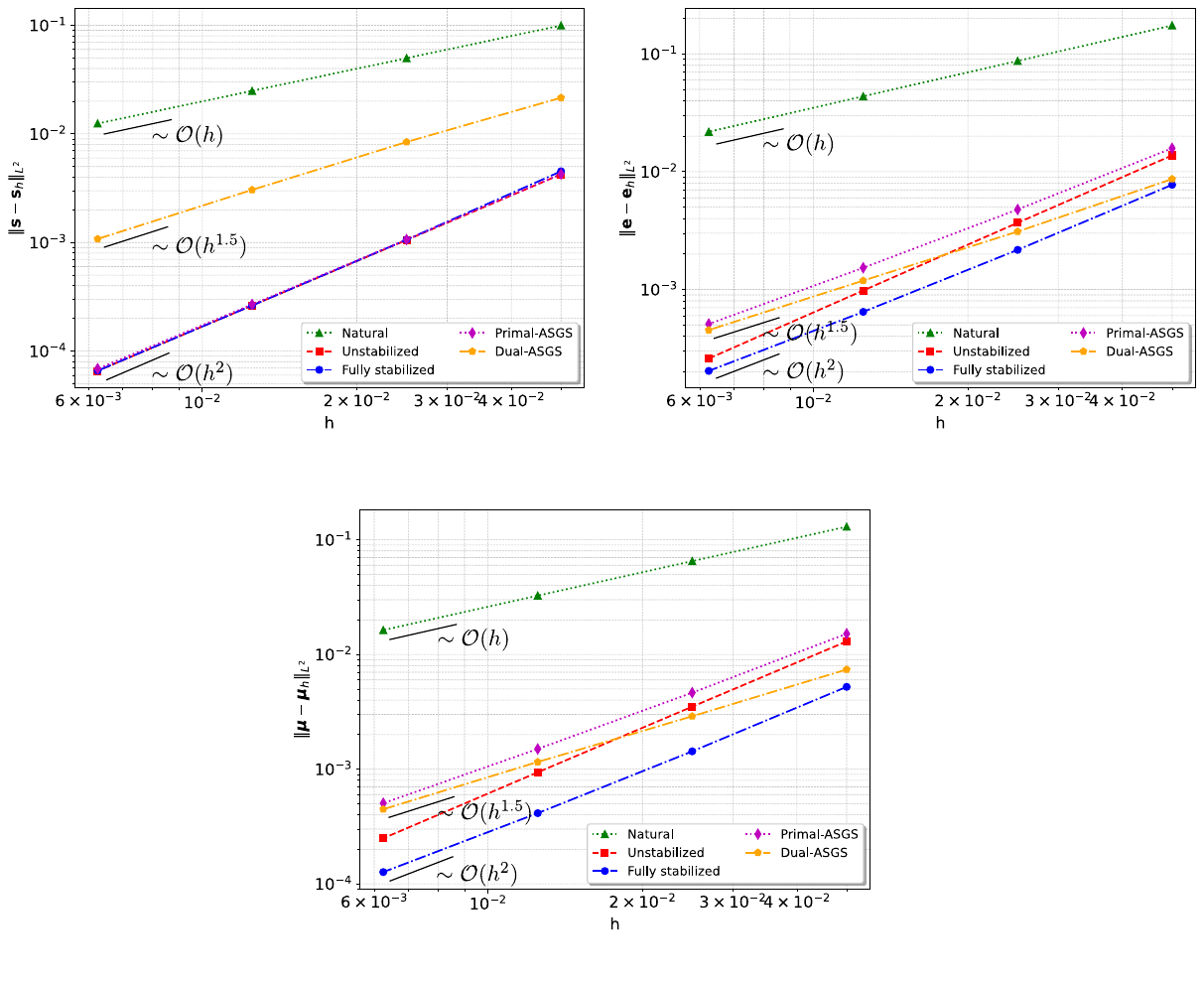}}%
  \end{picture}%
\endgroup%

	\label{fig:convergence_vector_linear_art02}
\end{figure}

\begin{figure}[h!]
	\centering
	\caption{Convergence of vector fields $\sF$, $\eF$ and $\mF$ with $k=1$, evaluated in the $\Hdiv$-norm for equal-order \textit{Unstabilized}, equal-order \textit{Fully stabilized}, \textit{Primal-ASGS} and \textit{Dual-ASGS} formulations.
	}
	~\\
	\def\svgwidth{1.0\columnwidth}
\begingroup%
  \makeatletter%
  \providecommand\color[2][]{%
    \errmessage{(Inkscape) Color is used for the text in Inkscape, but the package 'color.sty' is not loaded}%
    \renewcommand\color[2][]{}%
  }%
  \providecommand\transparent[1]{%
    \errmessage{(Inkscape) Transparency is used (non-zero) for the text in Inkscape, but the package 'transparent.sty' is not loaded}%
    \renewcommand\transparent[1]{}%
  }%
  \providecommand\rotatebox[2]{#2}%
  \newcommand*\fsize{\dimexpr\f@size pt\relax}%
  \newcommand*\lineheight[1]{\fontsize{\fsize}{#1\fsize}\selectfont}%
  \ifx\svgwidth\undefined%
    \setlength{\unitlength}{588.57786644bp}%
    \ifx\svgscale\undefined%
      \relax%
    \else%
      \setlength{\unitlength}{\unitlength * \real{\svgscale}}%
    \fi%
  \else%
    \setlength{\unitlength}{\svgwidth}%
  \fi%
  \global\let\svgwidth\undefined%
  \global\let\svgscale\undefined%
  \makeatother%
  \begin{picture}(1,0.80682394)%
    \lineheight{1}%
    \setlength\tabcolsep{0pt}%
    \put(0.50409535,0.00263656){\color[rgb]{0,0,0}\makebox(0,0)[lt]{\lineheight{1.25}\smash{\begin{tabular}[t]{l}(c)\end{tabular}}}}%
    \put(0.2710218,0.41417973){\color[rgb]{0,0,0}\makebox(0,0)[lt]{\lineheight{1.25}\smash{\begin{tabular}[t]{l}(a)\end{tabular}}}}%
    \put(0.77646493,0.41417002){\color[rgb]{0,0,0}\makebox(0,0)[lt]{\lineheight{1.25}\smash{\begin{tabular}[t]{l}(b)\end{tabular}}}}%
    \put(0,0){\includegraphics[width=\unitlength,page=1]{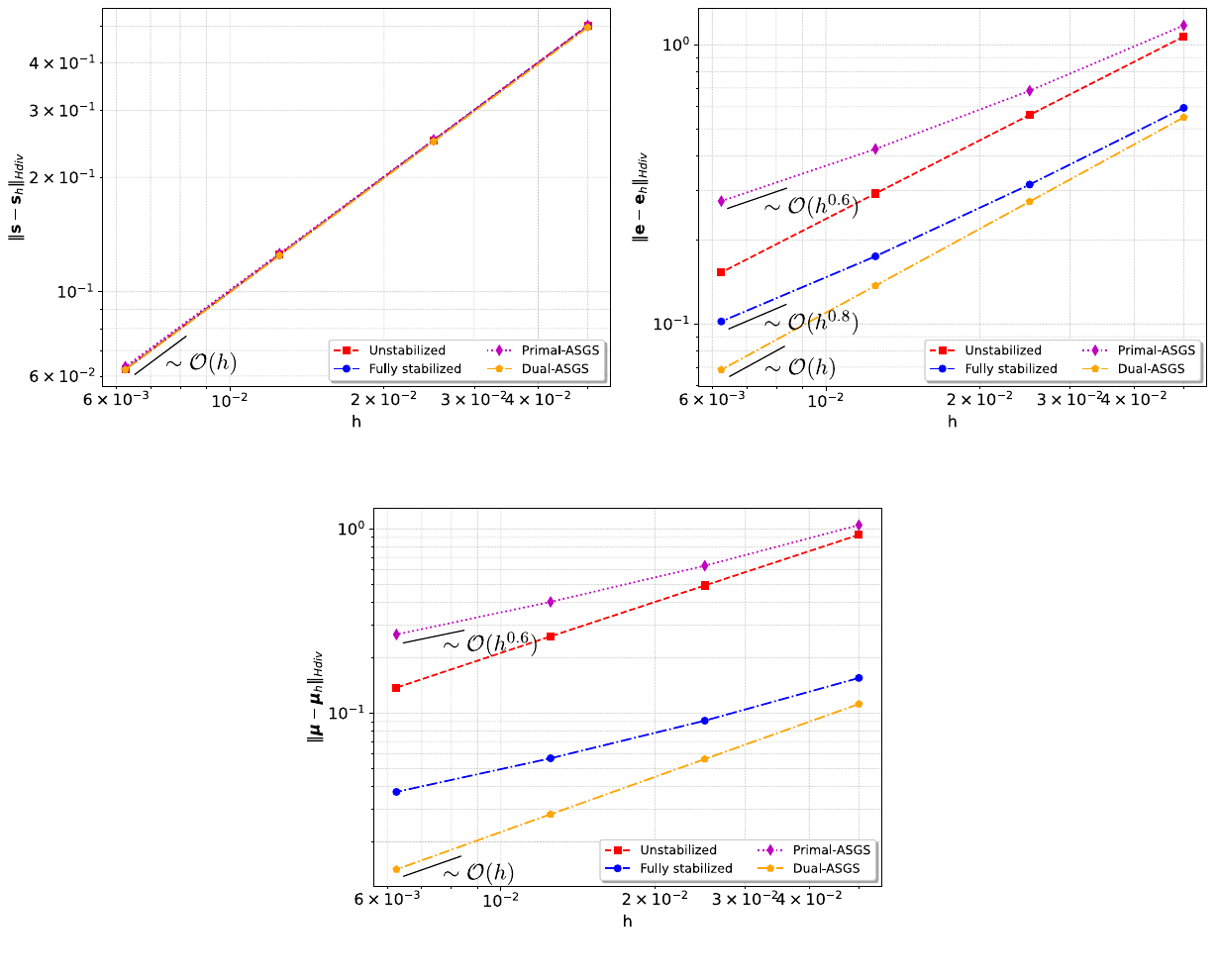}}%
  \end{picture}%
\endgroup%

	\label{fig:convergence_vector_linear_hdiv_art02}
\end{figure}

\begin{figure}[h!]
	\centering
	\caption{Convergence of scalar fields $\uF$ and $\lF$ for $k = 2$. Figures (a--b) show the results for $\uF$, and (c--d) for $\lF$, considering all formulations: \textit{Natural}, equal-order \textit{Unstabilized}, equal-order \textit{Fully stabilized}, \textit{Primal-ASGS} and \textit{Dual-ASGS}. Errors are evaluated in the $\Ltwo$- and $\Hone$-norms, as indicated. For the \textit{Dual-ASGS} formulation, only the $\Ltwo$-norm is assessed.}
	~\\
	\def\svgwidth{1.0\columnwidth}
\begingroup%
  \makeatletter%
  \providecommand\color[2][]{%
    \errmessage{(Inkscape) Color is used for the text in Inkscape, but the package 'color.sty' is not loaded}%
    \renewcommand\color[2][]{}%
  }%
  \providecommand\transparent[1]{%
    \errmessage{(Inkscape) Transparency is used (non-zero) for the text in Inkscape, but the package 'transparent.sty' is not loaded}%
    \renewcommand\transparent[1]{}%
  }%
  \providecommand\rotatebox[2]{#2}%
  \newcommand*\fsize{\dimexpr\f@size pt\relax}%
  \newcommand*\lineheight[1]{\fontsize{\fsize}{#1\fsize}\selectfont}%
  \ifx\svgwidth\undefined%
    \setlength{\unitlength}{584.02463898bp}%
    \ifx\svgscale\undefined%
      \relax%
    \else%
      \setlength{\unitlength}{\unitlength * \real{\svgscale}}%
    \fi%
  \else%
    \setlength{\unitlength}{\svgwidth}%
  \fi%
  \global\let\svgwidth\undefined%
  \global\let\svgscale\undefined%
  \makeatother%
  \begin{picture}(1,0.83075002)%
    \lineheight{1}%
    \setlength\tabcolsep{0pt}%
    \put(0.25747726,0.0026571){\color[rgb]{0,0,0}\makebox(0,0)[lt]{\lineheight{1.25}\smash{\begin{tabular}[t]{l}(c)\end{tabular}}}}%
    \put(0.25684248,0.42916598){\color[rgb]{0,0,0}\makebox(0,0)[lt]{\lineheight{1.25}\smash{\begin{tabular}[t]{l}(a)\end{tabular}}}}%
    \put(0.76340207,0.0026571){\color[rgb]{0,0,0}\makebox(0,0)[lt]{\lineheight{1.25}\smash{\begin{tabular}[t]{l}(d)\end{tabular}}}}%
    \put(0.76340207,0.42915946){\color[rgb]{0,0,0}\makebox(0,0)[lt]{\lineheight{1.25}\smash{\begin{tabular}[t]{l}(b)\end{tabular}}}}%
    \put(0,0){\includegraphics[width=\unitlength,page=1]{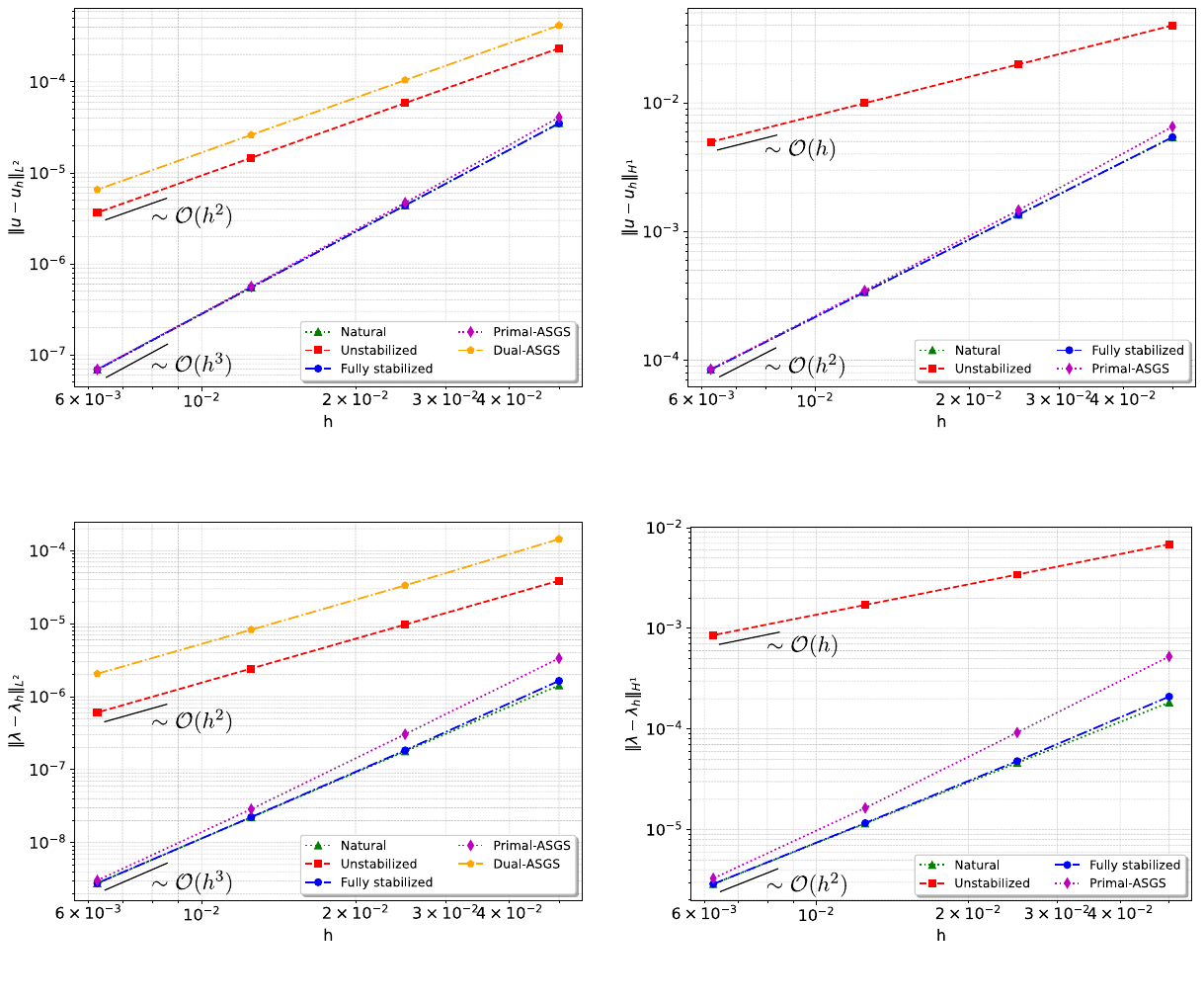}}%
  \end{picture}%
\endgroup%

	\label{fig:convergence_scalar_quadratic_art02}
\end{figure}

\begin{figure}[h!]
	\centering
	\caption{Convergence of vector fields $\sF$, $\eF$ and $\mF$ with $k=2$, evaluated in the $\Ltwo$-norm considering all formulations: \textit{Natural}, equal-order \textit{Unstabilized}, equal-order \textit{Fully stabilized}, \textit{Primal-ASGS} and \textit{Dual-ASGS}.}
	~\\
	\def\svgwidth{1.0\columnwidth}
\begingroup%
  \makeatletter%
  \providecommand\color[2][]{%
    \errmessage{(Inkscape) Color is used for the text in Inkscape, but the package 'color.sty' is not loaded}%
    \renewcommand\color[2][]{}%
  }%
  \providecommand\transparent[1]{%
    \errmessage{(Inkscape) Transparency is used (non-zero) for the text in Inkscape, but the package 'transparent.sty' is not loaded}%
    \renewcommand\transparent[1]{}%
  }%
  \providecommand\rotatebox[2]{#2}%
  \newcommand*\fsize{\dimexpr\f@size pt\relax}%
  \newcommand*\lineheight[1]{\fontsize{\fsize}{#1\fsize}\selectfont}%
  \ifx\svgwidth\undefined%
    \setlength{\unitlength}{581.83230423bp}%
    \ifx\svgscale\undefined%
      \relax%
    \else%
      \setlength{\unitlength}{\unitlength * \real{\svgscale}}%
    \fi%
  \else%
    \setlength{\unitlength}{\svgwidth}%
  \fi%
  \global\let\svgwidth\undefined%
  \global\let\svgscale\undefined%
  \makeatother%
  \begin{picture}(1,0.81617798)%
    \lineheight{1}%
    \setlength\tabcolsep{0pt}%
    \put(0.49834611,0.0026672){\color[rgb]{0,0,0}\makebox(0,0)[lt]{\lineheight{1.25}\smash{\begin{tabular}[t]{l}(c)\end{tabular}}}}%
    \put(0.26257038,0.41898169){\color[rgb]{0,0,0}\makebox(0,0)[lt]{\lineheight{1.25}\smash{\begin{tabular}[t]{l}(a)\end{tabular}}}}%
    \put(0.7738737,0.41897181){\color[rgb]{0,0,0}\makebox(0,0)[lt]{\lineheight{1.25}\smash{\begin{tabular}[t]{l}(b)\end{tabular}}}}%
    \put(0,0){\includegraphics[width=\unitlength,page=1]{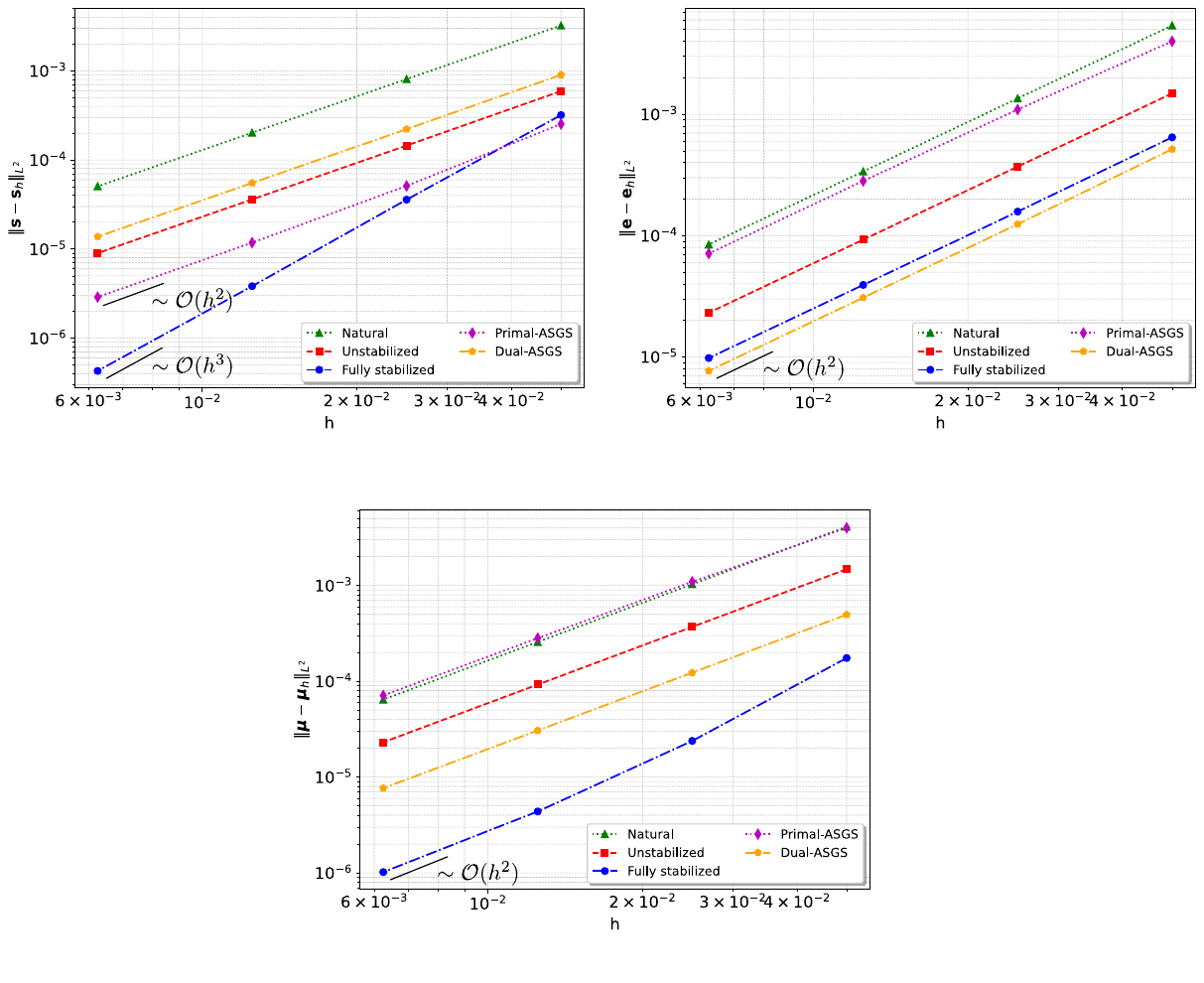}}%
  \end{picture}%
\endgroup%

	\label{fig:convergence_vector_quadratic_art02}
\end{figure}

\begin{figure}[h!]
	\centering
	\caption{Convergence of vector fields $\sF$, $\eF$ and $\mF$ with $k=2$, evaluated in the $\Hdiv$-norm for equal-order \textit{Unstabilized}, equal-order \textit{Fully stabilized}, \textit{Primal-ASGS} and \textit{Dual-ASGS} formulations.}
	~\\
	\def\svgwidth{1.0\columnwidth}
\begingroup%
  \makeatletter%
  \providecommand\color[2][]{%
    \errmessage{(Inkscape) Color is used for the text in Inkscape, but the package 'color.sty' is not loaded}%
    \renewcommand\color[2][]{}%
  }%
  \providecommand\transparent[1]{%
    \errmessage{(Inkscape) Transparency is used (non-zero) for the text in Inkscape, but the package 'transparent.sty' is not loaded}%
    \renewcommand\transparent[1]{}%
  }%
  \providecommand\rotatebox[2]{#2}%
  \newcommand*\fsize{\dimexpr\f@size pt\relax}%
  \newcommand*\lineheight[1]{\fontsize{\fsize}{#1\fsize}\selectfont}%
  \ifx\svgwidth\undefined%
    \setlength{\unitlength}{581.83239073bp}%
    \ifx\svgscale\undefined%
      \relax%
    \else%
      \setlength{\unitlength}{\unitlength * \real{\svgscale}}%
    \fi%
  \else%
    \setlength{\unitlength}{\svgwidth}%
  \fi%
  \global\let\svgwidth\undefined%
  \global\let\svgscale\undefined%
  \makeatother%
  \begin{picture}(1,0.81617786)%
    \lineheight{1}%
    \setlength\tabcolsep{0pt}%
    \put(0.49834605,0.00266712){\color[rgb]{0,0,0}\makebox(0,0)[lt]{\lineheight{1.25}\smash{\begin{tabular}[t]{l}(c)\end{tabular}}}}%
    \put(0.26257034,0.41898148){\color[rgb]{0,0,0}\makebox(0,0)[lt]{\lineheight{1.25}\smash{\begin{tabular}[t]{l}(a)\end{tabular}}}}%
    \put(0.77387357,0.41897167){\color[rgb]{0,0,0}\makebox(0,0)[lt]{\lineheight{1.25}\smash{\begin{tabular}[t]{l}(b)\end{tabular}}}}%
    \put(0,0){\includegraphics[width=\unitlength,page=1]{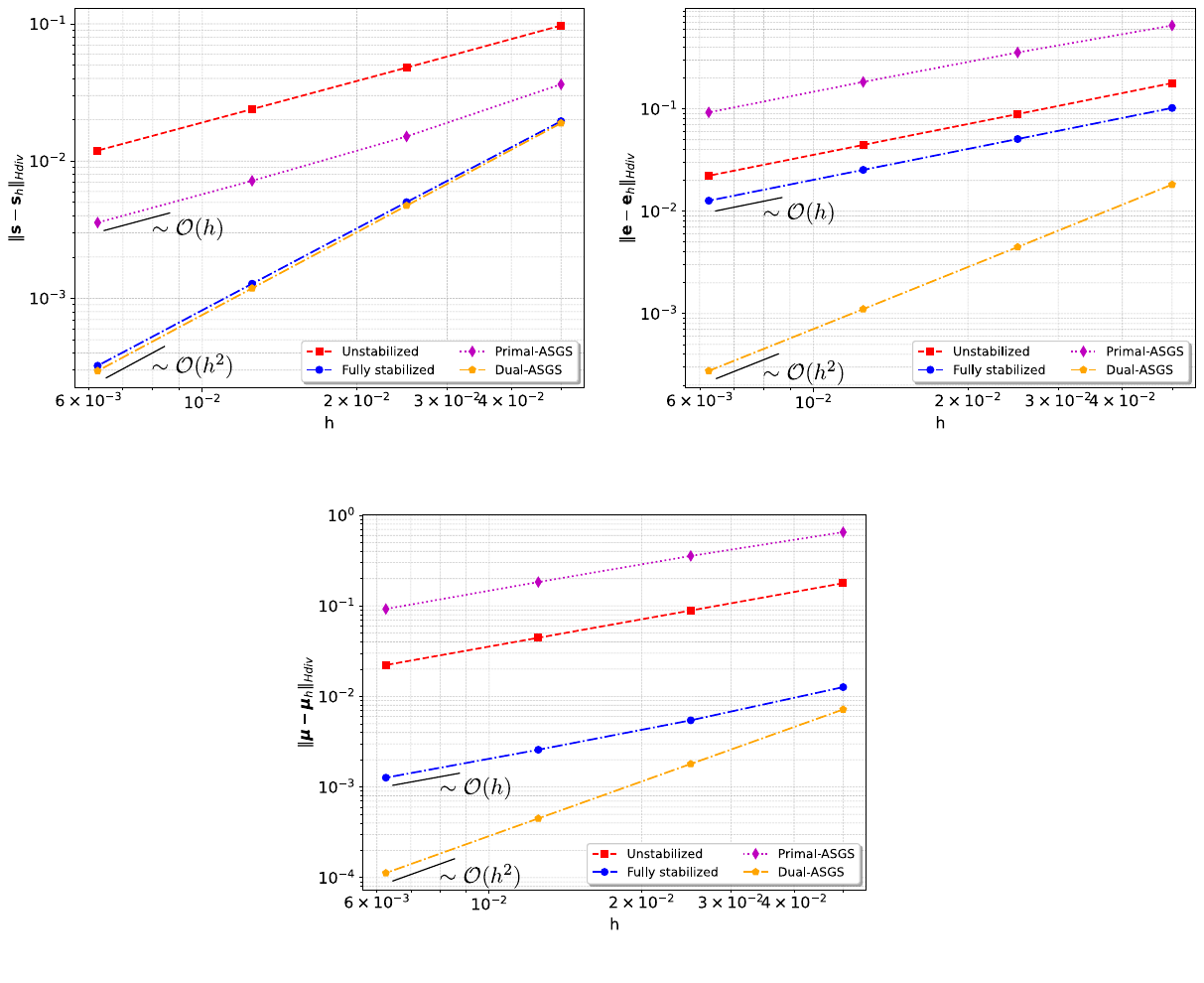}}%
  \end{picture}%
\endgroup%

	\label{fig:convergence_vector_quadratic_hdiv_art02}
\end{figure}
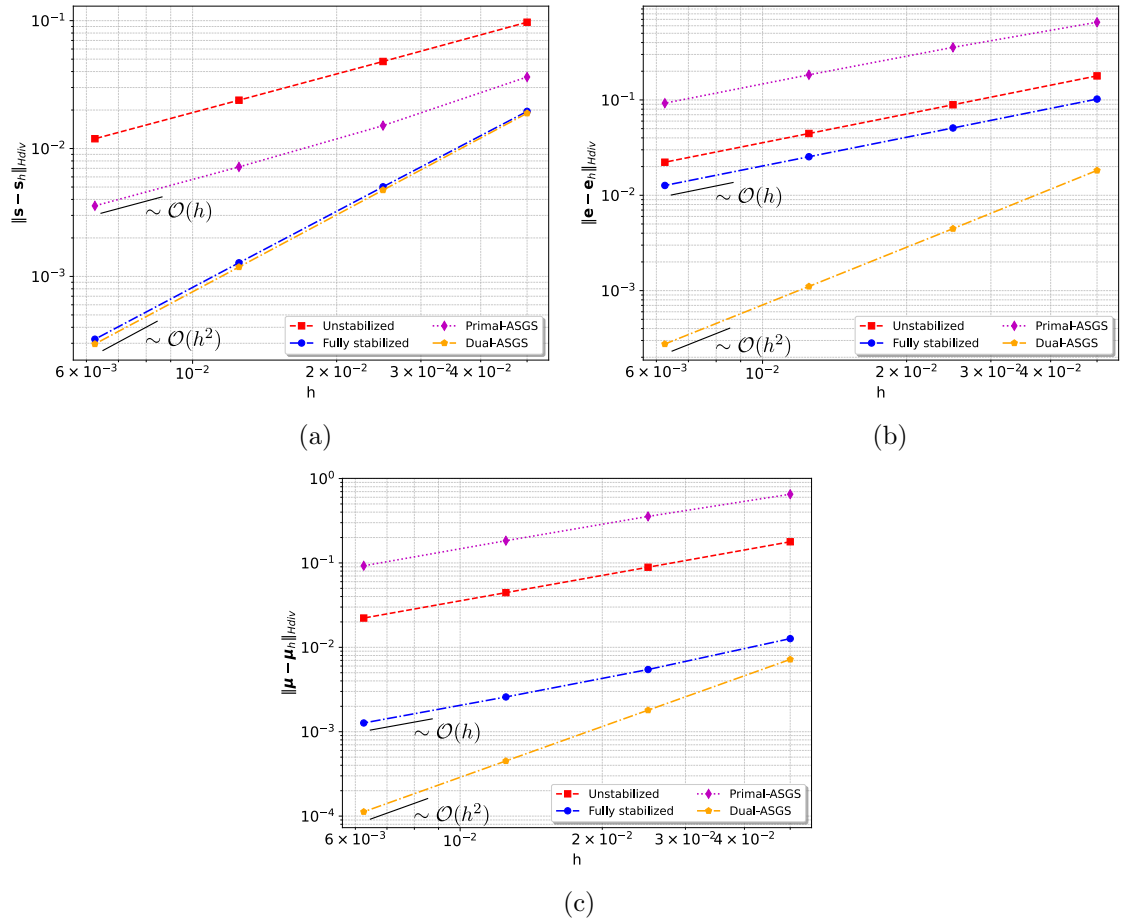

We now numerically assess the convergence properties of the discrete formulations. To this end, Figures~\ref{fig:convergence_scalar_linear_art02} to \ref{fig:convergence_vector_linear_hdiv_art02} present the error norms for all fields using linear elements ($k=1$), while Figures~\ref{fig:convergence_scalar_quadratic_art02} to \ref{fig:convergence_vector_quadratic_hdiv_art02} show the corresponding errors obtained with quadratic elements ($k=2$). Recall that we are also comparing our results against the different formulations recently evaluated in \cite{Bazon2025}: namely, a (stable) Galerkin-based \textit{Natural} formulation in which $\nabla{U}_h \subset M_h$, $\nabla{L}_h \subset S_h$, and $E_h = M_h$; a \textit{Fully stabilized} GLS formulation permitting equal-order interpolation; and finally, an \textit{Unstabilized} equal-order plain Galerkin formulation.

For linear elements, see Figures~\ref{fig:convergence_scalar_linear_art02} to \ref{fig:convergence_vector_linear_hdiv_art02}, all formulations achieve optimal convergence rates: second-order for the scalar fields in the $\Ltwo$-norm and first-order in the $\Hone$-norm. However, the errors for the \textit{Dual-ASGS} formulation are significantly larger, highlighting the superior accuracy of the primal formulation to approximate the scalar fields.
Note also that, in this case, only the $\Ltwo$-norm is assessed for the \textit{Dual-ASGS} formulation. For the vector fields, see Figure~\ref{fig:convergence_vector_linear_art02}, all formulations under consideration, except for the \textit{Natural} formulation, exhibit superlinear convergence in the $\Ltwo$-norm, nearly achieving $\mathcal{O}(h^2)$ despite a theoretical estimate of $\mathcal{O}(h)$.
Results for the $\Hdiv$-norm are shown in Figure~\ref{fig:convergence_vector_linear_hdiv_art02} for all formulations except the \textit{Natural} one. Recall that only the dual formulation guarantees convergence for the vector fields $\sF$ and $\mF$ in the $\Hdiv$-norm. Nevertheless, for linear elements, all cases exhibit $\mathcal{O}(h)$ convergence, with very similar errors for $\sF$. For $\mF$ (yellow lines), the dual formulation yields significantly better results and is the only one achieving linear convergence.

Moving to quadratic elements ($k=2$), the primal formulation achieved $\mathcal{O}(h^2)$ convergence in the $\Hone$-norm for the scalar fields $\uF$ and $\lF$, which is the theoretical optimal rate, see Figure~\ref{fig:convergence_scalar_quadratic_art02}. Its results are essentially identical to those of the \textit{Natural} and \textit{Fully stabilized} GLS formulations. The \textit{Unstabilized} formulation exhibited suboptimal convergence rates, underscoring the importance of stabilization, that was not evident for linear elements. For these fields in the $\Ltwo$-norm, the dual formulation shows the largest errors among all formulations since it only attains $\mathcal{O}(h^2)$ in contrast to
the primal and GLS formulations that achieved $\mathcal{O}(h^3)$.

For the vector fields $\sF$, $\mF$ and $\eF$, all formulations achieved second-order convergence in the $\Ltwo$-norm, see Figure~\ref{fig:convergence_vector_quadratic_art02}. Notably, the \textit{Fully stabilized} GLS formulation exhibited near-cubic ($\mathcal{O}(h^3)$) superconvergence for $\sF$. This behavior may, however, be specific to this particular example. As anticipated, the primary advantage of the dual formulation is evident in the $\Hdiv$-norm errors, see Figure~\ref{fig:convergence_vector_quadratic_hdiv_art02}, where it is the only method to achieve second-order convergence for $\mF$. Although theoretical estimates for $\eF$ in this norm are unavailable, we note that the dual formulation yields significantly lower errors for $\eF$ compared to all other methods.
These numerical results confirm the well posedness of the discrete formulations and validate our implementation of both the primal and dual methods introduced earlier.


\clearpage
\subsection{Example 2: Performance assessment in the presence of large gradients}

This final example aims to further assess the performance of the primal and dual formulations on a more challenging problem characterized by large internal gradients, typical of reaction-diffusion systems. We introduce a standard problem that serves a dual purpose: it will generate the synthetic data fields $\sFt$ and $\eFt$ required for the data-driven formulation while also providing a reference solution for comparative evaluation.

\begin{figure}[h]
	\centering
	\caption{Computational domain and the associated finite element mesh. Figures (b--d) illustrate the corresponding reference (exact) physical fields: the potential $\uF$, the flux $\sF$, and the gradient $\eF$, respectively.}
	~\\
	\fontsize{12pt}{10pt}\selectfont 
	\def\svgwidth{0.9\columnwidth}
	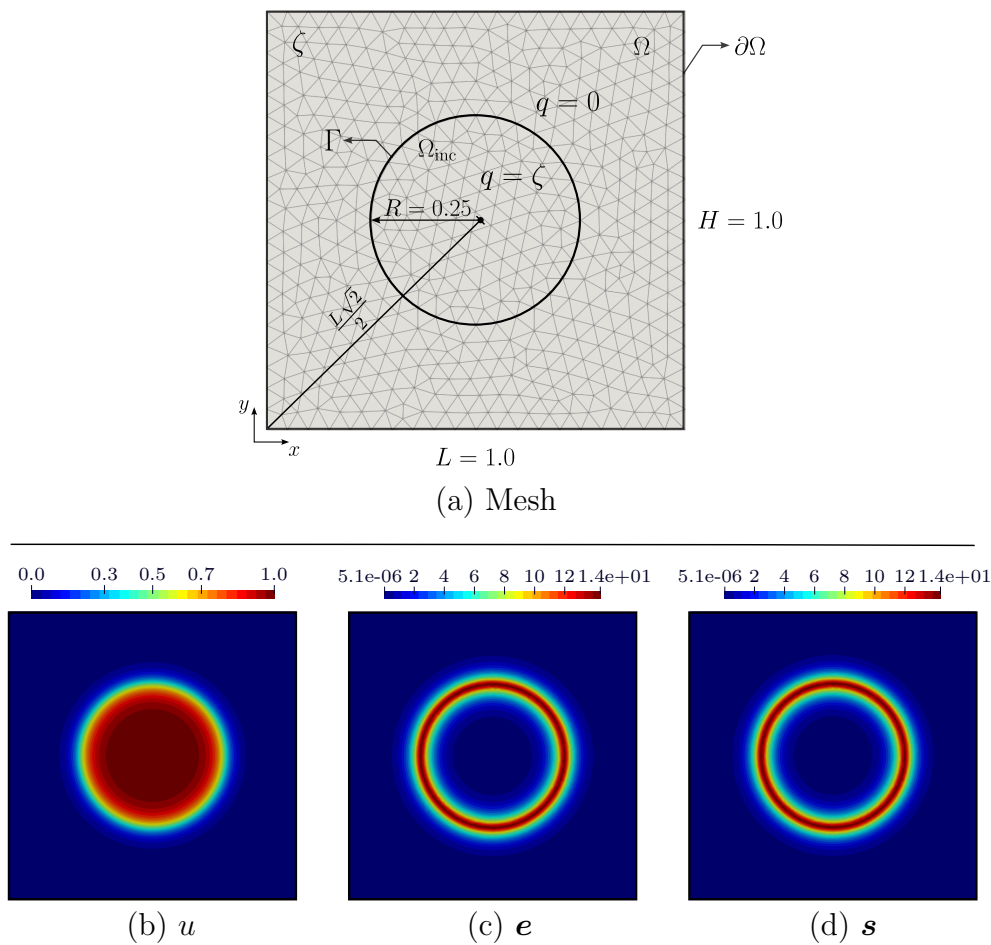
	\label{fig:inclusion_mesh_exact_fields}
\end{figure}

Let  $\Omega = (0,1)\times(0,1)$ be the computational domain.
We consider a circular inclusion immersed in $\Omega$, centered at ${\bm x}_{c} = (\frac12,\,\frac12)$ and defined by
$$
\Omega_{\mathrm{inc}} = \{\, {\bm x}\in \Omega : \|{\bm x}-{\bm x}_{c}\| < R \,\},
\qquad R=0.25,
$$
whose boundary is denoted by $\Gamma$. 
Figure~\ref{fig:inclusion_mesh_exact_fields}~(a) presents a schematic of the computational domain along with the finite element mesh. As shown, the interface $\Gamma$ is aligned with the element edges. 

The standard problem therefore consists of solving the following diffusion–reaction equation:
Find $\uF:\Omega \to \mathbb{R}$ such that
\begin{equation}
	\left\lbrace
	\begin{array}{rcll}
		-\nabla \cdot (\gamma \,\nabla \uF)+ \zeta\, \uF &=& q
		&\quad \text{in } \Omega, \\[4pt]
		\nabla \uF \cdot \normal &=& g 
		&\quad \text{on } \partial\Omega,
	\end{array}
	\right.
	\label{eq:diff_reaction}
\end{equation}
where the reaction coefficient $\zeta > 0$ is uniform and the source term $q$ 
is piecewise constant, i.e.,
$$
q({\bm x}) =
\begin{cases}
	\zeta, & {\bm x}\in \Omega_{\mathrm{inc}}, \\[3pt]
	0, & {\bm x}\in \Omega\setminus \Omega_{\mathrm{inc}}.
\end{cases}
$$
Assuming uniform conductivity $\gamma$, symmetry implies the solution depends only on the radial coordinate, i.e., $\uF = \uF(r)$ with $r = \sqrt{(x-\frac12)^{2} + (y-\frac12)^{2}}$. Under this assumption, the problem reduces to an ordinary differential equation in $r$, whose solution is
\begin{equation}
	u(r) =
	\begin{cases}
		\displaystyle
		1 - \alpha R\,K_{1}(\alpha R)\,I_{0}(\alpha r),
		& 0 \le r \le R, \\[10pt]
		\displaystyle
		\alpha R\,I_{1}(\alpha R)\,K_{0}(\alpha r),
		& r > R,
	\end{cases}
	\qquad
	\alpha = \sqrt{\frac{\zeta}{\gamma}},
\end{equation}
Here, $I_{0}$ and $I_{1}$ are the modified Bessel functions of the first kind (orders 0 and 1, respectively), and $K_{0}$ and $K_{1}$ are the corresponding modified Bessel functions of the second kind. Figure~\ref{fig:inclusion_mesh_exact_fields}~(b) shows the exact potential field $\uF$, while Figures~\ref{fig:inclusion_mesh_exact_fields}~(c--d) illustrate the corresponding gradient and flux fields, $\eF$ and $\sF$, respectively. Several observations are noteworthy: first, note that this solution satisfies the following physical conditions across the interface $\Gamma$:
$$
\text{at } r=R:\qquad
\begin{cases}
	u(R^{+}) = u(R^{-}), \\[3pt]
	\gamma\,u'(R^{+}) = \gamma\,u'(R^{-}),
\end{cases}
$$
ensuring continuity of $u$ and of the normal diffusive flux. 
Second, the solution decays to zero at infinity, and its boundedness at the origin implies $u'(0) = 0$. Third, the sharpness of the transition region is controlled by adjusting the ratio $\alpha$ of the reaction to the diffusion coefficient. In our numerical experiments, we set $\gamma = 1$ and $\zeta = 1000$, which yields a transition region with a characteristic length of approximately $0.03$. This behavior is clearly illustrated in Figure~\ref{fig:inclusion_mesh_exact_fields}~(b) Also, Figures~\ref{fig:inclusion_mesh_exact_fields}~(c--d) show the very sharp gradients in the vector fields $\eF$ and $\sF$.

Before presenting the numerical results, we specify the following computational setup:
\begin{enumerate}
\item[(i)] A finite element mesh with characteristic size $h = 0.025$;
\item[(ii)] Stabilization parameters set as
\[
\ku = \kl = 1, \qquad \ks = \ke = \km = \tfrac{1}{8}
\]
for the \textit{Primal-ASGS} formulation, and
\[
\ku = \kl = \ke = \tfrac{1}{4}, \qquad \ks = 1, \qquad \km = \tfrac{1}{8}
\]
for the \textit{Dual-ASGS} formulation;

\item[(iii)] Strong imposition of boundary conditions: for the primal formulation, $\uF = 0$ and $\lF = 0$; for the dual formulation, $\sF \cdot \normal = 0$ and $\mF \cdot \normal = 0$. These
homogenous values are safe since the solution has already sufficiently decayed at the boundary.

\item[(iv)] The data fields $\sFt$ and $\eFt$ are defined as elementwise constants. 
They are generated by sampling a finite element solution of problem (\ref{eq:diff_reaction}),
which is computed using a standard primal Galerkin formulation with continuous piecewise
linear elements on the same mesh used for the data-driven problem. Unlike the previous example, where data fields were treated as continuously available, here the data are provided only at discrete locations, that is, one data-pair per element. This brings the test closer to 
realistic data-driven mechanics applications.
\end{enumerate}

\begin{figure}[ht!!]
	\centering
	\caption{Figures (a--b) illustrate the field $\uF$ for the \textit{Primal}- and \textit{Dual-ASGS} formulations. Likewise, Figures (c--d) show the field $\sF$ for the \textit{Primal}- and \textit{Dual-ASGS} formulations.}
    ~\\
	\fontsize{12pt}{10pt}\selectfont 
	\def\svgwidth{0.95\columnwidth}
	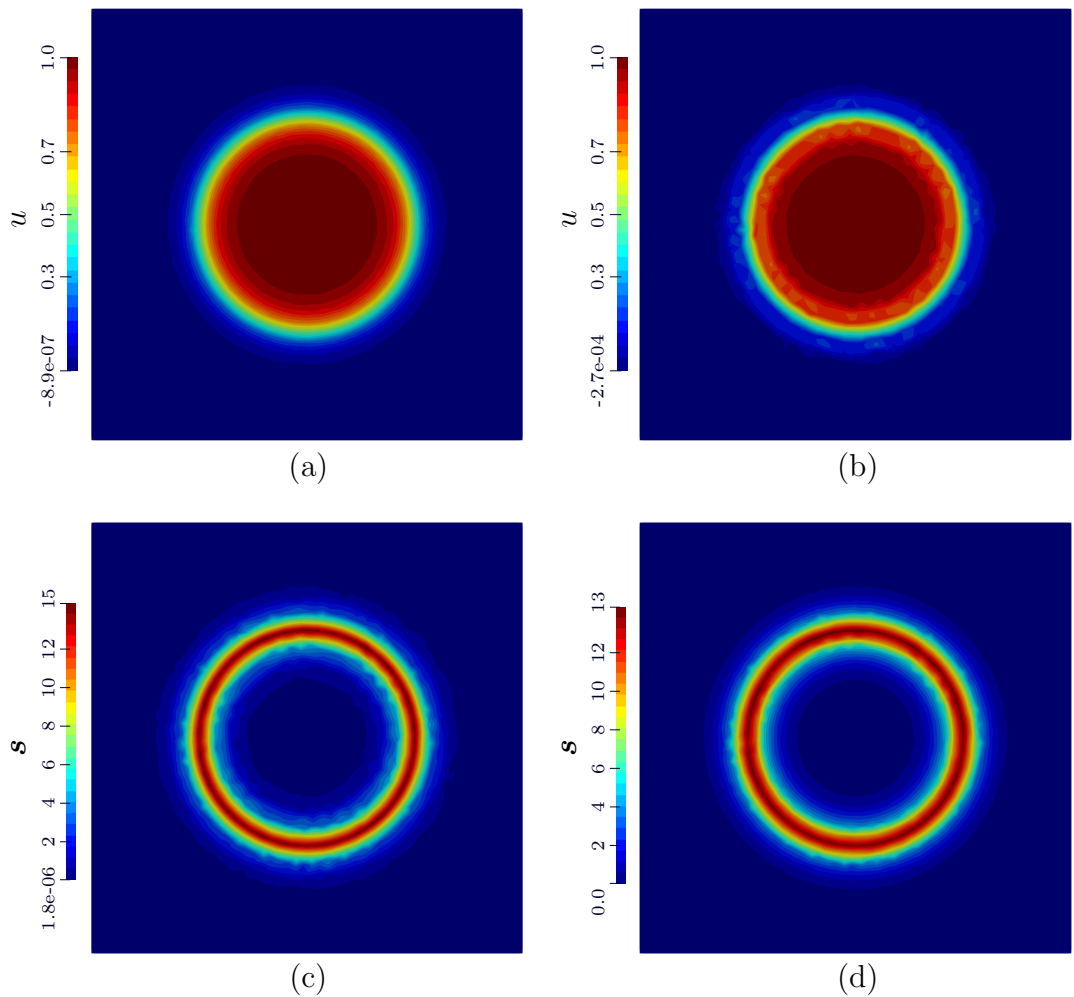
	\label{fig:inclusion_u_s_primal_dual_ASGS}
\end{figure}

\begin{figure}[ht!!]
	\centering
	\caption{Figure (a) displays the $\uF$ profile along the coordinate diagonal for the exact reference solution, together with the corresponding profiles obtained from the \textit{Primal}- and \textit{Dual-ASGS} formulations. Likewise, Figure (b) displays the $\sF$ profile along the coordinate diagonal for the exact reference solution and for the \textit{Primal}- and \textit{Dual-ASGS} formulations.}
	\fontsize{12pt}{10pt}\selectfont 
	\def\svgwidth{0.8\columnwidth}
\begingroup%
  \makeatletter%
  \providecommand\color[2][]{%
    \errmessage{(Inkscape) Color is used for the text in Inkscape, but the package 'color.sty' is not loaded}%
    \renewcommand\color[2][]{}%
  }%
  \providecommand\transparent[1]{%
    \errmessage{(Inkscape) Transparency is used (non-zero) for the text in Inkscape, but the package 'transparent.sty' is not loaded}%
    \renewcommand\transparent[1]{}%
  }%
  \providecommand\rotatebox[2]{#2}%
  \newcommand*\fsize{\dimexpr\f@size pt\relax}%
  \newcommand*\lineheight[1]{\fontsize{\fsize}{#1\fsize}\selectfont}%
  \ifx\svgwidth\undefined%
    \setlength{\unitlength}{425.19685039bp}%
    \ifx\svgscale\undefined%
      \relax%
    \else%
      \setlength{\unitlength}{\unitlength * \real{\svgscale}}%
    \fi%
  \else%
    \setlength{\unitlength}{\svgwidth}%
  \fi%
  \global\let\svgwidth\undefined%
  \global\let\svgscale\undefined%
  \makeatother%
  \begin{picture}(1,1.44437683)%
    \lineheight{1}%
    \setlength\tabcolsep{0pt}%
    \put(0.48985262,0.73495356){\color[rgb]{0,0,0}\makebox(0,0)[lt]{\lineheight{1.25}\smash{\begin{tabular}[t]{l}(a)\end{tabular}}}}%
    \put(0.48946998,0.00459162){\color[rgb]{0,0,0}\makebox(0,0)[lt]{\lineheight{1.25}\smash{\begin{tabular}[t]{l}(b)\end{tabular}}}}%
    \put(0,0){\includegraphics[width=\unitlength,page=1]{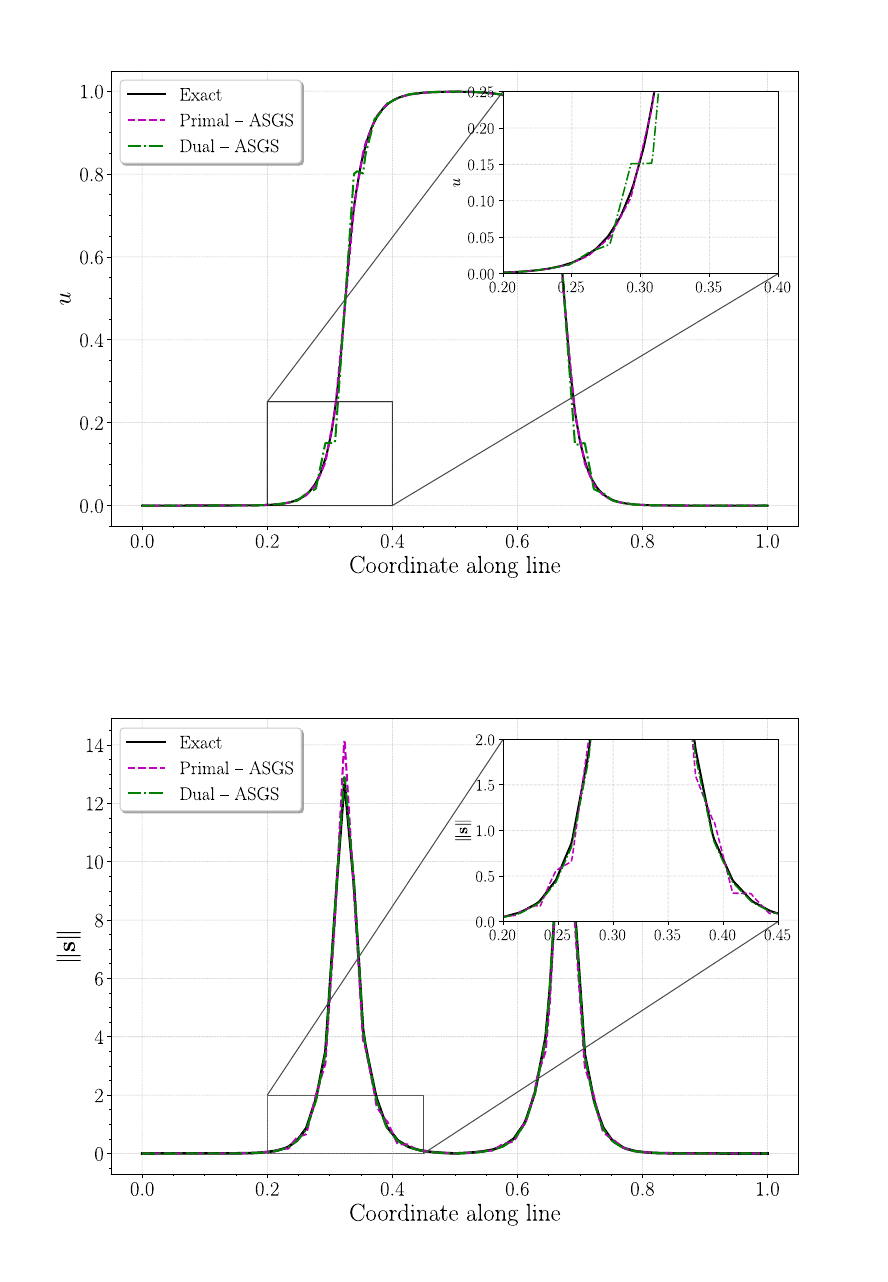}}%
  \end{picture}%
\endgroup%

	\label{fig:plot_over_line_u_s_primal_dual_ASGS}
\end{figure}

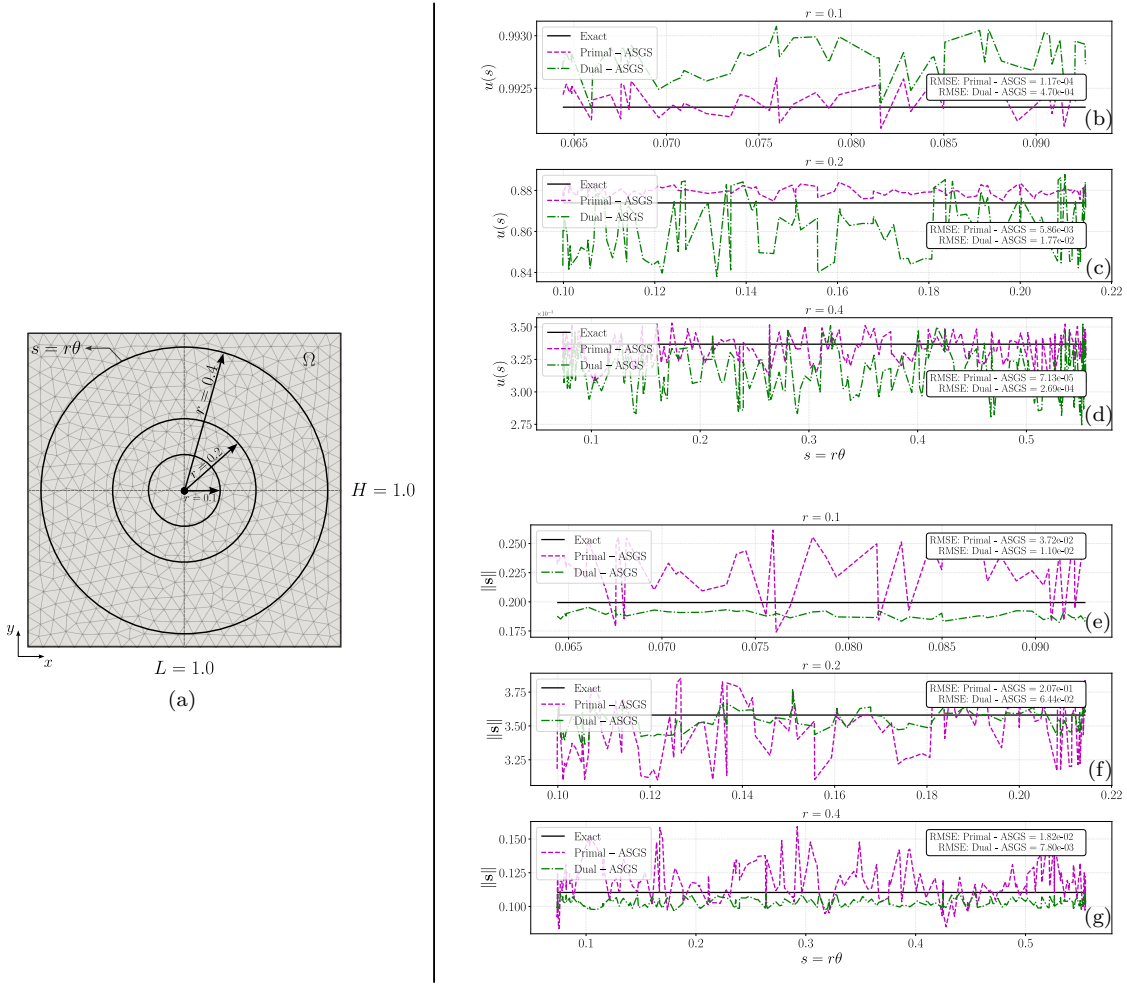
\begin{figure}[ht!!]
	\centering
	\caption{Figures (b--d) present the $\uF$ profiles along three circumferences of radius $r = 0.1$, $r = 0.2$ and $r = 0.4$, illustrated schematically in Figure (a). Each plot compares the exact reference solution with the corresponding profiles obtained from the \textit{Primal}- and \textit{Dual-ASGS} formulations. Likewise, Figures (e--g) display the $\|\sF\|$ profiles along the same radii for both the exact reference solution and the \textit{Primal}- and \textit{Dual-ASGS} formulations.}
    ~\\
	\fontsize{8pt}{10pt}\selectfont 
	\def\svgwidth{1.0\columnwidth}
\begingroup%
  \makeatletter%
  \providecommand\color[2][]{%
    \errmessage{(Inkscape) Color is used for the text in Inkscape, but the package 'color.sty' is not loaded}%
    \renewcommand\color[2][]{}%
  }%
  \providecommand\transparent[1]{%
    \errmessage{(Inkscape) Transparency is used (non-zero) for the text in Inkscape, but the package 'transparent.sty' is not loaded}%
    \renewcommand\transparent[1]{}%
  }%
  \providecommand\rotatebox[2]{#2}%
  \newcommand*\fsize{\dimexpr\f@size pt\relax}%
  \newcommand*\lineheight[1]{\fontsize{\fsize}{#1\fsize}\selectfont}%
  \ifx\svgwidth\undefined%
    \setlength{\unitlength}{745.70238489bp}%
    \ifx\svgscale\undefined%
      \relax%
    \else%
      \setlength{\unitlength}{\unitlength * \real{\svgscale}}%
    \fi%
  \else%
    \setlength{\unitlength}{\svgwidth}%
  \fi%
  \global\let\svgwidth\undefined%
  \global\let\svgscale\undefined%
  \makeatother%
  \begin{picture}(1,0.87464599)%
    \lineheight{1}%
    \setlength\tabcolsep{0pt}%
    \put(0.14344082,0.24695774){\color[rgb]{0,0,0}\makebox(0,0)[lt]{\lineheight{1.25}\smash{\begin{tabular}[t]{l}(a)\end{tabular}}}}%
    \put(0,0){\includegraphics[width=\unitlength,page=1]{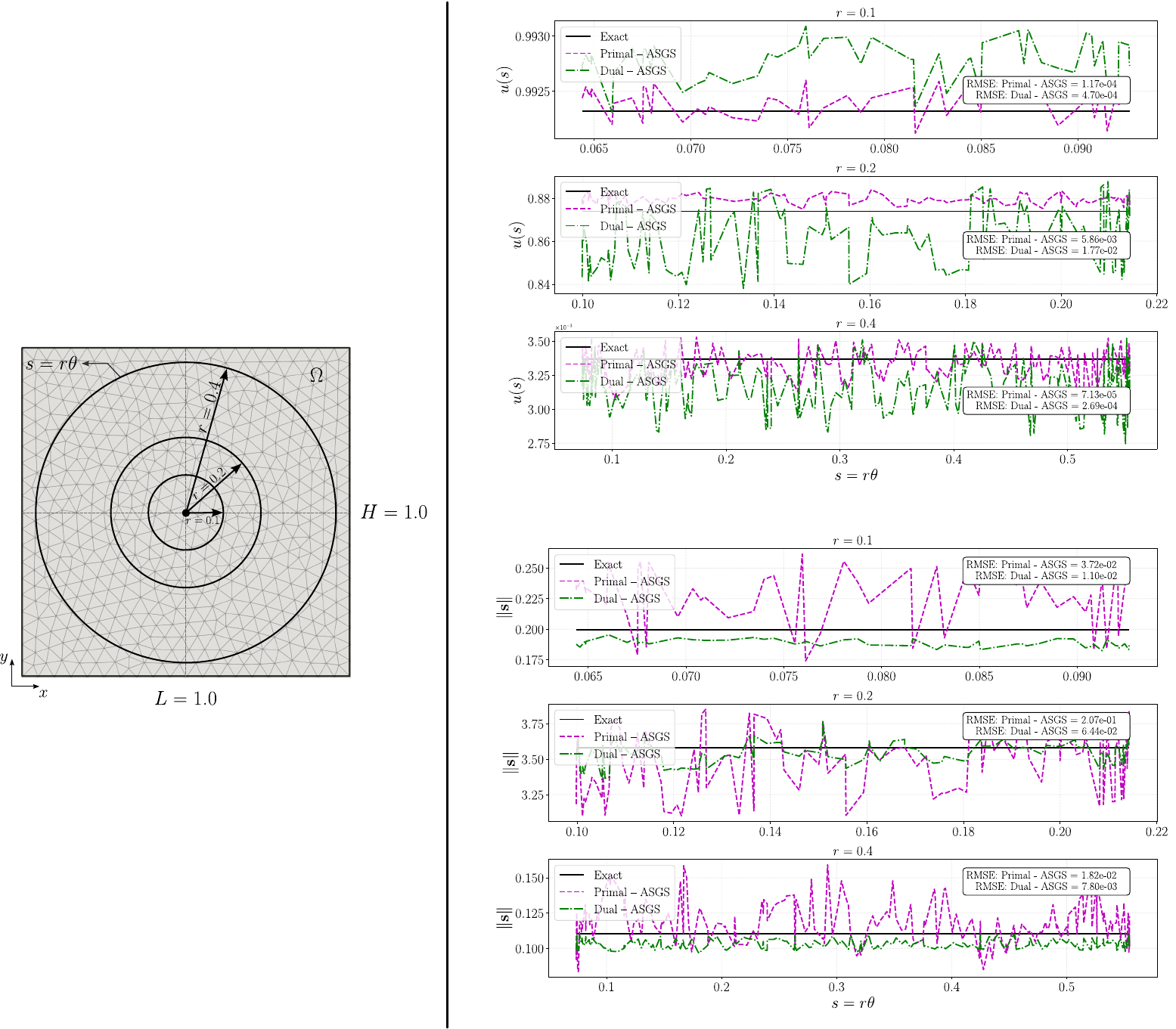}}%
    \put(0.95960293,0.76455834){\color[rgb]{0,0,0}\makebox(0,0)[lt]{\lineheight{1.25}\smash{\begin{tabular}[t]{l}(b)\end{tabular}}}}%
    \put(0.96052101,0.6324506){\color[rgb]{0,0,0}\makebox(0,0)[lt]{\lineheight{1.25}\smash{\begin{tabular}[t]{l}(c)\end{tabular}}}}%
    \put(0.95960293,0.50134837){\color[rgb]{0,0,0}\makebox(0,0)[lt]{\lineheight{1.25}\smash{\begin{tabular}[t]{l}(d)\end{tabular}}}}%
    \put(0.95981395,0.31659254){\color[rgb]{0,0,0}\makebox(0,0)[lt]{\lineheight{1.25}\smash{\begin{tabular}[t]{l}(e)\end{tabular}}}}%
    \put(0.96265785,0.18448473){\color[rgb]{0,0,0}\makebox(0,0)[lt]{\lineheight{1.25}\smash{\begin{tabular}[t]{l}(f)\end{tabular}}}}%
    \put(0.95960293,0.05338265){\color[rgb]{0,0,0}\makebox(0,0)[lt]{\lineheight{1.25}\smash{\begin{tabular}[t]{l}(g)\end{tabular}}}}%
  \end{picture}%
\endgroup%

	\label{fig:plot_over_curve_u_s_primal_dual_ASGS}
\end{figure}

\begin{figure}[ht!!]
	\centering
	\caption{Figures (a--b) present the error field $\uF - \uFh$ for the \textit{Primal}- and \textit{Dual-ASGS} formulations. In the same way, Figures (c--d) present the error field $\sF - \sFh$ for the \textit{Primal}- and \textit{Dual-ASGS} formulations.}
    ~\\
	\fontsize{12pt}{10pt}\selectfont 
	\def\svgwidth{1.0\columnwidth}
	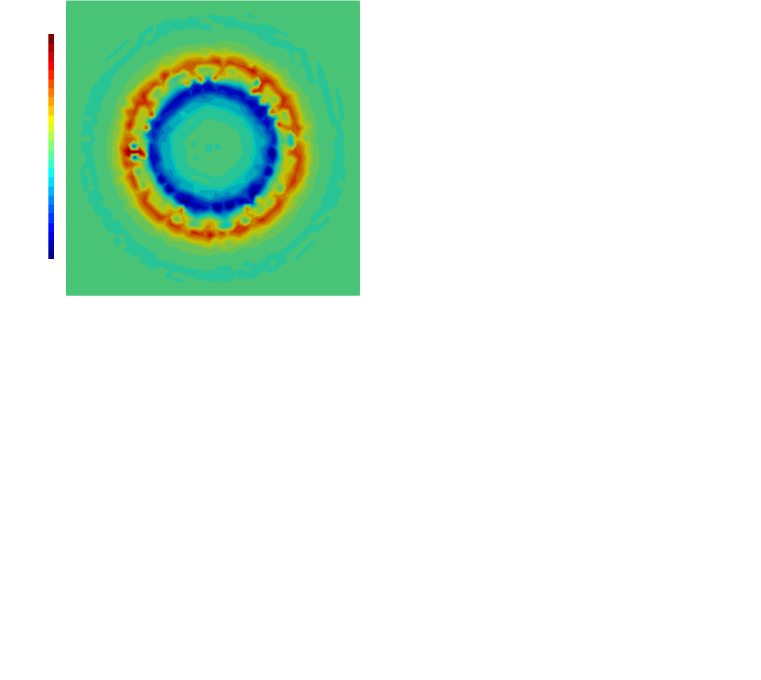
	\label{fig:inclusion_udif_sdif_primal_dual_ASGS}
\end{figure}

Turning to the numerical results, Figure~\ref{fig:inclusion_u_s_primal_dual_ASGS} illustrates
the physical fields $\uF$ and $\sF$ obtained using both the \textit{Primal}- and \textit{Dual-ASGS} formulations.
Figures~\ref{fig:inclusion_u_s_primal_dual_ASGS}~(a--b) demonstrate that the primal
formulation provides a superior approximation for the scalar potential field, whereas the dual 
formulation shows considerable distortion compared to the reference solution plotted in Figure ~\ref{fig:inclusion_mesh_exact_fields}. On the other hand, the dual formulation approximates the 
flux field more effectively. This is evident in Figures~\ref{fig:inclusion_u_s_primal_dual_ASGS}~(c--d), where the largest distortions  are observed in the primal result. Unlike the previous test, where improvements were confined to the $\Hdiv$-norm of the error, here the advantage of the dual formulation is already visible by 
simple inspection. This trend is further evident in Figures~\ref{fig:plot_over_line_u_s_primal_dual_ASGS} and~\ref{fig:plot_over_curve_u_s_primal_dual_ASGS}, 
which show profiles of $\uF$ and $\sF$ along the domain diagonal and along circumferential curves 
at different radii indicated in Figure~\ref{fig:plot_over_curve_u_s_primal_dual_ASGS}~(a), namely, $r=0.1,~0.2$ 
and $0.4$. Figure~\ref{fig:plot_over_line_u_s_primal_dual_ASGS}~(a) and the top 
part of Figure~\ref{fig:plot_over_curve_u_s_primal_dual_ASGS} reveal
that the primal formulation (magenta lines) clearly outperforms the
dual formulation (green lines) in capturing $\uF$. Conversely, Figure~\ref{fig:plot_over_line_u_s_primal_dual_ASGS}~(b) and the bottom part
of Figure~\ref{fig:plot_over_curve_u_s_primal_dual_ASGS} reveal the opposite situation for $\sF$. 
A closer view of these differences is provided in the insets of Figure~\ref{fig:plot_over_line_u_s_primal_dual_ASGS}.
See also the root mean square error reported for the circumferential profiles for each formulation
in Figure~\ref{fig:plot_over_curve_u_s_primal_dual_ASGS}. 
For further clarity, it is instructive to examine the pointwise error distribution across the
computational domain. Figure~\ref{fig:inclusion_udif_sdif_primal_dual_ASGS} presents these error
fields: parts (a--b) show the error in $\uF - \uF_h$ while parts (c--d) show the
error in the flux $\sF - \sF_h$. Results for the \textit{Primal-ASGS} formulation appear in the left column; those for the
\textit{Dual-ASGS} formulation appear on the right.

\subsubsection*{Numerical results in the presence of noisy data}

\begin{figure}[ht!!]
	\centering
	\caption{Data field $\sFt$ sampled from a finite element solution of problem
		(\ref{eq:diff_reaction}) to which uniform noise has been added. Different noise levels are considered: (a) $\delta=0.25$, (b) $\delta=0.5$ and (c) $\delta=1$.}
	~\\
	\fontsize{12pt}{10pt}\selectfont 
	\def\svgwidth{1.0\columnwidth}
	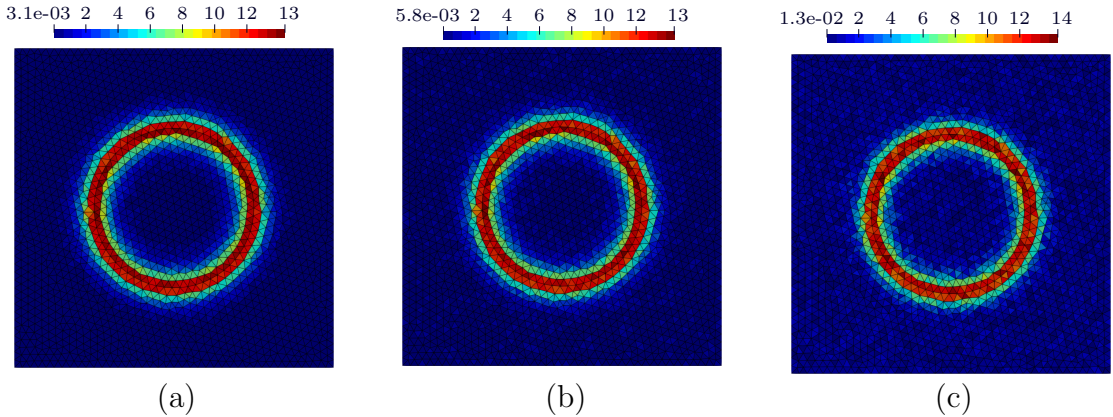
	\label{fig:st_with_noise}
\end{figure}

\begin{figure}[h!]
	\centering
	\caption{Figures (a--c) present the noise-free data pairs $(\eFt,\sFt)$ together with the corresponding noisy pairs for three different values of $\delta$, namely $\delta = 0.25$, $\delta = 0.5$, and $\delta = 1.0$, respectively. Figures (d--f) show the corresponding approximations of the scalar field $\uF$ obtained with the \textit{Primal-ASGS} formulation. Similarly, Figures (g--i) display the approximations of the same field computed using the \textit{Dual-ASGS} formulation.}
	~\\
	\def\svgwidth{1.0\columnwidth}
	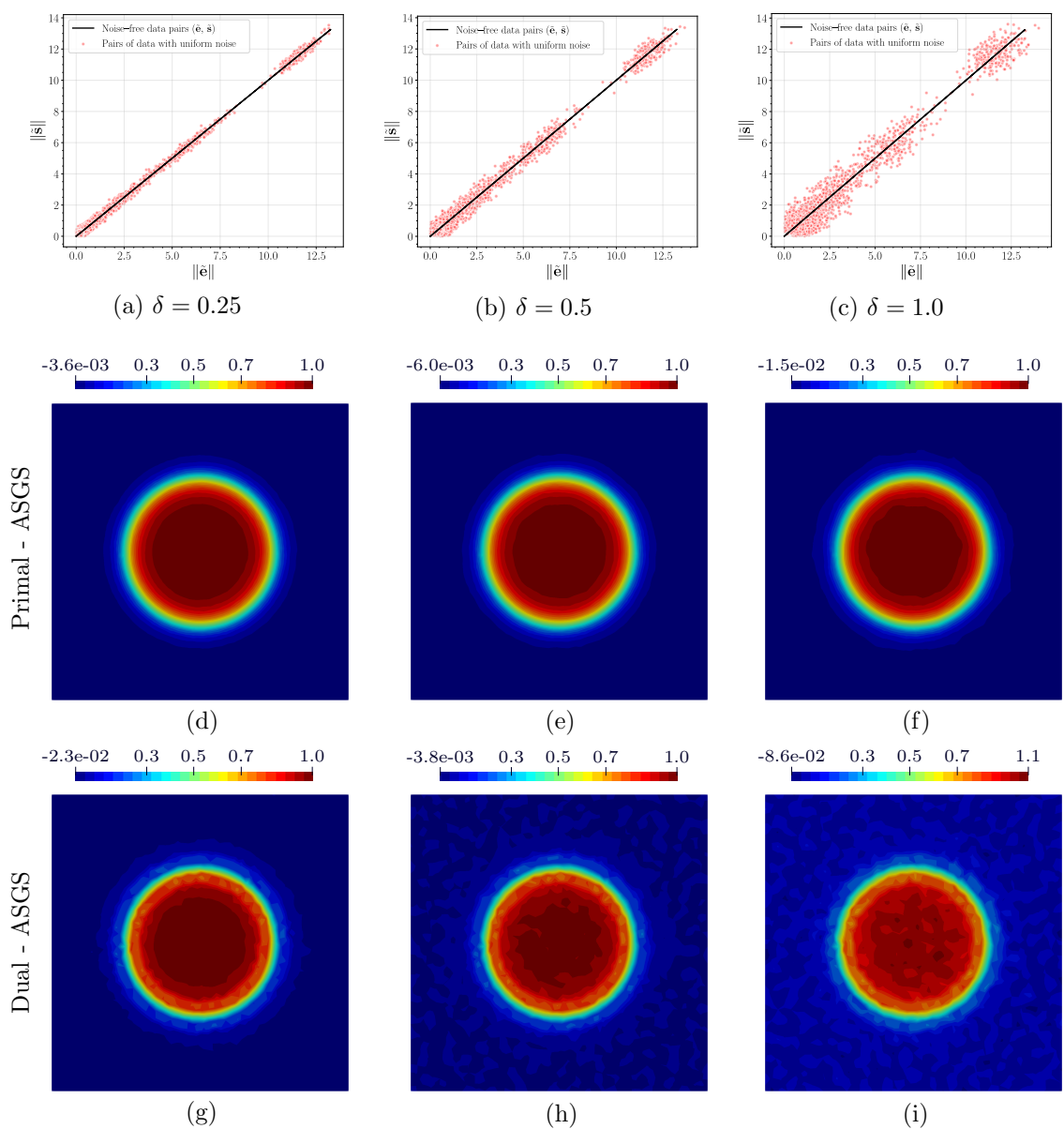
	\label{fig:noise_graph_and_u}
\end{figure}

\begin{figure}[h!]
	\centering
	\caption{Approximations of the vector field $\sF$ obtained with the \textit{Primal-ASGS} formulation (parts (a)--(c)) and with the \textit{Dual-ASGS} formulation (parts (d)--(f)): $\delta=0.25$ (left), $\delta=0.5$ (middle), $\delta=1$ (right)}
	~\\
	\def\svgwidth{1.0\columnwidth}
	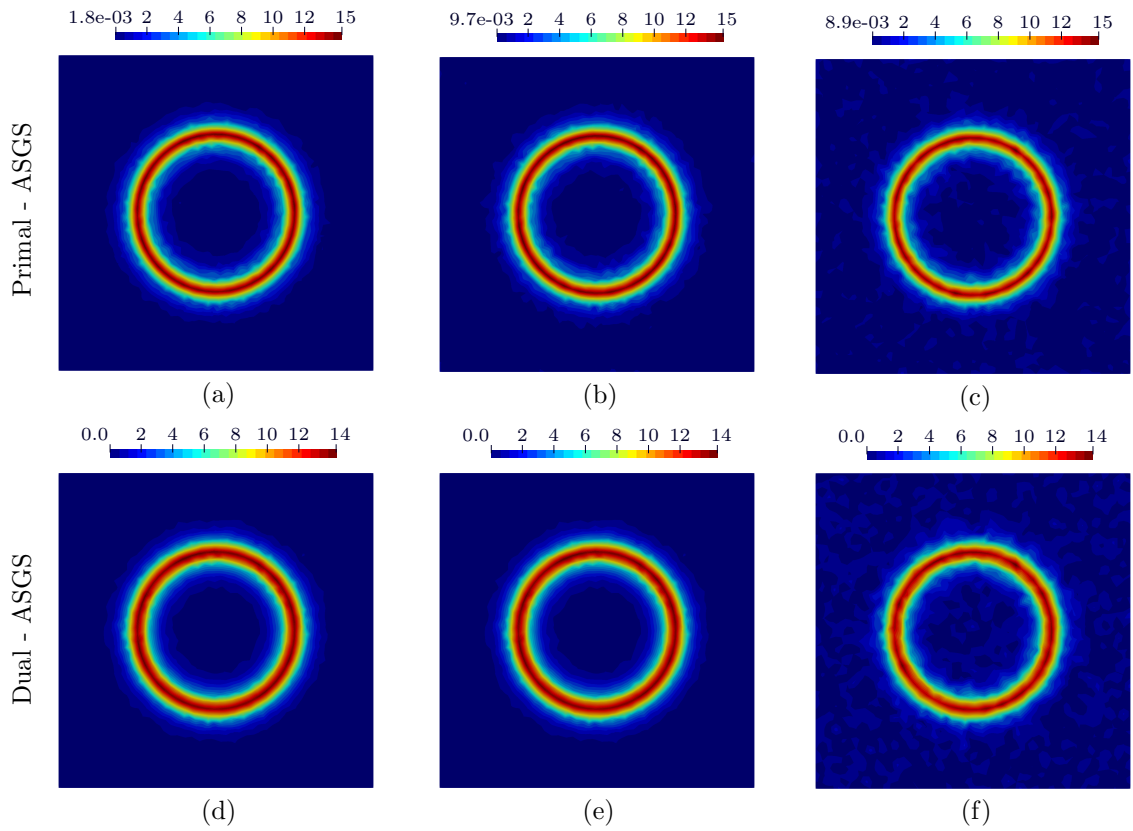
	\label{fig:noise_s}
\end{figure}

This experiment is well suited for examining the methods' behavior under conditions more representative of data-driven applications. To this end, we introduce an additional ingredient: incorporating noise into the data fields $\sFt$ and $\eFt$, which were originally sampled from a finite element solution of problem (\ref{eq:diff_reaction}). Specifically, we perturb these fields by adding uniform noise:
\begin{equation}
	\eFt^{\,\delta} = \eFt + \delta\,\boldsymbol{\xi}_{\eF},
	\qquad
	\sFt^{\,\delta} = \sFt + \delta\,\boldsymbol{\xi}_{\sF},
\end{equation}
where $\delta > 0$ represents the noise level and $\boldsymbol{\xi}_{\eF}$ and $\boldsymbol{\xi}_{\sF}$ are random vector fields with components independently drawn from a uniform distribution:
\begin{equation}
	\boldsymbol{\xi}_{\eF}, \boldsymbol{\xi}_{\sF} \sim \mathcal{U}(-1,1).
\end{equation}
This corresponds to adding a uniform perturbation in the interval $[-\delta, \delta]$ to each component of the gradient and flux fields. Figure~\ref{fig:st_with_noise} show the data field $\sFt^{\,\delta}$ for different noise levels, namely, (a) $\delta=0.25$, (b) $\delta=0.5$ and
(c) $\delta=1$. Figure~\ref{fig:noise_graph_and_u}~(a--c) further illustrates
both, the noise-free data pairs $(\eFt,\sFt)$, which form a straight line, since the flux and the gradient are related by a linear constitutive law, and the corresponding perturbed data pairs sampled from the finite element mesh.
The effect of increasing noise on the numerical approximation of $\uF$ is shown in Figure~\ref{fig:noise_graph_and_u}~(d--i) for both the \textit{Primal}- and \textit{Dual-ASGS} formulations. Corresponding results for $\sF$ appear in Figure~\ref{fig:noise_s}. These results reinforce the trends observed earlier: the primal formulation remains more robust for approximating the scalar field, while the dual formulation continues to better capture the vector field, even at high noise levels as seen in Figure~\ref{fig:noise_graph_and_u}~(f) and Figure~\ref{fig:noise_s}~(f). 

\begin{figure}[h!]
	\centering
	\caption{Contours of $\sF \cdot \eF$ to illustrate the violation of thermodynamic consistency across the domain for different levels of noise in the data fields $\eFt$ and $\sFt$: $\delta=0.25$ (left), $\delta=0.5$ (middle), $\delta=1$ (right).}
	~\\
	\def\svgwidth{1.0\columnwidth}
	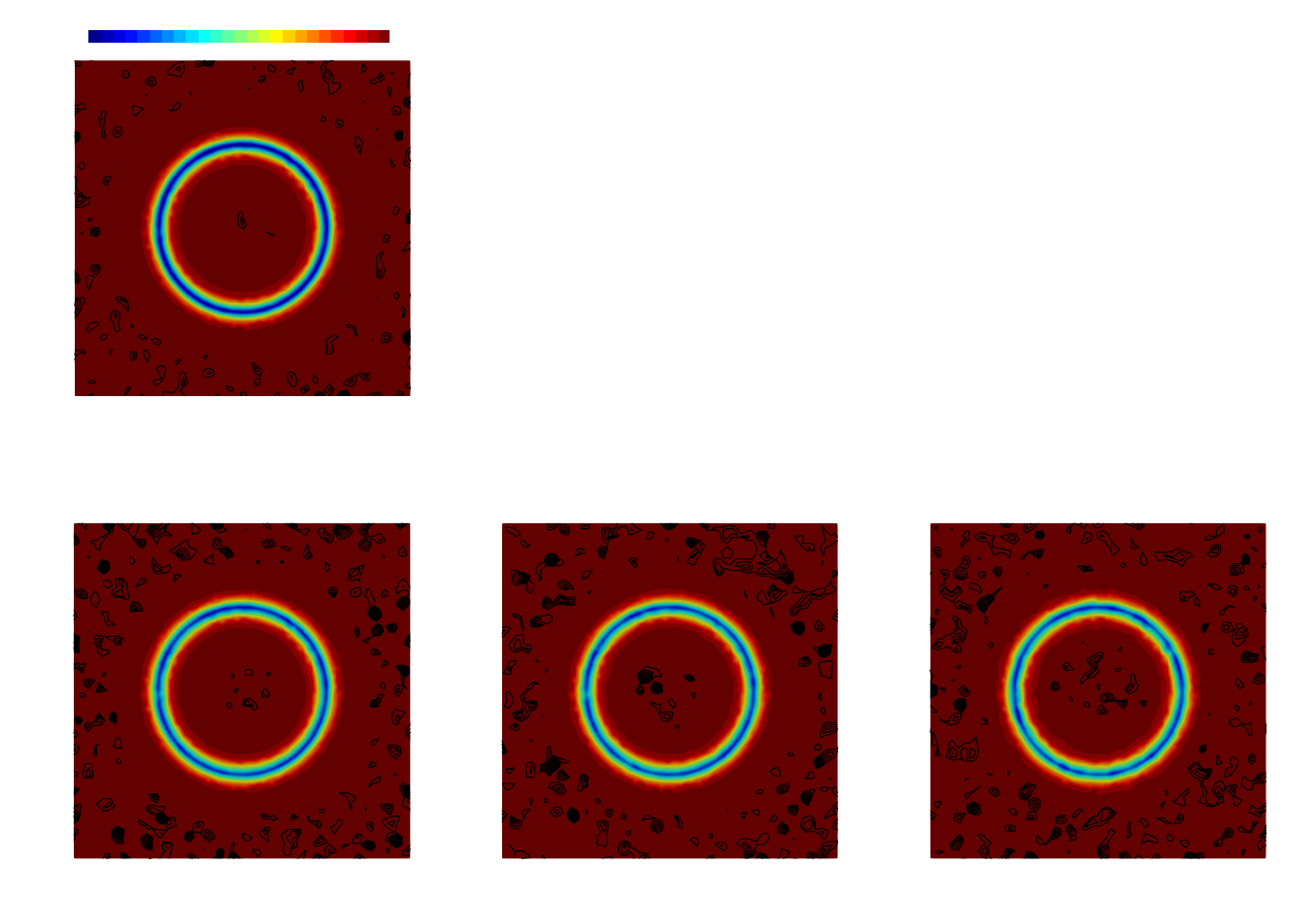
	\label{fig:noise_s_dot_e}
\end{figure}

The presence of noise introduces an additional consideration: the loss of thermodynamic consistency between $\eF$ and $\sF$, an effect present in both formulations. This is visualized 
in Figure~\ref{fig:noise_s_dot_e} via contours of the scalar field $\sF \cdot \eF$.
The onset of violation is indicated by the zero contour, shown as a black line.
The data near the origin exhibit a larger relative noise, leading to a higher concentration
of elements violating the thermodynamic consistency in regions where the solution gradient is 
low, that is, far from the interface $\Gamma$. Nothing in the formulation itself precludes 
these violations, but their occurrence underscores a central issue for reliable data-driven
computations in practical applications. The development of numerical methods capable of filtering potential experimental errors from the data fields and ensuring thermodynamic consistency is essential but left for future work.

\section{Conclusion}

In this work, we have proposed a finite element approximation of a PDE-constrained optimization problem arising in reaction-diffusion problems within the framework of Data-Driven Computational Mechanics (DDCM) based on the Variational Multiscale (VMS) concept.

We have established that both the primal and dual forms of the optimality conditions are well-posed in the continuous setting. This has been achieved by proving that the necessary requirements for mixed saddle-point problems are satisfied.

Regarding the finite element approximation, to overcome the limitations of standard Galerkin methods---which would require complex, problem-specific inf-sup conditions for the finite element spaces---we adopted the VMS framework. This approach allowed us to derive stable and consistent finite element approximations by systematically modeling the sub-grid scales. Specifically, we established the well-posedness of both the Algebraic Sub-Grid Scale (ASGS) and the Orthogonal Sub-Grid Scale (OSGS) methods.

For the VMS approximation, we have demonstrated that one can transition between the primal and the dual form of the problem simply by adjusting the design of the stabilization parameters, providing a unified view of the problem's discrete structure.

The numerical results presented have provided a robust qualitative and comparative assessment, confirming the stability and consistency of our stabilized finite element schemes. These results underscore the effectiveness of our approach in handling the inherent challenges of data-driven models.

While this work represents a major improvement in applying DDCM to diffusion-reaction problems, we acknowledge certain simplifications. In this study, we linearized the problem by omitting the explicit enforcement of the Second Law of Thermodynamics, though we verified its fulfillment a posteriori. We believe that the inclusion of this law directly into future optimization schemes will lead to even more robust formulations.

At the numerical level, we have focused on continuous finite element functions and quasi-uniform partitions; the extension of this VMS approach to discontinuous spaces and arbitrary mesh refinements remains a promising avenue for future research. Discontinuous (but conforming) interpolations could be considered by introducing the concept of sub-grid scales on the element boundaries \cite{Codina2009}, while non-uniform finite refinements could be addressed by using the techniques in \cite{codina2008}. Likewise, $hp$-convergence in the norm of the continuous problem and improved convergence estimates in $L^2(\Omega)$ represent important topics for subsequent research.

\cleardoublepage

\section*{CRediT authorship contribution statement}
\textbf{Ramon Codina:} Writing – original draft, Validation, Supervision, Methodology, Investigation, Formal analysis, Conceptualization;

\noindent \textbf{Roberto A. Ausas:} Writing – original draft, Software, Validation, Supervision, Methodology, Investigation, Formal analysis, Conceptualization;

\noindent \textbf{Pedro B. Bazon:} Writing – original draft, Software, Validation, Visualization, Investigation;

\noindent \textbf{Cristian G. Gebhardt:} Writing – original draft, Validation, Supervision, Methodology, Investigation, Formal analysis, Conceptualization;

\section*{Acknowledgments}
RC gratefully acknowledges the support received from the ICREA Acad\`emia Program, from the Catalan Government.

RFA thankfully acknowledges financial support from the São Paulo Research Foundation (FAPESP) under grants 2013/07375-0 and 2025/00975-9, and from the National Council for Scientific and Technological Development - CNPq under grant 308704/2022-3.

PBB thankfully acknowledges the support provided by the Coordenação de Aperfeiçoamento de Pessoal de Nível Superior – Brazil (CAPES), Finance Code 88887.975792/2024-00.

CGG gratefully acknowledges the financial support from the European Research Council through the ERC Consolidator Grant ``DATA-DRIVEN OFFSHORE'' (Project ID 101083157).

\bibliographystyle{plain}

\bibliography{bibliography}

\end{document}